\documentclass{amsart}
\usepackage{amssymb}

\newtheorem{thm}{Theorem}[section]
\newtheorem{prop}[thm]{Proposition}
\newtheorem{cor}[thm]{Corollary}
\newtheorem{lem}[thm]{Lemma}
\newtheorem{rema}[thm]{Remark}
\newtheorem{defn}[thm]{Definition}
\newtheorem{ex}[thm]{Example}

\numberwithin{equation}{section}

\newcommand{\dtheta}[0]{\frac{\partial}{\partial
\theta}}
\newcommand{\dz}[0]{\frac{\partial}{\partial z}}
\newcommand{\Z}{\mathbb{Z}_+}
 
\newcommand{\twoLow}[0]{2w \frac{\partial}{\partial w} + \rho
\frac{\partial}{\partial \rho}} 
\newcommand{\Lw}[0]{w^{j + 1} \frac{\partial}{\partial w} + (\frac{j +
1}{2})\rho w^j \frac{\partial}{\partial \rho}} 
\newcommand{\Gw}[0]{w^j \left( \frac{\partial}{\partial \rho} -
\rho \frac{\partial}{\partial w}\right)}

\newcommand{\asqrt}[0]{a_{\scriptscriptstyle{\square}}}

\newcommand{\Lx}[0]{x^{j + 1}
\frac{\partial}{\partial x} + (\frac{j +
1}{2})\varphi x^j \frac{\partial}{\partial
\varphi}}
\newcommand{\Gx}[0]{x^j \left( \frac{\partial}{\partial \varphi} -
\varphi \frac{\partial}{\partial x}\right)}

\begin{document}

\title[Change of variables for $N=1$ NS-VOSAs]{Superconformal change of
variables for $N=1$ Neveu-Schwarz vertex operator superalgebras}

\author{Katrina Barron}
\address{Department of Mathematics, University of Notre Dame,
Notre Dame, IN 46556}
\email{kbarron@nd.edu}
\thanks{The author was supported by an NSF Mathematical Sciences Postdoctoral
Research Fellowship.}


\date{July 2, 2003.}


\begin{abstract} Superconformal change of variables formulas for $N=1$
Neveu-Schwarz vertex operator superalgebras are presented for general
invertible superconformal changes of variables.   Using the underlying
worldsheet supergeometry of propagating superstrings, geometric proofs of
the change of variables formulas are given for the case of convergent
superconformal changes of variables.  More general formal algebraic proofs
of the change of variables formulas in the case of formal superconformal
changes of variables are then given.  Finally, isomorphic families of $N=1$
NS-VOSAs are derived {}from the superconformal change of variables formulas.
\end{abstract}

\maketitle

\section{Introduction} 

In this paper, we prove superconformal change of variables formulas for
an $N=1$ Neveu-Schwarz vertex operator superalgebras ($N=1$ NS -VOSAs)
and for general invertible superconformal changes of variables.  Using the
correspondence between the worldsheet supergeometry of propagating
superstrings in $N=1$ superconformal field theory and the algebraic
structure of $N=1$ NS-VOSA developed in \cite{B thesis}, \cite{B memoirs}
and \cite{B iso-thm}, we give geometric proofs of the change of variables
formulas for the case of convergent superconformal changes of variables. 
More general formal algebraic proofs of the change of variables formulas in
the case of formal superconformal changes of variables are then given.  In
addition, we introduce isomorphic families of $N=1$ NS-VOSAs using the
superconformal change of variables formulas. 

This paper is organized as follows.  In Section \ref{preliminaries}, we
give preliminary definitions including the notions of $N=1$ Neveu-Schwarz
Lie superalgebra and superconformal superfunction.   We recall the
characterization of invertible superconformal functions vanishing at zero
given in \cite{B thesis} and \cite{B memoirs} as exponentials of
superderivations which give a representation of the $N=1$ Neveu-Schwarz
algebra.  This characterization is crucial for the development of the
change of variables formulas.  In Section \ref{Theta-section}, we introduce
two series of supernumbers, $\Theta_j^{(1)}$ and $\Theta_j^{(2)}$, for $j
\in \frac{1}{2} \mathbb{N}$, which are determined by local coordinates
vanishing at zero and infinity, respectively, and which will appear in the
change of variables formulas.  

In Section \ref{NS-VOSA-section}, we recall the notion of $N = 1$
Neveu-Schwarz  vertex operator superalgebra over a Grassmann algebra
($N=1$ NS-VOSA) given in \cite{B thesis} and \cite{B vosas} and  recall some
of the consequences of this notion which we will need later.  

In Section \ref{change-section}, we present the change of variables
formulas for two different types of changes of variables.  First, in Section
\ref{change-zero}, we consider a change of variables by a formal invertible
superconformal function vanishing at zero.  Second, in Section
\ref{change-infinity}, we consider the change of variables associated to a
formal invertible superconformal function vanishing at infinity. In
both cases, we only state the results; we leave the proofs for Sections
\ref{geom} and \ref{alg-section}.  In addition, in Remark
\ref{non-zero-remark} we indicate how to obtain the change of variables
formula for an invertible superconformal change of variables vanishing at a
non-zero point {}from the change of variables for an invertible
superconformal function vanishing at zero.  

In Section \ref{geom} we give geometric proofs of the change of variables
formulas in the case of convergent superconformal changes of variables
vanishing at zero and infinity.  These geometric proofs rely on the
worldsheet supergeometry of propagating superstrings in $N=1$
superconformal field theory as developed in \cite{B memoirs}, the notion of
$N=1$ supergeometric vertex operator superalgebra ($N=1$ SG-VOSA) as
introduced in \cite{B iso-thm}, and the isomorphism between the category of
$N=1$ SG-VOSAs and the category of $N=1$ NS-VOSAs as proved in \cite{B
iso-thm}.  In addition to giving alternate proofs to the formal algebraic
proofs we will present in Section \ref{alg-section}, this section gives the
motivation behind the change of variables formulas.  

Sections \ref{moduli}-\ref{iso} recall the basic structures and results
involved in the correspondence between the algebraic and geometric aspects
of $N=1$ NS-VOSAs.   In Section \ref{moduli}, we recall from \cite{B
memoirs} the moduli space of superspheres with tubes and the sewing
operation and provide two important examples, Examples \ref{example1} and
\ref{example2}, of pairs of  superspheres with tubes and the sewing
operation.  These examples are used in Sections \ref{section-geom1} and
\ref{section-geom2}, respectively, to prove the change of variables
formulas.  In Section
\ref{N=1 SG-VOSA-section}, we recall the notion of $N=1$ SG-VOSA {}from
\cite{B iso-thm}.  In Section \ref{iso}, we recall {}from \cite{B iso-thm}
how to construct an $N=1$ NS-VOSA {}from an $N=1$ SG-VOSA, and how to
construct an
$N=1$ SG-VOSA {}from an $N=1$ NS-VOSA, and then we recall the Isomorphism
Theorem {}from \cite{B iso-thm} which states that the two notions are
equivalent because their respective categories are isomorphic. 
Finally in Sections \ref{section-geom1} and \ref{section-geom2}, we use
Examples \ref{example1} and \ref{example2} and the Isomorphism Theorem to
prove the change of variables formulas for superconformal changes of
variables convergent in a neighborhood of zero and infinity, respectively.

In Section \ref{alg-section}, we give proofs of the change of variables
formulas for formal algebraic, not necessarily convergent, changes of
variables via a formal algebraic argument first employed in
\cite{B iso-thm} to prove the Isomorphism Theorem.  However, in the proof
of the Isomorphism Theorem we assumed that all local changes of variables
were convergent in some neighborhood.  In this section, and in the change
of variables formulas presented in Section \ref{change-section}, we do not
assume the convergence of the changes of variables and consequently cannot
use the Isomorphism Theorem arising {}from the correspondence between the
geometry and algebra as we do in Section \ref{geom}.  However, enough of
the algebraic formalism used in the proof of the Isomorphism Theorem in
\cite{B iso-thm} carries over to the more general setting of formal
not necessarily convergent changes of variables, and in Section
\ref{alg-section} we  retrace that part of the construction  used in \cite{B
iso-thm} in this more general setting to obtain formal algebraic proofs of
the change of variables formulas.

In Section \ref{iso-families-section}, we use the results of the change of
variables formulas presented in Section \ref{change-section} to obtain
families of isomorphic $N=1$ NS-VOSAs. Finally, in Section \ref{annulus}, we
combine the results of Sections \ref{change-section} and
\ref{iso-families-section} to give the change of variables formulas and
family of isomorphic $N=1$ NS-VOSAs for an invertible superconformal change
of variables defined on a superannulus.

These results are super-extensions of the change of variables formulas
developed by Huang in \cite{H book}.  However, there are some typos
in \cite{H book} for defining the isomorphic families of VOAs associated to
the change of variable formulas in a neighborhood of infinity and in an
annulus.  We point out those mistakes in Remarks \ref{H-mistake} and
\ref{H-mistake2}.

\section{Preliminaries}\label{preliminaries}

\subsection{Grassmann algebras and the $N=1$ Neveu-Schwarz algebra}
\label{Grassmann}

Let $\mathbb{Z}_2$ denote the integers modulo two.  For a 
$\mathbb{Z}_2$-graded vector space $V = V^0 \oplus V^1$, define the 
{\it sign function} $\eta$ on the homogeneous subspaces of $V$ by 
$\eta(v) = i$ for $v \in V^i$, $i = 0,1$.  If $\eta(v) = 0$, we say 
that $v$ is {\it even}, and if $\eta(v) = 1$, we say that $v$ is 
{\it odd}.  A {\it superalgebra} is an (associative) algebra $A$ (with
identity $1 \in A$), such that: (i) $A$ is a $\mathbb{Z}_2$-graded algebra;
(ii) $ab = (-1)^{\eta(a) \eta(b)} ba$ for $a,b$ homogeneous in $A$.

The exterior algebra over a vector space $W$, denoted $\bigwedge(W)$, 
has the structure of a superalgebra.  Fix $W_L$ to be an 
$L$-dimensional vector space over $\mathbb{C}$ for $L \in \mathbb{N}$ 
with fixed basis $\{ \zeta_1, ..., \zeta_L \}$ such that $W_L \subset
W_{L+1}$.  We denote $\bigwedge(W_L)$ by $\bigwedge_L$ and call this the
{\it Grassmann algebra on $L$ generators}.  Note that $\bigwedge_L 
\subset \bigwedge_{L+1}$, and taking the direct limit as $L \rightarrow
\infty$, we have the {\it infinite Grassmann algebra} denoted by
$\bigwedge_\infty$.  We use the notation $\bigwedge_*$ to denote a
Grassmann algebra, finite or infinite.    

The $\mathbb{Z}_2$-grading of $\bigwedge_*$ is given explicitly by
\begin{eqnarray*}
\mbox{$\bigwedge_*^0$} \! &=& \! \Bigl\{a \in \mbox{$\bigwedge_*$} \; 
\big\vert \; a = \sum_{(i) \in I_*} a_{(i)}\zeta_{i_{1}}\zeta_{i_{2}} 
\cdots \zeta_{i_{2n}}, \; a_{(i)} \in \mathbb{C}, \; n \in \mathbb{N}
\Bigr\}\\ 
\mbox{$\bigwedge_*^1$} \! &=& \! \Bigl\{a \in \mbox{$\bigwedge_*$} \; 
\big\vert \; a = \sum_{(j) \in J_*} a_{(j)}\zeta_{j_{1}}\zeta_{j_{2}} 
\cdots \zeta_{j_{2n + 1}}, \; a_{(j)} \in \mathbb{C}, \; n \in \mathbb{N}
\Bigr\},
\end{eqnarray*}
where
\begin{eqnarray*}
I_* \! \!  &=&  \! \! \bigl\{ (i) = (i_1, i_2, \ldots, i_{2n}) \; | \; i_1 <
i_2 < \cdots < i_{2n}, \; i_l \in \{1, 2, ..., *\}, \; n \in
\mathbb{N} \bigr\}, \\ 
J_*  \! \! &=&  \! \! \bigl\{(j) = (j_1, j_2, \ldots, j_{2n + 1}) \; | \;
j_1 < j_2 < \cdots < j_{2n + 1}, \; j_l \in \{1, 2, ..., *\}, \;  n \in
\mathbb{N} \bigr\},
\end{eqnarray*}
with $\{1,2,...,*\}$ denoting $\{1,2,...,L\}$ if $\bigwedge_*$ is the
finite-dimensional Grassmann algebra $\bigwedge_L$ and denoting the 
positive integers $\mathbb{Z}_+$ if $\bigwedge_* = \bigwedge_\infty$. We
can also decompose $\bigwedge_*$ into {\it body}, $(\bigwedge_*)_B = 
\{ a_{(\emptyset)} \in \mathbb{C} \}$,  and {\it soul} 
\[(\mbox{$\bigwedge_*$})_S \; = \; \Bigl\{a \in \mbox{$\bigwedge_*$} \; \big\vert \; 
a = \! \! \! \sum_{ \begin{scriptsize} \begin{array}{c}
(k) \in I_* \cup J_*\\
k \neq (\emptyset)
\end{array} \end{scriptsize}} \! \! \!
a_{(k)} \zeta_{k_1} \zeta_{k_2} \cdots \zeta_{k_n}, \; a_{(k)} \in \mathbb{C} 
\Bigr\}\] 
subspaces such that $\bigwedge_* = (\bigwedge_*)_B \oplus 
(\bigwedge_*)_S$.  For $a \in \bigwedge_*$, we write $a = a_B + 
a_S$ for its body and soul decomposition.

A $\mathbb{Z}_2$-graded vector space $\mathfrak{g}$ is said to be a
{\it Lie superalgebra} if it has a bilinear operation $[\cdot,\cdot]$
such that for $u,v$ homogeneous in $\mathfrak{g}$: (i) $[u,v] \in
\mathfrak{g}^{(\eta(u) + \eta(v)) \mathrm{mod} \; 2}$; (ii)
skew-symmetry holds $[u,v] = -(-1)^{\eta(u)\eta(v)}[v,u]$; (iii) the
Jacobi identity holds $(-1)^{\eta(u)\eta(w)}[[u,v],w] +
(-1)^{\eta(v)\eta(u)}[[v,w],u] +  (-1)^{\eta(w)\eta(v)}[[w,u],v] = 0$.

For any $\mathbb{Z}_2$-graded associative algebra $A$ and for $u,v \in A$ 
of homogeneous sign, we can define $[u,v] = u v - (-1)^{\eta(u)\eta(v)} v 
u$, making $A$ into a Lie superalgebra.  The algebra of endomorphisms of 
$A$, denoted $\mbox{End} \; A$, has a natural $\mathbb{Z}_2$-grading 
induced {}from that of $A$, and defining $[X,Y] = X Y - (-1)^{\eta(X)
\eta(Y)} Y X$ for $X,Y$ homogeneous in $\mbox{End} \; A$, this gives 
$\mbox{End} \; A$ a Lie superalgebra structure.  An element $D \in 
(\mbox{End} \; A)^i$, for $i \in \mathbb{Z}_2$, is called a {\it 
superderivation of sign $i$} (denoted $\eta(D) = i$) if $D$ satisfies the 
super-Leibniz rule $D(uv) = (Du)v + (-1)^{\eta(D) \eta(u)} uDv$
for $u,v \in A$ homogeneous.

The {\it Virasoro algebra} is the Lie algebra with basis consisting of the
central element $d$ (called the {\it central charge}) and $L_n$, for $n \in
\mathbb{Z}$, satisfying commutation relations  
\begin{equation}\label{Virasoro relation}
[L_m ,L_n] = (m - n)L_{m + n} + \frac{1}{12} (m^3 - m) \delta_{m + n 
, 0} \; d ,
\end{equation}
for $m, n \in \mathbb{Z}$.  The {\it $N=1$ Neveu-Schwarz Lie superalgebra}
is a super-extension of the Virasoro algebra by the odd elements 
$G_{n + 1/2}$, for $n \in \mathbb{Z}$, with supercommutation relations 
\begin{eqnarray}
\left[ G_{m + \frac{1}{2}},L_n \right] &=& \Bigl(m - \frac{n - 1}{2} \Bigr) G_{m
+ n + \frac{1}{2}}  \label{Neveu-Schwarz relation1} \\  
\left[ G_{m + \frac{1}{2}} , G_{n - \frac{1}{2}} \right] &=& 2L_{m +
n} + \frac{1}{3} (m^2 + m) \delta_{m + n , 0} \; d \label{Neveu-Schwarz relation2}
\end{eqnarray}
in addition to (\ref{Virasoro relation}).  We denote the $N=1$
Neveu-Schwarz algebra by $\mathfrak{ns}$.

Let $V = \coprod_{k \in (1/2) \mathbb{Z}} V_{(k)}$ be a module for
$\mathfrak{ns}$, and let $L(n), G(n - 1/2) \in \mbox{End} \; V$ and $c \in
\mathbb{C}$ be the representation images of $L_n$, $G_{n - 1/2}$ and $d$,
respectively, such that $L(0)v = kv$ for $v \in V_{(k)}$.  If in addition,
$V_{(k)} = 0$ for $k$ sufficiently small, we call $V$ a {\it positive
energy representation} of $\mathfrak{ns}$.

\subsection{Superanalytic and superconformal superfunctions}

Let $U$ be a subset of $\bigwedge_*$, and write $U = U^0 \oplus U^1$
for the decomposition of $U$ into even and odd subspaces.  Let $z$ 
be an even variable in $U^0$ and $\theta$ an odd variable in $U^1$.  We
call $H: U \longrightarrow \bigwedge_*$, mapping $(z,\theta) \mapsto
H(z,\theta)$, a {\it $\bigwedge_*$-superfunction on $U$ in
(1,1)-variables} or just a {\it superfunction}.

Let $z_B$ be a complex variable and $h(z_B)$ a complex
analytic function in some open set $U_B \subset \mathbb{C}$.  For $z$ a
variable in $\bigwedge_*^0$, we define $h(z)$ to be the Taylor expansion
about the body of $z = z_B + z_S$. Then $h(z)$ is well defined (i.e.,
convergent) in the open neighborhood $\{z = z_B + z_S \in \bigwedge_*^0 \; |
\; z_B \in U_B \} = U_B \times(\bigwedge_*^0)_S \subseteq
\bigwedge_*^0$.   Since $h(z)$ is algebraic in each $z_{(i)}$, for $(i)
\in I_*$, it follows that $h(z)$ is complex analytic in each of the
complex variables $z_{(i)}$.

For $n \in \mathbb{N}$, we introduce the notation $\bigwedge_{*>n}$ to
denote a finite Grassmann algebra $\bigwedge_L$ with $L > n$ or an 
infinite Grassmann algebra.  We will use the corresponding index 
notations for the corresponding indexing sets $I_{*>n}$ and $J_{*>n}$.
A {\it superanalytic $\bigwedge_{*>0}$-superfunction in
(1,1)-variables} $H$ is a superfunction in of the form
\begin{eqnarray*}
H(z, \theta) \! \! &=& \! \! (f(z) + \theta \xi(z), \psi(z) + \theta
g(z)) \\ 
&=& \! \! \!  \Biggl(\sum_{(i) \in I_{ * - 1}} \!
\! \! f_{(i)} (z) \zeta_{i_1} \zeta_{i_2} \cdots \zeta_{i_n} + \;
\theta \! \! \! \sum_{(j) \in J_{ * - 1}} \! \! \!
\xi_{(j)} (z) \zeta_{j_1} \zeta_{j_2} \cdots \zeta_{j_{2n + 1}},
\Biggr. \\ 
& & \qquad \quad \Biggl. \sum_{(j) \in J_{ * - 1}} \!  \! \!
\psi_{(j)} (z) \zeta_{j_1} \zeta_{j_2} \cdots \zeta_{j_{2n + 1}} + \;
\theta \! \! \! \sum_{(i) \in I_{ * - 1}} \! \! \!
g_{(i)} (z) \zeta_{i_1} \zeta_{i_2} \cdots \zeta_{i_n} \Biggr)
\end{eqnarray*}
where $f_{(i)}(z_B)$, $g_{(i)}(z_B)$, $\xi_{(j)}(z_B)$, and
$\psi_{(j)}(z_B)$ are all complex analytic in some non-empty open
subset $U_B \subseteq \mathbb{C}$.  If each $f_{(i)}(z_B)$,
$g_{(i)}(z_B)$, $\xi_{(j)}(z_B)$,  and $\psi_{(j)}(z_B)$ is complex
analytic in $U_B \subseteq \mathbb{C}$, then $H(z,\theta)$ is well
defined (i.e., convergent) for $\{(z,\theta) \in \bigwedge_{*>0} \; | \;
z_B \in U_B \} = U_B \times (\bigwedge_{*>0})_S$.  Consider the topology
on $\bigwedge_{*}$ given by the product of the usual topology on 
$(\bigwedge_{*})_B = \mathbb{C}$ and the trivial topology on 
$(\bigwedge_{*})_S$.  This topology on $\bigwedge_{*}$ is called 
the {\it DeWitt topology}.  The natural domain of any superanalytic
$\bigwedge_{*>0}$-superfunction is an open set in the DeWitt topology 
on $\bigwedge_{*>0}$.

Since $1/a = \sum_{n \in \mathbb{N}} (-1)^n a_S^n/a_B^{n + 1}$ is well
defined if and only if $a_B \neq 0$, the set of invertible elements in
$\bigwedge_*$, denoted $\bigwedge_*^\times$, is given by
$\bigwedge_*^\times = \{a \in \mbox{$\bigwedge_*$} \; | \; a_B \neq 0
\}$.

We define the (left) partial derivatives $\partial / \partial z$ and
$\partial / \partial \theta $ acting on superfunctions which are
superanalytic in some DeWitt open neighborhood $U$ of $(z,\theta) \in
\bigwedge_{*>0}$ by
\begin{eqnarray*}
\Delta z \left(\dz H(z,\theta) \right) + O((\Delta z)^2) &=& H(z +
\Delta z, \theta) - H(z,\theta) \\
\Delta \theta \left(\dtheta H(z,\theta) \right) &=& H(z, \theta +
\Delta \theta) - H(z,\theta) 
\end{eqnarray*}
for all $\Delta z \in \bigwedge_{*>0}^0$ and $\Delta \theta \in
\bigwedge_{*>0}^1$ such that $z + \Delta z \in U^0 = U_B \times
(\bigwedge_{*>0}^0)_S$ and $\theta + \Delta \theta \in U^1 =
\bigwedge_{*>0}^1$.  Note that $\partial / \partial z$ and $\partial /
\partial \theta$ are endomorphisms of the superalgebra of superanalytic
superfunctions, and in fact, are  even and odd superderivations,
respectively.

If $h(z_B)$ is complex analytic in an open neighborhood of the complex
plane, then $h(z_B)$ has a Laurent series expansion in $z_B$, given by 
$h(z_B) = \sum_{l \in \mathbb{Z}} c_l z_B^l$, for $c_l \in \mathbb{C}$, 
and we have $h(z) = \sum_{n \in \mathbb{N}} (z_S^n/n!) h^{(n)}(z_B) =
\sum_{l \in \mathbb{Z}} c_l (z_B + z_S)^l = \sum_{l \in \mathbb{Z}} c_l z^l$
where $(z_B + z_S)^l$, for $l \in \mathbb{Z}$, is always understood to 
mean expansion in positive powers of the second variable, in this
case $z_S$.  Thus if $H$ is a $\bigwedge_{*>0}$-superfunction in 
$(1,1)$-variables which is superanalytic in a (DeWitt) open 
neighborhood, $H$ can be expanded as
\begin{equation}\label{Laurent series}
H(z, \theta) = \biggl(\sum_{l \in \mathbb{Z}} a_l z^l + \theta \sum_{l
\in \mathbb{Z}} n_l z^l, \sum_{l \in \mathbb{Z}} m_l z^l + \theta
\sum_{l \in \mathbb{Z}} b_l z^l \biggr)
\end{equation}   
for $a_l, b_l \in \bigwedge_{ * - 1}^0$ and $m_l, n_l \in
\bigwedge_{ * - 1}^1$.  

Define $D$ to be the odd superderivation $D = \partial / \partial \theta +
\theta \partial / \partial z$ acting on superanalytic superfunctions.  Then
$D^2 = \partial / \partial z$, and if $H(z,\theta) =
(\tilde{z},\tilde{\theta})$ is superanalytic in some  DeWitt open subset,
then $D$ transforms under $H(z,\theta)$ by $D = (D\tilde{\theta})\tilde{D} +
(D\tilde{z} - \tilde{\theta} D \tilde{\theta})\tilde{D}^2$.  We define a
{\it superconformal $(1,1)$-superfunction} on a DeWitt open  subset $U$ of
$\bigwedge_{*>0}$ to be a superanalytic superfunction $H$ under which $D$
transforms homogeneously of degree one.  Thus a superanalytic function
$H(z,\theta) =  (\tilde{z}, \tilde{\theta})$ is superconformal if and only
if, in  addition to being superanalytic, $H$ satisfies $D\tilde{z} -
\tilde{\theta} D\tilde{\theta} = 0$, for $D \tilde{\theta}$ not identically
zero, thus transforming $D$ by $D = (D\tilde{\theta})\tilde{D}$.  

Let $R$ be a superalgebra, let $x$ be an even formal variable, and let 
$\varphi$ be an odd formal variable.  By this we mean that $x$ commutes 
with all formal variables and all superalgebra elements and $\varphi$ 
anticommutes with all odd formal variables and all odd superalgebra 
elements but commutes with even elements.  We use the following notational
conventions for $V$ a vector space:
\begin{eqnarray*}
V[x] &=& \bigl\{  \sum_{j \in \mathbb{N}} v_j x^j \; | \; v_j \in V, \;
\mbox{all but finitely many $v_j = 0$} \bigr\} \\
V[x,x^{-1}] &=& \bigl\{  \sum_{j \in \mathbb{Z}} v_j x^j \; | \; v_j \in V,
\; \mbox{all but finitely many $v_j = 0$} \bigr\} \\
V[[x]] &=& \bigl\{  \sum_{j \in \mathbb{N}} v_j x^j \; | \; v_j \in
V,\bigr\} \\
V[[x,x^{-1}]] &=& \bigl\{  \sum_{j \in \mathbb{Z}} v_j x^j \; |
\; v_j \in V, \bigr\} 
\end{eqnarray*}
\begin{eqnarray*}
V((x)) &=& \bigl\{  \sum_{j \in \mathbb{Z}} v_j x^j \; |
\; v_j \in V, \; \mbox{$v_j = 0$ for $j$ sufficiently small} \bigr\} \\
V[\varphi] &=& V \oplus \varphi V.
\end{eqnarray*}
We will use analogous notation for several variables.  Note that
$R((x))[\varphi]$ and $R((x^{-1}))[\varphi]$ are superalgebras.  

For $j \in \mathbb{Z}$, define the following even and odd
superderivations, respectively, 
\begin{eqnarray}
L_j(x,\varphi) &=& - \biggl( \Lx \biggr) \label{L notation} \\
G_{j -\frac{1}{2}} (x,\varphi) &=& - \Gx \label{G notation}
\end{eqnarray}
in $\mbox{Der} (R[[x,x^{-1}]]))[\varphi])$.  It is easy to check that
$L_j(x,\varphi)$ and $G_{j-1/2}(x,\varphi)$ give a representation of the
$N=1$ Neveu-Schwarz superalgebra with central charge zero.

For any formal power series $H \in R[[x,x^{-1}]][\varphi]$, we say that 
$H(x,\varphi) = (\tilde{x},\tilde{\varphi})$ is {\it superconformal} if
$D\tilde{x} = \tilde{\varphi}D\tilde{\varphi}$ for
$D = \partial/\partial \varphi + \varphi \partial/\partial x$. 

Let $(R^0)^\infty$ be the set of all sequences $A = \{A_j\}_{j \in \Z}$ of
even  elements in $R$, let $(R^1)^\infty$ be the set of all sequences 
$M = \{M_{j - 1/2}\}_{j \in \Z}$ of odd elements in $R$, and let $R^\infty =
(R^0)^\infty \oplus (R^1)^\infty$.  Let $(R^0)^\times$ denote the set of
invertible elements in $R^0$.  Consider the set of formal superconformal
power series $H \in xR[[x]] \oplus \varphi R[[x]]$ with invertible even
coefficient of $\varphi$, i.e., of the form
\begin{multline}\label{formal Laurent series}
H(x, \varphi) = \biggl( \asqrt^2 \Bigl( x + \sum_{j \in \Z} a_j
x^{j+1} + \varphi \sum_{j \in \Z} n_j x^j \Bigr), \\
\asqrt \Bigl( \sum_{j \in \Z} m_{j- \frac{1}{2}} x^j + \varphi ( 1 + \sum_{j
\in \Z} b_j x^j) \Bigr) \biggr) 
\end{multline}  
with $\asqrt \in (R^0)^\times$, $(a,m) = \{(a_j, m_{j - \frac{1}{2}}) \}_{j
\in \Z} \in R^\infty$ and satisfying $D\tilde{x} = \tilde{\varphi} D
\tilde{\varphi}$.  If $R  = \bigwedge_*$ and $(x,\varphi) = (z,\theta)$
then these are the invertible superconformal local coordinates vanishing at
zero.  Note that if $H$ is given by (\ref{formal Laurent series}), then the
condition that $H$ be superconformal means that $n_j$ and $b_j$, for $j \in
\Z$, are completely determined by $\asqrt$, the $a_j$'s and the $m_{j -
\frac{1}{2}}$'s, for $j \in \Z$.   In \cite{B thesis} and \cite{B memoirs},
we show that there is a bijection between the set of formal superconformal
power series $H \in xR[[x]] \oplus \varphi R[[x]]$ with invertible even
coefficient of $\varphi$ and the set of formal power series of the form
\begin{equation}
\exp\Biggl( \! - \! \sum_{j \in \Z} \Bigl( A_j L_j(x,\varphi) + M_{j -
\frac{1}{2}} G_{j -\frac{1}{2}} (x,\varphi) \Bigr) \! \Biggr) \cdot
(\asqrt^2 x, \asqrt \varphi)  
\end{equation}
with $\asqrt \in (R^0)^\times$, and $(A,M) = \{(A_j, M_{j - \frac{1}{2}})
\}_{j \in \Z} \in R^\infty$.  As is shown in \cite{B memoirs}, this defines
a bijection 
\begin{eqnarray}
E : R^\infty &\longrightarrow& R^\infty\\
(A,M) &\mapsto& (a,m). \nonumber
\end{eqnarray}
which allows us to define a map $\tilde{E}$ {}from $R^\infty$ to  the set of
all formal superconformal power series in $xR[[x]] \oplus \varphi R[[x]]$
with leading even coefficient of $\varphi$ equal to one, by defining
\begin{eqnarray}
\varphi \tilde{E}^0(A,M)(x,\varphi) &=& \varphi \biggl( x + \sum_{j \in
\Z} E_j(A,M) x^{j + 1} \biggr), \label{Etilde1} \\
\varphi \tilde{E}^1(A,M)(x,\varphi) &=& \varphi \sum_{j \in \Z} E_{j -
\frac{1}{2}}(A,M) x^j , \label{Etilde2}
\end{eqnarray}
and letting $\tilde{E}(A,M) (x, \varphi)$ be the unique formal
superconformal power series with even coefficient of $\varphi$ equal to one 
such that the even and odd components of $\tilde{E}$ satisfy (\ref{Etilde1})
and (\ref{Etilde2}), respectively.  Similarly, we define a map $\hat{E}$
{}from $(R^0)^\times \times R^\infty$ to the set of all formal
superconformal power series in $xR[[x]] \oplus \varphi R[[x]]$ with
invertible leading even coefficient of $\varphi$, by defining 
\begin{eqnarray}
\hat{E}(\asqrt,A,M)(x,\varphi) &=&  (\hat{E}^0(\asqrt, A, M)
(x,\varphi), \hat{E}^1(\asqrt, A, M) (x,\varphi)) \\ \nonumber
&=& (\asqrt^2 \tilde{E}^0(A,M)(x,\varphi), \asqrt
\tilde{E}^1(A,M)(x,\varphi).  
\end{eqnarray}
Thus we have the following proposition which is proved in 
\cite{B memoirs}.

\begin{prop}\label{above}(\cite{B memoirs})
The map $\hat{E}$ {}from $(R^0)^{\times} \times R^{\infty}$ to the
set of all formal superconformal power series $H(x, \varphi) \in
xR[[x]] \oplus \varphi R[[x]]$ of the form
\begin{equation}\label{Ehat} 
\varphi H(x,\varphi) = \varphi \Biggl(\asqrt^2 \Bigl( x + \sum_{j \in \Z}
a_j x^{j + 1} \Bigr), \asqrt \sum_{j \in \Z} m_{j - \frac{1}{2}} x^j
\Biggr)  
\end{equation} 
and with even coefficient of $\varphi$ equal to $\asqrt$ 
for $(\asqrt,a,m) \in (R^0)^{\times} \times R^\infty$, is a
bijection.

The map $\tilde{E}$ {}from $R^{\infty}$ to the set of
formal superconformal power series of the form (\ref{Ehat}) with
$\asqrt = 1$ and even coefficient of $\varphi$ equal to one is also a 
bijection. 

In particular, we have inverses $\tilde{E}^{-1}$ and $\hat{E}^{-1}$.
\end{prop}

Note: In \cite{B memoirs}, the space $xR[[x]] \oplus \varphi R[[x]]$ was
erroneously denoted $xR[[x]][\varphi]$.

\begin{rema}\label{corollary to above} {\em Let $I(x,\varphi)  = (1/x,
i\varphi/x)$.  Then $I$ is superconformal and vanishing at $(x,\varphi) =
(\infty,0) = \infty$.  Note that the composition of two superconformal
functions is again superconformal.  Thus if $H \in x^{-1}R[[x^{-1}]]
[\varphi]$ is superconformal, vanishing at infinity and with even
coefficient of $\varphi x^{-1}$ equal to $i$, then $H \circ I^{-1}$ is
superconformal, vanishing at zero and with even coefficient of $\varphi$
equal to one.  Then by Proposition \ref{above}, $H\circ I^{-1} (x,\varphi) =
\tilde{E}(A,M)(x,\varphi)$, and therefore 
\begin{eqnarray*}
H(x,\varphi) &=& \tilde{E}(A,M) \circ I (x,\varphi) \\
&=& \Bigl. \exp\Biggl( \! - \! \sum_{j \in \Z} \Bigl( A_j L_j(x,\varphi) +
M_{j - \frac{1}{2}} G_{j -\frac{1}{2}} (x,\varphi) \Bigr) \! \Biggr) \cdot
(x,\varphi) \Bigr|_{(x,\varphi) = I(x,\varphi)} \\
&=&  \exp\Biggl( \sum_{j \in \Z} \Bigl( A_j L_{-j}(x,\varphi) + i
M_{j - \frac{1}{2}} G_{-j +\frac{1}{2}} (x,\varphi) \Bigr) \! \Biggr) \cdot
(x,\varphi) .
\end{eqnarray*}
}
\end{rema}

\subsection{The $\{\Theta^{(1)}_j,\Theta^{(1)}_{j- 1/2})\}_{j \in
\mathbb{Z}}$ and $\{\Theta^{(2)}_j,\Theta^{(2)}_{j-1/2})\}_{j \in
\mathbb{Z}}$ series} \label{Theta-section}

In this section we recall two series first introduced
in \cite{B memoirs} which will play a central role in the change of
variables formulas.  Though rather technical as introduced here, a
geometric motivation for these series will be given in Section \ref{geom}.

Let $\alpha_0^{1/2}$ and $\mathcal{A}_j$, for $j \in \Z$, be even formal
variables and let $\mathcal{M}_{j - \frac{1}{2}}$, for $j \in \Z$, be
odd formal variables.  Let
\begin{eqnarray*}
H^{(1)}_{\alpha_0^{1/2},\mathcal{A}, \mathcal{M}} (x, \varphi)  &=& 
\exp \Biggl(\!  - \! \sum_{j \in \Z}
\biggl( \mathcal{A}_j L_j(x,\varphi)  +  \mathcal{M}_{j - \frac{1}{2}}
G_{j - \frac{1}{2}}(x,\varphi) \biggr) \! \Biggr) \cdot \\
& & \hspace{2in} \cdot (\alpha_0^{1/2})^{-2L_0(x,\varphi)}\cdot (x,
\varphi ) \\ 
&=& \hat{E}(\alpha_0^{1/2},\mathcal{A}, \mathcal{M} ) (x,\varphi).
\end{eqnarray*}
Let 
\[(\tilde{x}, \tilde{\varphi}) = (H^{(1)}_{\alpha_0^{1/2},\mathcal{A}, 
\mathcal{M}})^{-1} (x,\varphi) \in (\alpha_0^{-1} x,\alpha_0^{-\frac{1}{2}}
\varphi) + x\mathbb{C}[x,\varphi][\alpha_0^{-\frac{1}{2}}][[\mathcal{A}]]
[\mathcal{M}] . \]
Let $w$ be another even formal variable and $\rho$ another odd formal
variable and define the superconformal shift
\begin{eqnarray}
s_{(x,\varphi)} (w,\rho) =  (w - x -\rho\varphi,\rho - \varphi).
\end{eqnarray}
Then $s_{(x,\varphi)} \circ H^{(1)}_{\alpha_0^{1/2},\mathcal{A},
\mathcal{M}}
\circ s_{(\tilde{x}, \tilde{\varphi})}^{-1} (\alpha_0^{-1} w, 
\alpha_0^{-1/2}\rho)$ is in 
\[w\mathbb{C} [x, \varphi] [\alpha_0^{\frac{1}{2}}, \alpha_0^{-\frac{1}{2}}] 
[[\mathcal{A}]] [\mathcal{M}][[w]] \oplus \rho \mathbb{C}[x, \varphi] 
[\alpha_0^{\frac{1}{2}}, \alpha_0^{-\frac{1}{2}}][[\mathcal{A}]][\mathcal{M}][[w]] 
,\] 
is superconformal in $(w,\rho)$, and the even coefficient of the monomial
$\rho$ is an element in $(1 + x\mathbb{C}[x][\alpha_0^{1/2},
\alpha_0^{-1/2}] [[\mathcal{A}]][\mathcal{M}] \oplus \varphi \mathbb{C}
[x][\alpha_0^{1/2}, \alpha_0^{-1/2}] [[\mathcal{A}]][\mathcal{M}])$.

Let $\Theta^{(1)}_j = \Theta^{(1)}_j(\alpha_0^{1/2},\mathcal{A}, 
\mathcal{M}, (x, \varphi)) \in \mathbb{C}[x, \varphi] [\alpha_0^{1/2},
\alpha_0^{-1/2}][[\mathcal{A}]][\mathcal{M}]$, for $j \in \frac{1}{2}
\mathbb{N}$, be defined by
\begin{multline}\label{define first Theta}
\Bigl(\exp(\Theta^{(1)}_0(\alpha_0^{\frac{1}{2}}, \mathcal{A}, \mathcal{M}, 
(x, \varphi)), \Bigl\{\Theta^{(1)}_j(\alpha_0^{\frac{1}{2}}, \mathcal{A}, 
\mathcal{M}, (x, \varphi)), \Bigr. \Bigr. \\
\Bigl. \Bigl. \Theta^{(1)}_{j - \frac{1}{2}} (\alpha_0^{\frac{1}{2}}, 
\mathcal{A}, \mathcal{M}, (x, \varphi)) \Bigr\}_{j \in \Z} \Bigr)
\end{multline}  
\[ = \;\hat{E}^{-1}(s_{(x,\varphi)} \circ H^{(1)}_{\alpha_0^{1/2},\mathcal{A}, 
\mathcal{M}} \circ s_{(\tilde{x},\tilde{\varphi})}^{-1}
(\alpha_0^{-1}w,\alpha_0^{-\frac{1}{2}}\rho)) .  
\hspace{1.5in} \]
In other words, the $\Theta_j^{(1)}$'s are determined uniquely by
\begin{multline*}
s_{(x,\varphi)} \circ H^{(1)}_{\alpha_0^{1/2},\mathcal{A}, \mathcal{M}}
\circ s_{(\tilde{x}, \tilde{\varphi})}^{-1} (\alpha_0^{-1} w, 
\alpha_0^{-\frac{1}{2}}\rho) \\  
= \exp \Biggl( \sum_{j \in  \Z} \biggl( \Theta^{(1)}_j \! \left( \Lw \right)
+ \Theta^{(1)}_{j - \frac{1}{2}} \Gw \biggr) \! \Biggr) \! \cdot \\
\exp \left(\! \Theta^{(1)}_0 \left( \twoLow \right) \! \right) \! \cdot
(w,\rho) . 
\end{multline*}  

This formal power series in $(w,\rho)$ gives the formal local 
superconformal coordinate at a puncture of the canonical supersphere 
obtained {}from the sewing together of two particular canonical 
superspheres with punctures as we shall in Example \ref{example1} in
Section \ref{moduli}.

The following proposition states results contained in Corollary 3.44 and
Proposition 3.45 in \cite{B memoirs}.  Corollary 3.44 and Proposition 3.45
are in turn based on Proposition 3.30 of \cite{B memoirs}.  However, there
is a typo on line six on p.72 in \cite{B memoirs} of the proof of
Proposition 3.30 and hence in the statement of Proposition 3.30 and in the
statements of Corollary 3.44 and Proposition 3.45 which are based on
Proposition 3.30.  Line six on p.72 of the proof of Proposition 3.30 in
\cite{B memoirs} should read 
\[\exp \left( -( \alpha_0 \tilde{x} - x - \alpha_0^{\frac{1}{2}}
\tilde{\varphi}\varphi)\frac{\partial}{\partial w}  - (\alpha^{\frac{1}{2}}
\tilde{\varphi} - \varphi) \Bigl(\frac{\partial}{\partial \rho} - \rho
\frac{\partial}{\partial w} \Bigr) \right)  \]
rather than
\[\exp \left( -(\tilde{x} - \alpha_0^{-1} x)\frac{\partial}{\partial w}  -
(\tilde{\varphi} - \alpha_0^{-\frac{1}{2}} \varphi) \Bigl(
\frac{\partial}{\partial \rho} - \rho \frac{\partial}{\partial w} \Bigr)
\right) .\] 
We will need the corrected results {}from Corollary 3.44 and
Proposition 3.45 in \cite{B memoirs}, and these results are given in the
proposition below.  It is the first line of the righthand side of equation
(\ref{last first Theta equation}) below that contains the corrected
terms.

\begin{prop}(\cite{B memoirs}) \label{first Theta identity in End}
The formal series $\Theta^{(1)}_j (t^{ - 1/2} \alpha_0^{1/2}, \mathcal{A}, 
\mathcal{M}, (x,\varphi))$, for $j \in \frac{1}{2} \mathbb{N}$, are in 
$\mathbb{C} [x, \varphi][\mathcal{A}][\mathcal{M}][\alpha_0^{-1/2}]
[[t^{1/2}]]$.  Thus for $\asqrt \in (\bigwedge_*^0)^\times$, and $(A,M) \in
\bigwedge_*^\infty$, the series $\Theta^{(1)}_j (t^{ - 1/2} \asqrt,
A,M, (x,\varphi))$,  for $j \in \frac{1}{2}\mathbb{N}$, are well defined
and belong to $\bigwedge_*  [x, \varphi][[t^{1/2}]]$. 

Let $(\tilde{x}(t^{1/2}),\tilde{\varphi} (t^{1/2})) = 
(H^{(1)}_{ t^{ - 1/2}\asqrt,A,M})^{-1}(x,\varphi)$, and let $V$ be a
positive-energy module for $\mathfrak{ns}$.  Then in
$((\mbox{\em End} \; (\bigwedge_* \otimes_\mathbb{C} V))
[[t^{\frac{1}{2}}]] [[x,x^{-1}]][\varphi])^0$, the following identity
holds for $\Theta^{(1)}_j(t^{1/2}) = \Theta^{(1)}_j (t^{-1/2}\asqrt,
A,M,(x,\varphi))$. 
\begin{multline}\label{last first Theta equation}
\exp \Biggl(\! - \! \! \sum_{m = -1}^{\infty} \sum_{j \in \Z}
\binom{j+1}{m+1} t^j \asqrt^{-2j} x^{j - m}  \\
\biggl(\!  \Bigl( A_j + 2\left(\frac{j-m}{j+1} \right)
t^{-\frac{1}{2}} \asqrt x^{-1} \varphi M_{j - \frac{1}{2}}
\Bigr)  L(m) \\
+ \; x^{-1} \Bigl( \left(\frac{j-m}{j+1} \right) t^{-\frac{1}{2}}
\asqrt M_{j - \frac{1}{2}} + \varphi \frac{(j-m)}{2} A_j
\Bigr)  G(m + \frac{1}{2}) \biggl) \Biggr) 
\end{multline}
\begin{multline*}
= \; \exp \left( \! (t^{-1} \asqrt^{2} \tilde{x}(t^\frac{1}{2}) -
x - t^{-\frac{1}{2}} \asqrt \tilde{\varphi}(t^{\frac{1}{2}}) \varphi) L(-1)
+ (t^{-\frac{1}{2}} \asqrt \tilde{\varphi} (t^\frac{1}{2}) - \varphi ) 
G(-\frac{1}{2}) \! \right) \cdot \\
\exp \Biggl( \! - \! \sum_{j \in \Z} \Bigl( \Theta^{(1)}_j (t^{\frac{1}{2}})
L(j) + \Theta^{(1)}_{j - \frac{1}{2}} (t^{\frac{1}{2}}) G(j - \frac{1}{2})
\Bigr) \! \Biggr) \! \cdot \exp \left( - 2\Theta^{(1)}_0 (t^{\frac{1}{2}})
L(0) \right) . 
\end{multline*}
\end{prop}

Similarly, let $\mathcal{B}_j$, for $j \in \Z$, be even formal variables and
let $\mathcal{N}_{j - \frac{1}{2}}$, for $j \in \Z$, be odd formal
variables.  Let 
\begin{eqnarray*}
H^{(2)} (x, \varphi) \! &=& \! \exp \Biggl(\sum_{j \in \Z} \biggl(
\mathcal{B}_j L_{-j}(x,\varphi) +  \mathcal{N}_{j - \frac{1}{2}} G_{-j +
\frac{1}{2}}(x,\varphi) \biggr) \! \Biggr) \cdot \Bigl(\frac{1}{x}, \frac{i
\varphi}{x} \Bigr) \\ 
&=& \tilde{E}(\mathcal{B},-i \mathcal{N}) \Bigl(\frac{1}{x},
\frac{i \varphi}{x} \Bigr) 
\end{eqnarray*}
and let
\[(\tilde{x}, \tilde{\varphi}) = (H^{(2)}_{\mathcal{B}, \mathcal{N}})^{-1}
\circ I (x,\varphi) \in (x,\varphi) + \mathbb{C}[x^{-1}, \varphi]
[[\mathcal{B}]][\mathcal{N}] . \]
Let $w$ be another even formal variable and $\rho$ another odd formal
variable.  Now write $s_{(x,\varphi)}(w,\rho) = (-x + w - \rho 
\varphi, \rho - \varphi)$.  We will use the convention that we should expand
$(-x + w -\rho \varphi)^j = (-x + w)^j - j\rho \varphi (-x + w)^{j-1}$ 
in positive powers of the second even variable $w$, for $j \in
\mathbb{Z}$. Then 
\[s_{(x,\varphi)} \circ I^{-1} \circ H^{(2)}_{\mathcal{B},\mathcal{N}}
\circ s_{(\tilde{x}, \tilde{\varphi})}^{-1} (w,\rho)  \in  
w\mathbb{C} [x^{-1} \!, \varphi][[\mathcal{B}]][\mathcal{N}][[w]] \oplus 
\rho \mathbb{C} [x^{-1} \!, \varphi][[\mathcal{B}]][\mathcal{N}][[w]] ,\]
is superconformal in $(w,\rho)$, and the even coefficient of the monomial
$\rho$ is an element in $(1 + x^{-1} \mathbb{C} [x^{-1},
\varphi][[\mathcal{B}]][\mathcal{N}])$.  

Let $\Theta^{(2)}_j = \Theta^{(2)}_j(\mathcal{B}, \mathcal{N}, (x, \varphi))
\in \mathbb{C}[x^{-1}, \varphi][[\mathcal{B}]][\mathcal{N}]$, 
for $j \in \frac{1}{2} \mathbb{N}$, be defined by
\begin{equation}\label{define second Theta}
\Bigl(\exp(\Theta^{(2)}_0(\mathcal{B}, \mathcal{N}, (x, \varphi)), \Bigl\{
\Theta^{(2)}_j(\mathcal{B},\mathcal{N}, (x, \varphi)),
\Theta^{(2)}_{j - \frac{1}{2}} (\mathcal{B},
\mathcal{N}, (x, \varphi)) \Bigr\}_{j \in \Z} \Bigr) 
\end{equation}
\[ = \; \hat{E}^{-1}(s_{(x,\varphi)} \circ I^{-1} \circ
H^{(2)}_{\mathcal{B}, \mathcal{N}} \circ s_{(\tilde{x},
\tilde{\varphi})}^{-1} (w,\rho)) . \hspace{1.6in}  \] 
In other words, the $\Theta_j^{(2)}$'s are determined uniquely by
\begin{multline*}
s_{(x,\varphi)} \circ I^{-1} \circ H^{(2)}_{\mathcal{B}, 
\mathcal{N}} \circ s_{(\tilde{x}, \tilde{\varphi})}^{-1} (w,\rho) \\ 
= \exp \Biggl( \sum_{j \in  \Z} \biggl( \Theta^{(2)}_j \! \left( \Lw \right)
+ \Theta^{(2)}_{j - \frac{1}{2}} \Gw \biggr) \! \Biggr) \! \cdot \\
\exp \left(\! \Theta^{(2)}_0 \left( \twoLow \right) \! \right) \! 
\cdot (w,\rho) . 
\end{multline*}

This formal power series in $(w,\rho)$ gives the formal local 
superconformal coordinate at a puncture of the canonical supersphere 
obtained {}from the sewing together of two particular canonical 
superspheres with punctures as we shall see in Example \ref{example2} in
Section \ref{moduli}. 

The following proposition states results contained in Corollary 3.47 and
Proposition 3.48 in \cite{B memoirs}.  Corollary 3.47 and
Proposition 3.48 are in turn based on Proposition 3.31 of \cite{B
memoirs}.  However, there is a typo on line five of the proof of
Proposition 3.31 (which appears on p.76 of \cite{B memoirs}) and hence in
the statement of Proposition 3.31 and in the statements of Corollary 3.47
and Proposition 3.48 which are based on Proposition 3.31.  Line five of the
proof of Proposition 3.31 on p.76 in \cite{B memoirs} should read 
\[\exp \left( -(\tilde{x} - x
-\tilde{\varphi}\varphi)\frac{\partial}{\partial w}  - (\tilde{\varphi} -
\varphi) \Bigl(\frac{\partial}{\partial \rho} - \rho
\frac{\partial}{\partial w} \Bigr) \right) \cdot \]
rather than
\[\exp \left( -(\tilde{x} - x)\frac{\partial}{\partial w}  -
(\tilde{\varphi} - \varphi) \Bigl(\frac{\partial}{\partial \rho} - \rho
\frac{\partial}{\partial w} \Bigr) \right) \cdot .\]
We will need the corrected results {}from Corollary 3.47 and
Proposition 3.48 in \cite{B memoirs}, and these results are given in the
proposition below. It is the first line of the righthand side of equation
(\ref{last second Theta equation}) below that contains the corrected term.

\begin{prop} (\cite{B memoirs})\label{second Theta identity in End}
The formal series $\Theta^{(2)}_j (\{t^k \mathcal{B}_k, t^{k - 1/2}
\mathcal{N}_{k - 1/2} \}_{k \in \Z}, (x,\varphi))$, for $j \in
\frac{1}{2} \mathbb{N}$, are in $\mathbb{C} [x^{-1}, \varphi]
[\mathcal{B}][\mathcal{N}][[t^{1/2}]]$.  Thus for $(B,N) \in 
\bigwedge_*^\infty$, the series 
\[\Theta^{(2)}_j (\{t^k B_k, t^{k - \frac{1}{2}} N_{k - \frac{1}{2}} \}_{k
\in \Z}, (x,\varphi))\] 
are well defined and belong to $\bigwedge_* [x^{-1}, \varphi] [[t^{1/2}]]$. 

Let $(\tilde{x}(t^{1/2}),\tilde{\varphi} (t^{1/2})) = (H_{ \{t^j B_j, 
t^{j - 1/2} N_{j - 1/2} \}_{j \in \Z}}^{(2)})^{-1} \circ I (x,\varphi)$, 
and let $V$ be a positive-energy module for $\mathfrak{ns}$.  Then in 
$((\mbox{\em End} \; (\bigwedge_* \otimes_\mathbb{C} V))
[[t^{1/2}]][[x,x^{-1}]][\varphi])^0$, the  following identity holds for
$\Theta^{(2)}_j (t^{1/2}) = \Theta^{(2)}_j (\{t^k  B_k, t^{k - 1/2} N_{k -
1/2} \}_{k \in \Z}, (x,\varphi))$. 
\begin{multline}\label{last second Theta equation}
\exp  \biggl( \sum_{m = -1}^{\infty} \sum_{j \in \Z} \binom{-j + 1}{m +
1} x^{-j - m} \biggl( \Bigl( t^j B_j +  2\varphi t^{j - \frac{1}{2}} N_{j
- \frac{1}{2}} \Bigr) L(m) \\
+ \; \Bigl( t^{j - \frac{1}{2}}  N_{j - \frac{1}{2}} + \varphi x^{-1}
\frac{(-j-m)}{2} t^j B_j \Bigr)  G(m + \frac{1}{2}) \biggr) 
\biggr)
\end{multline} 
\begin{multline*}
= \; \exp \left( (\tilde{x}(t^\frac{1}{2}) - x - \tilde{\varphi}
(t^\frac{1}{2}) \varphi) L(-1) + (\tilde{\varphi}(t^\frac{1}{2}) - \varphi)
G(-\frac{1}{2}) \right) \cdot \\  
\exp \Biggl( \! - \! \sum_{j \in \Z} \Bigl( \Theta^{(2)}_j (t^\frac{1}{2})
L(j) + \Theta^{(2)}_{j - \frac{1}{2}} (t^\frac{1}{2}) G(j - \frac{1}{2})
\Bigr) \! \Biggr) \! \cdot \exp \left( - 2\Theta^{(2)}_0 (t^\frac{1}{2})
L(0) \right) . 
\end{multline*} 
\end{prop}

\section[$N=1$ NS-VOSAs]{$N=1$ Neveu-Schwarz vertex operator
superalgebras}\label{NS-VOSA-section}

In this section, we recall the notion of {\it $N = 1$ Neveu-Schwarz 
vertex operator superalgebra over a Grassmann algebra and with odd 
formal variables ($N=1$ NS-VOSA)} given in \cite{B thesis} and \cite{B
vosas} and  recall some of the consequences of this notion which we
will need later.    

Let $x$, $x_0$, $x_1$ and $x_2$ be even formal variables, and
let $\varphi$, $\varphi_1$ and $\varphi_2$ be odd
formal variables.  For any formal Laurent series $f(x) \in  
\bigwedge_\infty[[x,x^{-1}]]$, we can define 
\begin{equation}\label{expansion}
f(x + \varphi_1 \varphi_2) = f(x) + \varphi_1
\varphi_2 f'(x) \; \in \mbox{$\bigwedge_\infty$} [[x,x^{-1}]] [\varphi_1]
[\varphi_2]. 
\end{equation}

Recall (cf. \cite{FLM}) the {\it formal $\delta$-function at $x=1$} given
by $\delta(x) = \sum_{n \in \mathbb{Z}} x^n$. 
As developed in \cite{B vosas}, we have the following $\delta$-function of 
expressions involving three even formal variables and two odd formal 
variables
\begin{eqnarray*}
\delta \biggl( \frac{x_1 - x_2 - \varphi_1 \varphi_2}{x_0} \biggr) &=&
\sum_{n \in \mathbb{Z}} (x_1 - x_2 - \varphi_1 \varphi_2)^n x_0^{-n} \\
&=& \delta \biggl( \frac{x_1 - x_2}{x_0} \biggr)  -
\varphi_1 \varphi_2 x_0^{-1} \delta' \biggl( \frac{x_1 - x_2}{x_0}
\biggr)
\end{eqnarray*} 
where $\delta'(x) = d/dx \; \delta (x) = \sum_{n \in \mathbb{Z}} n
x^{n-1}$, and we use the conventions that a function of even and odd
variables should be expanded about the even variables and any expression in
two even variables (such as $(x_1 - x_2)^n$, for $n \in \mathbb{Z}$) should
be expanded in positive powers of the second variable, (in this case $x_2$). 
        
{}From \cite{B vosas}, we have the following $\delta$-function
identity which will be used in the proofs of Lemmas \ref{bracket-Lemma} and 
\ref{bracket-Lemma2}.
\begin{equation}\label{delta 2 terms with phis}
x_1^{-1} \delta \biggl( \frac{x_2 + x_0 + \varphi_1 \varphi_2}{x_1}
\biggr) = x_2^{-1} \delta \biggl( \frac{x_1 - x_0 - \varphi_1
\varphi_2}{x_2} \biggr)  .
\end{equation}

\begin{defn}\label{VOSA definition}
{\em An} $N = 1$ NS-VOSA over $\bigwedge_*$ and with odd variables
{\em is a $\frac{1}{2} \mathbb{Z}$-graded (by weight)
$\bigwedge_\infty$-module which is  also $\mathbb{Z}_2$-graded (by sign)  
\begin{equation}\label{vosa1}
V = \coprod_{k \in \frac{1}{2} \mathbb{Z}} V_{(k)} = \coprod_{k
\in \frac{1}{2}\mathbb{Z}} V_{(k)}^0 \oplus \coprod_{k \in
\frac{1}{2}\mathbb{Z}} V_{(k)}^1 = V^0 \oplus V^1 
\end{equation}  
such that only $\bigwedge_* \subseteq \bigwedge_\infty$ acts
nontrivially, 
\begin{equation}\label{vosa2}
\dim V_{(k)} < \infty \quad \mbox{for} \quad k \in \frac{1}{2}
\mathbb{Z} , 
\end{equation}
\begin{equation}\label{positive energy}
V_{(k)} = 0 \quad \mbox{for $k$ sufficiently small} , 
\end{equation}
equipped with a linear map $V \otimes V \longrightarrow V[[x,x^{-1}]]
[\varphi]$, or equivalently,
\begin{eqnarray*} 
V &\longrightarrow&  (\mbox{End} \; V)[[x,x^{-1}]][\varphi] \\
v  &\mapsto&  Y(v,(x,\varphi)) = \sum_{n \in \mathbb{Z}} v_n x^{-n-1} +
\varphi \sum_{n \in \mathbb{Z}} v_{n - \frac{1}{2}} x^{-n-1}
\end{eqnarray*}
where $v_n \in (\mbox{End} \; V)^{\eta(v)}$ and $v_{n - 1/2} \in 
(\mbox{End} \; V)^{(\eta(v) + 1) \mbox{\begin{footnotesize} mod 
\end{footnotesize}} 2}$ for $v$ of homogeneous sign in $V$, $x$ is an 
even formal variable, and $\varphi$ is an odd formal variable, and 
where $Y(v,(x,\varphi))$ denotes the} vertex operator associated with 
$v$, {\em and equipped also with two distinguished homogeneous vectors 
$\mathbf{1} \in V_{(0)}^0$ (the {\em vacuum}) and $\tau \in V_{(3/2)}^1$ 
(the {\em Neveu-Schwarz element}).  The following conditions are assumed 
for $u,v \in V$:   
\begin{equation}\label{truncation}
u_n v = 0 \quad \mbox{for $n \in \frac{1}{2} \mathbb{Z}$ sufficiently
large;} 
\end{equation}
\begin{equation}\label{vacuum identity}
Y(\mathbf{1}, (x, \varphi)) = 1 \quad \mbox{(1 on the right being
the identity operator);} 
\end{equation}
the} creation property {\em holds:
\[Y(v,(x,\varphi)) \mathbf{1} \in V[[x]][\varphi] \qquad \mbox{and}
\qquad \lim_{(x,\varphi) \rightarrow 0} Y(v,(x,\varphi)) \mathbf{1} =
v ; \] 
the} Jacobi identity {\em holds:  
\begin{multline*}
x_0^{-1} \delta \biggl( \frac{x_1 - x_2 - \varphi_1 \varphi_2}{x_0}
\biggr) Y(u,(x_1, \varphi_1))Y(v,(x_2, \varphi_2)) \\
- (-1)^{\eta(u)\eta(v)} x_0^{-1} \delta \biggl( \frac{x_2 - x_1 + 
\varphi_1 \varphi_2}{-x_0} \biggr)Y(v,(x_2, \varphi_2))Y(u,(x_1,
\varphi_1)) \\
= \; x_2^{-1} \delta \biggl( \frac{x_1 - x_0 - \varphi_1
\varphi_2}{x_2} \biggr) Y(Y(u,(x_0, \varphi_1 - \varphi_2))v,(x_2,
\varphi_2)) , 
\end{multline*}
for $u,v$ of homogeneous sign in $V$; the $N=1$ Neveu-Schwarz algebra 
relations hold:
\begin{eqnarray*}
\left[L(m),L(n) \right] \! &=& \!(m - n)L(m + n) + \frac{1}{12} (m^3 - m)
\delta_{m + n , 0} (\mbox{rank} \; V) , \\ \label{V9}
\biggl[ G(m + \frac{1}{2}),L(n) \biggr] \! &=& \! (m - \frac{n - 1}{2} ) G(m
+ n + \frac{1}{2}) ,\\ \label{V10}
\biggl[ G (m + \frac{1}{2} ) , G(n - \frac{1}{2} ) \biggr] \! &=& \! 2L(m +
n) + \frac{1}{3} (m^2 + m) \delta_{m + n , 0} (\mbox{rank} \; V) , \label{V11}
\end{eqnarray*}
for $m,n \in \mathbb{Z}$, where 
\[G(n + \frac{1}{2}) = \tau_{n + 1}, \qquad \mbox{and} \qquad 2L(n) =
\tau_{n + \frac{1}{2}} \qquad \mbox{for} \; n \in \mathbb{Z} , \]
i.e., 
\begin{equation}\label{stress tensor}
Y(\tau,(x,\varphi)) = \sum_{n \in \mathbb{Z}} G (n + \frac{1}{2}) x^{-
n - \frac{1}{2} - \frac{3}{2}} \; + \; 2 \varphi \sum_{n \in \mathbb{Z}}
L(n) x^{- n - 2} ,
\end{equation}
and $\mbox{rank} \; V \in \mathbb{C}$; 
\begin{equation}\label{grading for vosa with}
L(0)v = kv \quad \mbox{for} \quad k \in \frac{1}{2} \mathbb{Z} \quad
\mbox{and} \quad v \in V_{(k)}; 
\end{equation}
and the} $G(-1/2)$-derivative property {\em holds:
\begin{equation}\label{G(-1/2)-derivative}
\biggl( \frac{\partial}{\partial \varphi} + \varphi
\frac{\partial}{\partial x} \biggr) Y(v,(x,\varphi)) =  Y(G(- 
\frac{1}{2})v,(x,\varphi)) . 
\end{equation} }
\end{defn}

The $N=1$ NS-VOSA just defined is denoted
by $(V,Y(\cdot,(x,\varphi)),\mathbf{1},\tau)$, or for simplicity by $V$.

The following are consequences of the definition of $N=1$ NS-VOSA with
odd formal variables which we will need in our proofs of the change of
variables formulas (cf. \cite{B vosas}).  We have 
\[ L(n) \mathbf{1} = G(n + \frac{1}{2}) \mathbf{1} = 0, \quad \mbox{for
$n \geq -1$}\]
and $\tau = G(-3/2) \mathbf{1}$. There exists $\omega = (1/2)G(-1/2) \tau
\in V_{(2)}$ such that 
\begin{equation}
Y(\omega,(x,\varphi)) = \sum_{n \in \mathbb{Z}} L(n) x^{-n-2} -
\frac{\varphi}{2} \sum_{n \in \mathbb{Z}} (n + 1) G(n -
\frac{1}{2}) x^{-n-2} ,
\end{equation}
and $\omega = L(-2) \mathbf{1}$.  The supercommutator formula is given by 
\begin{multline}\label{bracket relation for a vosa}
[ Y(u, (x_1,\varphi_1)), Y(v,(x_2,\varphi_2))] \\
= \mbox{Res}_{x_0} x_2^{-1} \delta \biggl( \frac{x_1
- x_0 - \varphi_1 \varphi_2}{x_2} \biggr) Y(Y(u,(x_0, \varphi_1 -
\varphi_2))v,(x_2, \varphi_2))  
\end{multline}
where $\mbox{Res}_{x_0}$ of a power series in $x_0$ is the coefficient of 
$x_0^{-1}$.  We have
\begin{eqnarray}
x_0^{2L(0)} Y(v, (x, \varphi)) x_0^{-2 L(0)} = Y(x_0^{2L(0)}v, (x_0^2 x,
x_0 \varphi)), \label{conjugate by L(0)}
\end{eqnarray}
and
\begin{eqnarray}
Y(e^{x_0 L(-1) + \varphi_0 G(-\frac{1}{2})}v, (x,\varphi)) &=& e^{x_0
\frac{\partial}{\partial x} + \varphi_0 \left(
\frac{\partial}{\partial \varphi} + \varphi \frac{\partial}{\partial 
x} \right)} Y(v,(x,\varphi)) \label{exponential L(-1) and G(-1/2)
property} \\  
&=& Y(v,(x + x_0 + \varphi_0 \varphi, \varphi_0 + \varphi)) . \nonumber
\end{eqnarray}
and
\begin{eqnarray}
& & \hspace{-2.2in} e^{x_0 L(-1) + \varphi_0 G(-\frac{1}{2})}
Y(v,(x,\varphi)) e^{- x_0 L(-1) - \varphi_0 G(-\frac{1}{2})} =  \nonumber\\
\hspace{1.9in} &=& Y(e^{x_0 L(-1) + \varphi_0
G(-\frac{1}{2}) - 2\varphi_0 \varphi L(-1)}v,(x,\varphi))  \nonumber \\
&=& Y(v,(x +x_0  + \varphi \varphi_0, \varphi + \varphi_0)). 
\label{L(-1) and G(-1/2) exp 1}
\end{eqnarray}

Let $(V_1, Y_1(\cdot,(x,\varphi)),\mathbf{1}_1,\tau_1)$ and $(V_2,
Y_2(\cdot,(x,\varphi)),\mathbf{1}_2,\tau_2)$ be two $N=1$ NS-VOSAs over
$\bigwedge_{*_1}$ and $\bigwedge_{*_2}$, respectively.  A {\it
homomorphism of $N=1$ NS-VOSAs} is a doubly graded
$\bigwedge_\infty$-module homomorphism $\gamma : V_1 \longrightarrow V_2
\;$ such that 
\[\gamma (Y_1(u,(x,\varphi))v) = Y_2(\gamma(u),(x,\varphi))\gamma(v)
\quad \mbox{for} \quad u,v \in V_1 ,\]
$\gamma(\mathbf{1}_1) = \mathbf{1}_2$, and $\gamma(\tau_1) =
\tau_2$.

\begin{rema}\label{homo remark}{\em If $(V_1, Y_1 (\cdot, (x,
\varphi)), \mathbf{1}_1, \tau_1)$ is an $N=1$ NS-VOSA and $V_2$ is a
$\bigwedge_\infty$-module which is isomorphic to $V_1$ for some
$\bigwedge_\infty$-module isomorphism $\gamma : V_1 \longrightarrow V_2$,
then defining $Y_2 : V_2 \longrightarrow (\mathrm{End} V_2) [[x, x^{-1}]]
[\varphi]$ by
\[Y_2 (u,(x,\varphi))v = \gamma(Y_1(\gamma^{-1}(u), (x,\varphi))
\gamma^{-1}(v)) \] 
for all $u,v \in V_2$, we have that $(V_2,Y_2(\cdot, (x,\varphi)),
\gamma(\mathbf{1}_1), \gamma(\tau_1))$, with $\frac{1}{2}
\mathbb{Z}$-grading and $\mathbb{Z}_2$-grading on $V_2$ induced by $V_1$, is
an $N=1$ NS-VOSA which is isomorphic to $(V_1,Y_1(\cdot,(x,\varphi)),
\mathbf{1}_1, \tau_1)$.}
\end{rema}

\begin{rema}{\em
Let $(V,Y(\cdot,(x,\varphi)), \mathbf{1}, \tau)$ be an $N=1$ NS-VOSA.  In
\cite{B vosas}, we also study the notion of $N=1$ NS-VOSA over a
Grassmann algebra without odd formal variables and show that
$(V,Y(\cdot,(x,0)), \mathbf{1}, \tau)$ is such an algebra.  Conversely,
given an $N=1$ NS-VOSA without odd formal variables $(V,Y(\cdot,x),
\mathbf{1}, \tau)$,  we can define $\tilde{Y}(v, (x,\varphi)) = Y(v,x) +
\varphi Y(G(-1/2) v,x)$,  and then $(V,\tilde{Y}(\cdot,(x,\varphi)),
\mathbf{1}, \tau)$ is an $N=1$ NS-VOSA with odd formal variables.  Using
this correspondence, in \cite{B vosas} we prove that the category of
$N=1$ NS-VOSAS over $\bigwedge_*$ with odd formal variables is isomorphic
to the category of $N=1$ NS-VOSAs over $\bigwedge_*$ without odd formal 
variables.   However, in including the odd formal variables the
correspondence with the geometry and the role of the operator $G(-1/2)$, as
in the $G(-1/2)$-derivative property, is more  explicit.  In addition, in a
related issue, the choice of the odd components, i.e., the $\varphi$
components of the vertex operators, determines the form of the
$G(-1/2)$-derivative property which in turn is directly related to the
choice of superconformal operator $D$ for the underlying worldsheet
geometry. }
\end{rema}

\section{Superconformal change of variables formulas for $N=1$
NS-VOSAs}\label{change-section}

Let $(V,Y(\cdot, (x,\varphi)), \mathbf{1}, \tau)$ be an $N=1$ NS-VOSA.
In this section, we present the change of variables formulas for two
superconformal changes of variables. First, in Section \ref{change-zero},
we consider the change of variables $(x,\varphi) \mapsto H(x,\varphi)$ where
$H$ is a formal superconformal function vanishing at zero with invertible
leading even coefficient of $\varphi$ in its power series expansion about
zero.  This last condition is equivalent to $H$ being bijective in a
neighborhood of zero if $H$ is convergent in a neighborhood of zero,
however we do not assume $H$ is convergent in a neighborhood of zero.

Second, in Section \ref{change-infinity}, we consider the change of
variables $(x,\varphi) \mapsto H \circ I (x,\varphi)$ where $H^{-1}$ is a
formal  superconformal function vanishing at infinity with leading
even coefficient of $\varphi x^{-1}$ equal to $i$ in its power series
expansion about infinity.  Thus $H \circ I$ is a formal invertible
superconformal function taking infinity to infinity.  

The results of this section will be proved in Sections \ref{geom} and
\ref{alg-section}.  In Section \ref{geom}, we give geometric proofs of the
change of variables formulas for changes of variables convergent in a
neighborhood of zero and infinity, respectively, and in Section
\ref{alg-section}, we give purely algebraic proofs in the more
general case of formal changes of variables not necessarily convergent.

\subsection{Superconformal change of variables at zero}\label{change-zero}

Let $(V,Y(\cdot, (x,\varphi)), \mathbf{1}, \tau)$ be an $N=1$ NS-VOSA and
let $H$ be a formal invertible superconformal change of variables
vanishing at zero.  By Proposition \ref{above}, $H^{-1}$ can be written
uniquely as 
\begin{equation}\label{inverse coordinate change}
H^{-1} (x,\varphi) =   \exp \biggl(- \! \! \sum_{j \in \mathbb{Z}_+}
\Bigl(A_j L_j(x,\varphi) + M_{j - \frac{1}{2}}  G_{j-\frac{1}{2}}(x,\varphi)
\Bigr)\biggr)  \cdot \asqrt^{-2L_0(x,\varphi)} \cdot (x,\varphi) 
\end{equation}
for some $\asqrt \in (\bigwedge_\infty^0)^\times$ and $(A_j, M_{j- 1/2}) \in
\bigwedge_\infty$ for $j \in \Z$, and by Proposition 3.21 in \cite{B
memoirs}
\begin{equation}\label{coordinate change}
H (x,\varphi) =  \asqrt^{2L_0(x,\varphi)} \cdot  \exp \biggl(\sum_{j \in
\mathbb{Z}_+} \Bigl(A_j L_j(x,\varphi) + M_{j - \frac{1}{2}} 
G_{j-\frac{1}{2}}(x,\varphi)
\Bigr)\biggr)  \cdot (x,\varphi) .
\end{equation}

Let $t^{1/2}$ be an even formal variable and let $H_{t^{1/2}}(x,\varphi) \in
\bigwedge_\infty [t^{1/2}] [[x]] [\varphi]$ be defined by
\[H_{t^{1/2}}(x,\varphi) = (t^{-\frac{1}{2}} \asqrt)^{2L_0(x,\varphi)} \cdot
\exp \biggl(\sum_{j \in \mathbb{Z}_+} \Bigl(A_j L_j(x,\varphi) + M_{j -
\frac{1}{2}}  G_{j-\frac{1}{2}}(x,\varphi) \Bigr) \! \biggr) \! \cdot
(x,\varphi) \]
so that by Proposition 3.13 in \cite{B memoirs}
\begin{multline*}
H_{t^{1/2}}^{-1}(x,\varphi) = \exp \biggl(- \! \! \sum_{j \in 
\mathbb{Z}_+} \Bigl(A_j L_j(x,\varphi) + M_{j - \frac{1}{2}} 
G_{j-\frac{1}{2}}(x,\varphi) \Bigr)\biggr)  \cdot \\
(t^{-\frac{1}{2}} \asqrt)^{-2L_0(x,\varphi)} \cdot (x,\varphi) . 
\end{multline*}

Let $V$ be a positive energy module for the $N=1$ Neveu-Schwarz algebra. 
Define
\begin{eqnarray*}
\gamma_{H, t^{1/2}} : V &\longrightarrow& V[t^{-\frac{1}{2}},
t^{\frac{1}{2}}]\\ v &\mapsto&  \exp \biggl(- \! \! \sum_{j \in
\mathbb{Z}_+} \Bigl(A_j L(j) + M_{j - \frac{1}{2}} G(j-\frac{1}{2})
\Bigr) \biggr) \cdot (t^{-\frac{1}{2}} \asqrt)^{-2L(0)} \cdot  v .
\end{eqnarray*}
Define $\gamma_H : V \longrightarrow V$, by $\gamma_H =
\gamma_{H,t^{1/2}} |_{t^{1/2} = 1}$.  Note that $\gamma_H$ is well-defined
and bijective with inverse $\gamma_H^{-1} = \gamma_{H^{-1}}$.  

Define
\begin{eqnarray*}
\gamma_{\Theta(H, t^{1/2}, (x,\varphi))} : V \! \! &\longrightarrow& \! \!
V[x,\varphi] [[t^{\frac{1}{2}}]] \\  
v \! \! &\mapsto& \! \! (t^{-\frac{1}{2}} \asqrt)^{-2L(0)} \cdot \exp
\biggl( - \! \! \sum_{j \in \mathbb{Z}_+} \Bigl( \Theta^{(1)}_j
(t^{\frac{1}{2}}) L(j) \Bigr. \biggr.
\\ & & \Bigl. \biggl. \hspace{.3in} + \; \Theta^{(1)}_{j - \frac{1}{2}}
(t^{\frac{1}{2}}) G(j-\frac{1}{2}) \Bigr) \! \biggr) \cdot \exp
(-\Theta_0^{(1)}(t^{\frac{1}{2}}) 2L(0)) \cdot v
\end{eqnarray*}
where the $\Theta^{(1)}_j (t^{1/2}) = \Theta^{(1)}_j (
t^{-1/2} \asqrt, A,M, (x,\varphi)) \in
\bigwedge_\infty[x,\varphi][[t^{1/2}]]$, for $j \in \frac{1}{2}
\mathbf{Z}_+$, are defined by (\ref{define first Theta}) with
$H^{(1)}_{\alpha_0^{1/2}, \mathcal{A}, \mathcal{M}} (x,\varphi) =
H^{-1}_{t^{1/2}} (x,\varphi)$.  That is
\begin{multline}\label{defining Thetas}
\Bigl(\exp(\Theta^{(1)}_0(t^{\frac{1}{2}} \asqrt, A,M,  (x, \varphi)),
\Bigl\{\Theta^{(1)}_j(t^{\frac{1}{2}} \asqrt, A, M, (x, \varphi)),
\Bigr. \Bigr. \\
\Bigl. \Bigl. \Theta^{(1)}_{j - \frac{1}{2}} (t^{\frac{1}{2}} \asqrt, 
A, M, (x, \varphi)) \Bigr\}_{j \in \Z} \Bigr)
\end{multline}  
\[ = \;\hat{E}^{-1}(s_{(x,\varphi)} \circ H^{-1}_{t^{1/2}} \circ
s_{H_{t^{1/2}}(x, \varphi)}^{-1} (t \asqrt^{-2} w,t^{\frac{1}{2}}
\asqrt^{-1} \rho)) .  \hspace{1.5in} \]
In other words, the $\Theta_j^{(1)} (t^{1/2})$'s are uniquely determined
by
\begin{multline*}
s_{(x,\varphi)} \circ H^{-1}_{t^{1/2}} \circ s_{H_{t^{1/2}}(x,
\varphi)}^{-1} (t \asqrt^{-2} w, t^{\frac{1}{2}} \asqrt^{-1} \rho) \\  
= \exp \Biggl( \sum_{j \in  \Z} \biggl( \Theta^{(1)}_j (t^{\frac{1}{2}})
\! \Bigl( \Lw \Bigr) + \Theta^{(1)}_{j - \frac{1}{2}}
(t^{\frac{1}{2}}) w^j \Bigl( \frac{\partial}{\partial \rho} - \rho
\frac{\partial}{\partial w} \Bigr) \! \biggr) \! \! \Biggr) \! \cdot \\
\exp \left(\! \Theta^{(1)}_0 (t^{\frac{1}{2}}) \Bigl( \twoLow \Bigr) \!
\right) \! \cdot (w,\rho) . 
\end{multline*}  
Finally, define 
\begin{eqnarray*}
\gamma_{\Theta(H, t^{1/2}, (x,\varphi))}^{-1} : V \! \! &\longrightarrow& \!
\! V[x,\varphi] [[t^{\frac{1}{2}}]] \\  
v \! \! &\mapsto& \! \! \exp (\Theta_0^{(1)}(t^{\frac{1}{2}}) 2L(0)) \cdot
\exp \biggl( \sum_{j \in \mathbb{Z}_+} \Bigl( \Theta^{(1)}_j
(t^{\frac{1}{2}}) L(j) \Bigr. \biggr. \\ 
& & \Bigl. \biggl. \hspace{.7in} + \; \Theta^{(1)}_{j - \frac{1}{2}}
(t^{\frac{1}{2}}) G(j-\frac{1}{2}) \Bigr) \! \biggr) \cdot (t^{-\frac{1}{2}}
\asqrt)^{2L(0)} \cdot  v .
\end{eqnarray*}

\begin{prop}\label{change prop}
Let $(V, Y(\cdot, (x,\varphi)), \mathbf{1}, \tau)$ be an $N=1$ NS-VOSA and
$H(x,\varphi)  \in x\bigwedge_\infty[[x]] \oplus \varphi
\bigwedge_\infty [[x]]$ superconformal and invertible.  Then
$H$ can be expressed uniquely as (\ref{coordinate change}), the map
\begin{equation}
\gamma_{\Theta(H, (x,\varphi))} = \gamma_{\Theta(H,t^{1/2}, (x,\varphi))}
|_{t^{1/2} = 1} : V \longrightarrow V[[x]][\varphi]
\end{equation}
is well defined, and in $V[[x^{-1},x]] [\varphi]$, we have the following
change of variables formula 
\begin{equation}\label{change formula}
\gamma_H (Y(u,(x,\varphi))v) =  Y(\gamma_{\Theta(H,
(x,\varphi))}(u),H(x,\varphi))
\gamma_H(v) ,
\end{equation}
for $u,v \in V$.  In particular, we have that $\gamma_{\Theta(H,
(x,\varphi))}^{-1} = \gamma_{\Theta(H, t^{1/2}, (x,\varphi))}^{-1}|_{t^{1/2}
= 1}$ is well-defined, and 
\begin{equation}\label{change formula rewritten}
Y(u,H(x,\varphi))v = \gamma_H (Y(\gamma_{\Theta(H, (x,\varphi))}^{-1}
(u),(x,\varphi)) \gamma_H^{-1} (v)) .
\end{equation}
Furthermore, if $H$ is convergent in a DeWitt open neighborhood of zero,
then for $(z,\theta)$ in the domain of convergence, both sides of
(\ref{change formula}) and of (\ref{change formula rewritten}) exist for
$(x,\varphi) = (z,\theta)$ and are equal.
\end{prop}

\begin{rema} {\em The odd components $M_{j - 1/2}$, for $j \in \Z$, and the
even soul components $(\asqrt)_S$ and $(A_j)_S$, for $j \in \Z$,  of the
change of variables related to $H$ act on both the even and odd components
of the vertex operator and thus appear even if we restrict ourselves to an
$N=1$ NS-VOSA without odd formal variable components.  However, if $V$ is an
$N=1$ NS-VOSA over $\mathbb{C}$ rather than over a Grassmann algebra
$\bigwedge_* \neq \mathbb{C}$ then the affect of the soul components
of the change of variables related to $H$, i.e., the affects of
$(\asqrt)_S$, $(A_j)_S$ and $M_{j - 1/2}$, for $j \in \Z$, will be trivial.
(Recall that an $N=1$ NS-VOSA over a Grassmann algebra
$\bigwedge_*$ is defined to be a $\bigwedge_\infty$-module where only the
subalgebra $\bigwedge_*$ acts nontrivially.)  This trivialization of the
soul components for $N=1$ NS-VOSAs over $\mathbb{C}$ occurs regardless
of whether the odd variable components of $V$ are included or not.  Thus,
if $V$ is an $N=1$ NS-VOSA over $\mathbb{C}$, to realize a superconformal
change of coordinates involving soul components, one must consider
$\bigwedge_* \otimes_\mathbb{C} V$ with $\bigwedge_* = \bigwedge_\infty$
or $\bigwedge_L$, where $L>0$. }
\end{rema}

\begin{rema}\label{non-zero-remark} 
{\em For invertible superconformal change of variables $F(x,\varphi)$ 
vanishing at $(z,\theta) \neq 0$, we can write $F(x,\varphi) = H \circ
s_{(z,\theta)} (x,\varphi)$ where $H$ is invertible and vanishing at zero,
and we can use the superconformal shift formula given by (\ref{L(-1) and
G(-1/2) exp 1}) which, when combined with the change of variables formulas
given in Proposition \ref{change prop} above, provide the appropriate
change of variables formulas. }
\end{rema}

\begin{rema} {\em Replacing the superconformal function $H(x,\varphi)$ by
the body portion of $H$, setting all odd variables and the soul portion
of supernumbers equal to zero and restricting $V$ to $V^0$ in Proposition
\ref{change prop}, we obtain the change of variables formula for a VOA in
the case of an invertible analytic function vanishing at zero as developed
in \cite{H book}.}
\end{rema} 

\subsection{Superconformal change of variables at
infinity}\label{change-infinity}
Let  $(V,Y(\cdot, (x,\varphi)), \mathbf{1},\tau)$ be an $N=1$ NS-VOSA.
In this section we consider a change of variables $(x,\varphi) \mapsto
H \circ I(x,\varphi)$ where $H^{-1}$ is formally superconformal vanishing
at infinity with even coefficient of $\varphi x^{-1}$ equal to $i$ in
its power series expansion about infinity.  By Remark \ref{corollary to
above}, $H^{-1}$ can be written uniquely as  
\begin{eqnarray}\label{coordinate change22}
\hspace{.3in} \; H^{-1} (x,\varphi) = \exp \Biggl( \sum_{j \in
\mathbb{Z}_+} \Biggl(B_j  L_{-j}(x,\varphi) + N_{j - \frac{1}{2}}
G_{-j+\frac{1}{2}}(x,\varphi) \Biggr)\Biggr) \cdot
\Bigl(\frac{1}{x},\frac{i\varphi}{x}\Bigr) ,
\end{eqnarray} 
for some $(B,N) \in \bigwedge_\infty^\infty$.  Let 
\[H^{-1}_{t^{1/2}}(x,\varphi) = \exp \Biggl( \sum_{j \in \mathbb{Z}_+}
\Biggl(t^j B_j  L_{-j}(x,\varphi) + t^{j - \frac{1}{2}} N_{j - \frac{1}{2}}
G_{-j+\frac{1}{2}}(x,\varphi) \Biggr)\Biggr) \cdot
\Bigl(\frac{1}{x},\frac{i\varphi}{x}\Bigr) . \]
Proposition 3.17 in \cite{B memoirs} implies that $I^{-1} \circ
H^{-1}_{t^{1/2}} (x,\varphi)  \in (x,\varphi) + \bigwedge_\infty [t^{1/2}]
[[x^{-1}]] [\varphi]$ is given by
\begin{equation}\label{coordinate change3}
I^{-1} \circ H^{-1}_{t^{1/2}}(x,\varphi) = \exp \Biggl(\sum_{j \in
\mathbb{Z}_+} 
\Biggl(t^j B_j L_{-j}(x,\varphi) + t^{j - \frac{1}{2}} N_{j - \frac{1}{2}}
G_{-j+\frac{1}{2}}(x,\varphi) \Biggr)\Biggr) \cdot (x,\varphi),
\end{equation}
and thus 
\begin{equation}\label{exp for H circ I}
 H_{t^{1/2}} \circ I (x,\varphi) = \exp \Biggl(- \! \! \sum_{j \in
\mathbb{Z}_+} \Biggl(t^j B_j L_{-j}(x,\varphi) + t^{j - \frac{1}{2}} N_{j -
\frac{1}{2}} G_{-j+\frac{1}{2}}(x,\varphi) \Biggr)\Biggr) \cdot (x,\varphi).
\end{equation}

Let $(V, Y(\cdot, (x,\varphi)), \mathbf{1}, \tau)$ be an $N=1$ NS-VOSA.  Let
$V' = \coprod_{k\in \frac{1}{2} \mathbb{Z}} V_{(k)}^*$ be the graded dual
space of $V$, and let $L'(j)$ and $G'(j - 1/2)$, for $j \in \Z$, be the
adjoint operators corresponding to $L(-j)$ and $G(-j + 1/2)$, respectively. 
Define 
\begin{eqnarray*}
\xi_{H \circ I, t^{1/2}} : V' &\longrightarrow& V'[t^\frac{1}{2}]\\
v' &\mapsto& \exp \Biggl( - \! \! \sum_{j \in \mathbb{Z}_+} \Biggl(t^j B_j
L'(j) + t^{j - \frac{1}{2}} N_{j - \frac{1}{2}} G'(j-\frac{1}{2}) \Biggr)
\Biggr) \cdot v'
\end{eqnarray*}
and define $\xi_{H \circ I} : V' \longrightarrow V'$, by $\xi_{H \circ I} =
\xi_{H \circ I,t^{1/2}} |_{t^{1/2} = 1}$.  Note that $\xi_{H \circ I}$ is
well defined and bijective.   Let $\bar{V} = \prod_{k\in \frac{1}{2}
\mathbb{Z}} V_{(k)} = (V')^*$ be the algebraic completion of $V$. Let
\begin{eqnarray*}
\xi_{H \circ I, t^{1/2}}^* : \bar{V} &\longrightarrow& \bar{V}
[t^\frac{1}{2}] \\ 
v &\mapsto& \exp \Biggl( - \! \! \sum_{j \in \mathbb{Z}_+} \Biggl(t^j B_j
L(-j) + t^{j - \frac{1}{2}} N_{j - \frac{1}{2}} G(-j+\frac{1}{2})
\Biggr) \Biggr) \cdot v
\end{eqnarray*} 
be the adjoint operator corresponding to $\xi_{H \circ I, t^{1/2}}$, and let
$\xi_{H \circ I}^* = \xi_{H \circ I, t^{1/2}}^*|_{t^{1/2} = 1}$ be the
adjoint operator corresponding to $\xi_{H \circ I}$.

Define
\begin{eqnarray*}
\xi_{\Theta(H \circ I, t^{1/2}, (x,\varphi))} : V &\longrightarrow&
V[x^{-1},\varphi] [[t^{\frac{1}{2}}]]\\ 
v &\mapsto& \exp \Biggl(  - \! \! \sum_{j \in \mathbb{Z}_+}
\Biggl(\Theta^{(2)}_j (t^{\frac{1}{2}}) L(j) + \Theta^{(2)}_{j -
\frac{1}{2}} (t^{\frac{1}{2}}) G(j-\frac{1}{2})  \Biggr) \Biggr) \cdot \\
& & \hspace{1.7in} \exp \left(- \Theta^{(2)}_0(t^{\frac{1}{2}}) 2L(0)
\right) \cdot v
\end{eqnarray*}
where the $\Theta^{(2)}_j(t^{1/2}) = \Theta^{(2)}_j (\{t^k B_k,
t^{k- 1/2} N_{k - 1/2} \}_{k \in \Z}, (x,\varphi)) \in
\bigwedge_\infty[x^{-1}, \varphi][[t^{1/2}]]$ are defined by (\ref{define
second Theta}) with $H^{(2)}_{\mathcal{B}, \mathcal{N}} (x, \varphi) =
H_{t^{1/2}}^{-1}(x,\varphi)$. That is  
\begin{multline}\label{defining Thetas2}
\Bigl(\exp(\Theta^{(2)}_0(\{t^k B_k,t^{k- 1/2} N_{k - 1/2} \}_{k \in \Z} , 
(x, \varphi)), \\
\Bigl\{ \Theta^{(2)}_j(\{t^k B_k, t^{k- 1/2} N_{k - 1/2} \}_{k \in \Z}, (x,
\varphi)),  \\
\Theta^{(2)}_{j - \frac{1}{2}} (\{t^k B_k,
t^{k- 1/2} N_{k - 1/2} \}_{k \in \Z}, (x, \varphi)) \Bigr\}_{j \in \Z} 
\Bigr) 
\end{multline}
\begin{eqnarray*}
&=& \; \hat{E}^{-1}(s_{(x,\varphi)} \circ I^{-1} \circ
H_{t^{1/2}}^{-1} \circ s_{H_{t^{1/2}} \circ I(x,\varphi)}^{-1}
(w,\rho)) . \hspace{1.6in}
\end{eqnarray*}
In other words, the $\Theta_j^{(2)} (t^{1/2})$'s are uniquely determined
by
\begin{multline*}
s_{(x,\varphi)} \circ I^{-1} \circ
H_{t^{1/2}}^{-1} \circ s_{H_{t^{1/2}} \circ I(x,\varphi)}^{-1}
(w,\rho) \\
= \exp \Biggl( \sum_{j \in  \Z} \biggl( \Theta^{(2)}_j (t^\frac{1}{2}) \!
\Bigl( \Lw \Bigr) + \Theta^{(2)}_{j - \frac{1}{2}} (t^\frac{1}{2}) w^j
\Bigl( \frac{\partial}{\partial \rho} - \rho
\frac{\partial}{\partial w} \Bigr) \!  \biggr) \! \Biggr) \! \cdot \\
\exp \left(\! \Theta^{(2)}_0 (t^\frac{1}{2}) \Bigl( \twoLow \Bigr) \!
\right) \!  \cdot (w,\rho) . 
\end{multline*}

Finally, define
\begin{eqnarray*}
\xi_{\Theta(H \circ I, t^{1/2}, (x,\varphi))}^{-1} : V &\longrightarrow&
V[x^{-1},\varphi] [[t^{\frac{1}{2}}]]\\ 
v &\mapsto& \exp \left(\Theta^{(2)}_0(t^{\frac{1}{2}}) 2L(0)
\right) \cdot \exp \Biggl( \sum_{j \in \mathbb{Z}_+} \Biggl(\Theta^{(2)}_j
(t^{\frac{1}{2}}) L(j) \Biggr. \Biggr. \\
& & \hspace{1.5in} \Biggl. \Biggl. + \; \Theta^{(2)}_{j - \frac{1}{2}}
(t^{\frac{1}{2}}) G(j-\frac{1}{2})  \Biggr) \Biggr) \cdot  v .
\end{eqnarray*}

\begin{prop}\label{change prop2}
Let $(V, Y(\cdot, (x,\varphi)), \mathbf{1}, \tau)$ be an $N=1$
NS-VOSA and $H(x,\varphi)  \circ I \in (x,\varphi) +
\bigwedge_\infty[[x^{-1}]][\varphi]$ formally superconformal.  Then $H \circ
I$ can be expressed uniquely as (\ref{exp for H circ I}) with $t^{1/2} =
1$,  the map
\begin{equation}
\xi_{\Theta(H \circ I,(x,\varphi))} = \xi_{\Theta(H \circ I, t^{1/2},
(x,\varphi))}|_{t^{1/2} = 1} : V \longrightarrow V[[x^{-1}]][\varphi]
\end{equation}
is well defined, and in $\bar{V} [[x^{-1},x]] [\varphi]$, we have the
following change of variables formula 
\begin{equation}\label{change formula2}
Y(u,(x,\varphi))\xi_{H \circ I}^*(v) = \xi_{H \circ I}^*(Y(\xi_{\Theta(H
\circ I, (x,\varphi))} (u), H \circ I (x,\varphi))  v),
\end{equation}
for $u,v \in V$.  In particular, we have that $\xi_{\Theta(H \circ I,
(x,\varphi))}^{-1} = \xi_{\Theta(H \circ I, t^{1/2}, (x,\varphi))}^{-1}
|_{t^{1/2} = 1}$ is well-defined, and 
\begin{equation}\label{change formula2 rewritten}
Y(u, H \circ I (x,\varphi))  v  = (\xi_{H \circ I}^*)^{-1}(Y( \xi_{\Theta(H
\circ I, (x,\varphi))}^{-1} (u), (x,\varphi))\xi_{H \circ I}^*(v)).
\end{equation}
Furthermore, if $H \circ I$ is convergent in a DeWitt open neighborhood of
infinity, then for  $(z,\theta)$ in the domain of convergence of $H \circ
I$ both sides of (\ref{change formula2}) and of (\ref{change formula2
rewritten}) exist when $(x,\varphi) = (z,\theta)$ and are equal.
\end{prop}

\begin{rema} {\em Equation (\ref{change formula2}) in Proposition
\ref{change prop2} above, is the superextension of the change
of variables formula given at the top of p. 181 in \cite{H book}.  This can
be seen by letting Huang's $f^{-1}(x)$ be the body of $H(x,\varphi)$ (and
thus $f^{-1} (1/x)$ is the body of $H \circ I (x,\varphi)$), setting all odd
variables and the soul portion of supernumbers equal to zero and
restricting $V$ to $V^0$. (Note that Huang's operator $\xi^*_{1/f}$ is the
body portion of our operator $\xi^*_{H\circ I}$ if $f^{-1}$ is the body of
$H$.) }
\end{rema} 

\section{Geometric proofs for convergent superconformal changes
of variables}\label{geom}

In this section, we give geometric proofs of Propositions \ref{change prop}
and \ref{change prop2} in the case of convergent superconformal changes of
variables.  These geometric proofs rely on the correspondence between the
worldsheet supergeometry of propagating superstrings in $N=1$
superconformal field theory and the algebraic structure of $N=1$ NS-VOSA
developed in \cite{B memoirs} and \cite{B iso-thm}.  More than just giving
alternate proofs to the formal algebraic proofs we will present in Section
\ref{alg-section}, this section gives the motivation behind the change of
variables formulas.

In Sections \ref{moduli}-\ref{iso} we recall some of the basic results and
structures of this geometric/algebraic correspondence.  For  more details,
the reader should refer to \cite{B memoirs} and \cite{B iso-thm}.  In
Sections \ref{section-geom1} and \ref{section-geom2}, we present the
actual proofs of Propositions \ref{change prop} and \ref{change prop2},
respectively.  

\subsection{The moduli space of $N=1$ superspheres with tubes and the
sewing operation}\label{moduli}

In \cite{B memoirs}, we define the moduli space $SK(n)$ of $N=1$
superspheres with $n$ incoming ordered tubes and one outgoing tube, for $n
\in \mathbb{N}$.  Incoming ``tubes" represent propagating incoming
superstrings and the one outgoing ``tube" represents an outgoing
propagating superstring in the worldsheet supergeometry of $N=1$
superconformal field theory.  These propagating superstrings sweep out a
supersurface or worldsheet in space time which has superconformal coordinate
transition functions.  Such a supersurface is called a {\it super-Riemann
surface}.  It is shown in \cite{B memoirs} that each incoming or outgoing
tube can be represented by an oriented puncture and local superconformal
coordinates vanishing at the puncture.  Thus the moduli space $SK(n)$ is the
set of equivalence classes, under superconformal equivalence, of genus zero
super-Riemann surfaces with $n$ incoming ordered punctures, one outgoing
puncture, and local superconformal coordinates vanishing at the punctures.

Furthermore the Uniformization Theorem for genus-zero super-Riemann surfaces
proved by Crane and Rabin in \cite{CR} states that any compact genus-zero
super-Riemann surface (i.e., supersphere) is superconformally equivalent to
the super-Riemann sphere denoted $S\hat{\mathbb{C}}$.  Just as the usual
Riemann sphere is identified with the complex plane and a point added at
infinity, the super-Riemann sphere $S\hat{\mathbb{C}}$ can be thought of as
the Grassmann algebra $\bigwedge_\infty$ together with a fiber at infinity
given by $\infty \times (\bigwedge_\infty)_S$.  Thus $S\hat{\mathbb{C}}$
has two coordinate charts $\bigwedge_\infty = S\hat{\mathbb{C}} \diagdown
(\infty \times (\bigwedge_\infty)_S)$ and $\bigwedge_\infty \cong
S\hat{\mathbb{C}} \diagdown (0 \times (\bigwedge_\infty)_S)$ with coordinate
transition function {}from the first coordinate chart to the second given
by $I(w,\rho) = (1/w, i \rho/w)$ for $(w,\rho) \in
(\bigwedge_\infty)^\times$; see \cite{B memoirs}.

Let
\begin{multline*}
\mathcal{H}  =  \bigl\{ (A,M) \in \mbox{$\bigwedge_\infty^\infty$} 
\; | \; \tilde{E}(A,M)(z,\theta) \; \mbox{is an absolutely convergent power }
\bigr. \\
\bigl. \mbox{series in some neighborhood of } (z,\theta) = 0 \bigr\},
\end{multline*} 
and for $n \in \Z$, let
\[ SM^{n - 1} = \bigl\{ \bigl((z_1, \theta_1),...,(z_{n-1}, \theta_{n-1})\bigr) \; | \;
(z_i, \theta_i) \in \mbox{$\bigwedge_\infty^\times$}, \; (z_i)_B \neq (z_j)_B , \;
\mbox{for} \; i \neq j \bigr\} . \] 
Note that for $n=1$, the set $SM^0$ has exactly one element. It is shown in
\cite{B memoirs}, that as a set 
\[SK(n) = SM^{n-1} \times \mathcal{H} \times
\left((\mbox{$\bigwedge_\infty^0$})^\times \times \mathcal{H} \right)^n\]
for $n \in \Z$, and $SK(0) = \bigl\{(A,M) \in \mathcal{H} \; | \; (A_1,
M_{1/2}) = (0,0) \bigr\}$.  Each element $Q \in SK(n)$ can be thought of as
the super-Riemann sphere with $n$ incoming punctures, one outgoing puncture
and local superconformal coordinates vanishing at the punctures, where this
super-Riemann sphere $Q$ is a canonical representative of the equivalence
class of superspheres with tubes superconformally equivalent to $Q$.  For
example, for $Q \in SK(n)$ given by
\[Q = \bigl((z_1, \theta_1),...,(z_{n-1}, \theta_{n-1});
(A^{(0)},M^{(0)}), (\asqrt^{(1)}, A^{(1)}, M^{(1)}),...,(\asqrt^{(n)},
A^{(n)}, M^{(n)})\bigr) \]
this is the canonical representative of the equivalence class
of superspheres superconformally equivalent to $Q$ where $Q$ is
the super-Riemann sphere with outgoing puncture at $(\infty, 0) \in
S\hat{\mathbb{C}}$ and local coordinate vanishing at infinity given by
$\tilde{E}(A^{(0)}, -iM^{(0)})(1/w, i\rho/w)$, the $i$-th incoming puncture
at $(z_i, \theta_i) \in S\hat{\mathbb{C}}$ and local coordinate vanishing
at the puncture given by $\hat{E}(\asqrt^{(i)}, A^{(i)}, M^{(i)}) (w - z_i
- \rho \theta_i, \rho - \theta_i)$, for $1 \leq i \leq n-1$, and last
incoming puncture at $(0,0) \in S\hat{\mathbb{C}}$ and local coordinate
vanishing at the puncture given by $\hat{E}(\asqrt^{(n)}, A^{(n)}, M^{(n)})
(w, \rho)$.

In \cite{B memoirs}, we define a (partial) sewing operation on $SK =
\bigcup_{n \in \mathbb{N}} SK(n)$, denoted by
\begin{eqnarray*}
 _i\infty_0 : SK(m) \times SK(n) &\rightarrow& SK(m + n -  1)\\
(Q_1, Q_2) &\mapsto& Q_1 \; _i\infty_0 \; Q_2 ,
\end{eqnarray*}
where $m \in \Z$, $n \in \mathbb{N}$, and $1 \leq i \leq m$.  This sewing
together of $Q_1$ and $Q_2$ is defined if there exists  $ r_1 > r_2 >
0$ such that there exists a (DeWitt) open ball of radius $r_2 >0$ about
zero which is contained in the image of the local coordinate map vanishing
at the $i$-th puncture of $Q_1$ but whose preimage contains no other
punctures, and a (DeWitt) open ball of radius $1/r_1$ about zero which is
contained in the image of the local coordinate map vanishing at the outgoing
puncture of $Q_2$ but whose preimage contains no other punctures.  The
sewn sphere is then obtained by removing the preimage of the ball of
radius $r_2$ {}from $Q_1$ and the preimage of the ball off radius $1/r_1$ {}from $Q_2$ and identifying the superannuli of the two boundaries using the
superconformal inversion $I(w,\rho) = (1/w, i\rho/w)$.  See
\cite{B memoirs} for details. 

The resulting sewn supersphere is then in a superconformal equivalence class
of superspheres whose canonical representative can be expressed as an
element of $SK(m + n -1)$.  Much of the details which we will leave out
here, but which are given in \cite{B memoirs}, concern the determination of
the superconformal uniformizing function that maps the resulting sewn
supersphere to its canonical representative in $SK(m+n-1)$.

Below we give two examples of the sewing together of two canonical
superspheres and the determination of the resulting canonical supersphere. 
These examples are the two that we will need for the  geometric proofs of
Propositions \ref{change prop} and \ref{change prop2} in the case of
convergent changes of variables.   

\begin{ex}\label{example1}
{\em Denote by $\mathbf{0}$ the infinite sequence in
$\bigwedge_\infty^\infty$ that consists of all zeros.  Let
\[Q_1 = (\mathbf{0}, (\asqrt, A, M)) \in SK(1),\]
i.e., $Q_1$ is the canonical supersphere representative with one outgoing
tube at infinity and local coordinate given by $(1/w,i\rho/w)$, and one
incoming puncture at zero and local coordinate given by
\[H^{-1}(w,\rho) = \exp \biggl(-\sum_{j \in \mathbb{Z}_+} \Bigl(A_j
L_j(w,\rho) + M_{j - \frac{1}{2}}  G_{j-\frac{1}{2}}(w,\rho) \Bigr)\biggr) 
\cdot \asqrt^{-2L_0(w,\rho)} \cdot (w,\rho) . \]
In otherwords, the local coordinate vanishing at zero is given by
(\ref{inverse coordinate change}) with $(x,\varphi) = (w,\rho)$ where we are
assuming that $H^{-1}$ is convergent in a neighborhood of zero, i.e, that
$(A,M) \in \mathcal{H}$. 

Let
\[Q_2 = ((z,\theta); \mathbf{0}, (1,\mathbf{0}), (1,\mathbf{0})) \in SK(2),
\]
i.e., $Q_2$ is the canonical supersphere representative with one outgoing
tube at infinity and local coordinate given by $(1/w,i\rho/w)$, the first
incoming puncture at $(z,\theta)$ and local coordinate given by
$s_{(z,\theta)}(w,\rho) = (w - z - \rho \theta, \rho - \theta)$, and
second incoming puncture at zero and local coordinate $(w,\rho)$. {}from the definition of sewing on $SK$, $Q_2$ can be sewn to $Q_1$ if there
exists a ball of radius $r>0$ in the image of $H^{-1}$, i.e., in the domain
of convergence of $H$, such that $r>|z_B|$.  In this case, we can form the
sewn supersphere $Q_1 \; _1\infty_0 \; Q_2$, and there exists a
uniformizing function $F: Q_1 \; _1\infty_0 \; Q_2 \rightarrow
S\hat{\mathbb{C}}$ which, restricted to $Q_1$, is given by $F_1 (w,\rho) =
F|_{Q_1} (w,\rho) = (w,\rho)$ and, restricted to $Q_2$, is given by
$F_2(w,\rho) = F|_{Q_2} (w,\rho) = H(w,\rho)$.  

Thus the resulting canonical supersphere has: outgoing puncture at
$F_1(\infty, 0) = (\infty, 0)$ with local coordinate vanishing at the
puncture given by $I \circ F_1^{-1}(w,\rho) = I(w,\rho) = (1/w,i\rho/w)$;
first incoming puncture at $F_2(z,\theta) = H(z,\theta)$ with local
coordinate vanishing at the puncture given by $s_{(z,\theta)} \circ
F_2^{-1}(w,\rho) = s_{(z,\theta)} \circ H^{-1}(w,\rho) = s_{(z,\theta)}
\circ H^{-1}  \circ s_{H(z,\theta)}^{-1} \circ s_{H(z,\theta)} (w,\rho)$;
and last incoming puncture at $F_2(0,0) = H(0,0) = (0,0)$ with local
coordinate vanishing at the puncture given by $F_2^{-1}(w,\rho) =
H^{-1}(w,\rho)$.  

In otherwords, letting $\Theta_j^{(1)} = \Theta^{(1)}_j (t^{-1/2} \asqrt,
A,M, (z,\theta))|_{t^{1/2} = 1} \in \bigwedge_\infty$ be defined by
(\ref{defining Thetas}), the resulting canonical supersphere is given by  
\begin{eqnarray}
Q_1 \; _1\infty_0 \; Q_2 
&=& (H (z,\theta); \mathbf{0}, \hat{E}^{-1}(s_{(z,\theta)} \circ H^{-1}
\circ s_{H(z,\theta)}^{-1} (w,\rho)), (\asqrt, A, M)) \nonumber \\
&=& (H (z,\theta); \mathbf{0}, (\asqrt e^{\Theta^{(1)}_0},
\{\asqrt^{2j} \Theta^{(1)}_j, \asqrt^{2j - 1} \Theta^{(1)}_{j
- \frac{1}{2}}\}_{j\in \Z}), (\asqrt, A, M)). 
\end{eqnarray}
Note in particular that this implies that if $H(w,\rho)$ is convergent in a
neighborhood of zero, i.e., if $(A,M) \in \mathcal{H}$, and $(z,\theta)$ is
in the radius of convergence, then $\Theta^{(1)}_j (t^{-1/2} \asqrt, A,M,
(x,\varphi))$ is convergent for $t^{1/2} = 1$ and $(x,\varphi) =
(z,\theta)$.}\end{ex}

\begin{ex}\label{example2} {\em
Let 
\[Q_1 = ((z,\theta); \mathbf{0}, (1,\mathbf{0}), (1,\mathbf{0})) \in SK(2).
\]
(which is described in detail in Example \ref{example1} above) and, 
\[Q_2  = ((B,N), (1,\mathbf{0})) \in SK(1), \]
i.e., $Q_2$ is the canonical supersphere representative with one outgoing
tube at infinity and local coordinate given by 
\[H^{-1}(w,\rho) = \exp \biggl(\sum_{j \in \mathbb{Z}_+} \Bigl(B_j
L_{-j}(w,\rho) + N_{j - \frac{1}{2}}  G_{-j+\frac{1}{2}}(w,\rho)
\Bigr)\biggr)  \cdot (\frac{1}{w},\frac{i\rho}{w}) ,\]
and one incoming puncture at zero and local coordinate given by $(w,\rho)$.
In otherwords, the local coordinate vanishing at infinity is given by
(\ref{coordinate change22}) with $(x,\varphi) = (w,\rho)$, where we are
assuming that $H^{-1}$ is convergent in a neighborhood of infinity, i.e,
that $(B,N) \in \mathcal{H}$. {}from the definition of sewing on $SK$, $Q_2$ can be sewn to the last
puncture of $Q_1$ if there exists a ball of radius $1/r$, for $r>0$, about
zero which lies in the  image of $H$ and such that $1/|z_B| < 1/r$, i.e.,
such that $(z,\theta)$ is in the domain of convergence of $H \circ I$.  In
this case, we can form the sewn supersphere $Q_1 \; _2\infty_0 \; Q_2$, and
there exists a uniformizing function $F : Q_1 \; _2\infty_0 \; Q_2
\rightarrow S\hat{\mathbb{C}}$ which, restricted to $Q_1$, is given by $F_1
(w,\rho) = F|_{Q_1} (w,\rho) = H \circ I (w,\rho)$ and, restricted to
$Q_2$, is given by $F_2(w,\rho) = F|_{Q_2} (w,\rho) = (w,\rho)$. 

Thus the resulting canonical supersphere has: outgoing puncture at
$F_1(\infty, 0) =  H \circ I (\infty, 0) = (\infty, 0)$ with local
coordinate vanishing at the puncture given by $I \circ F_1^{-1}(w,\rho) =
H^{-1}(w,\rho)$; first incoming puncture at $F_1(z,\theta) =
H \circ I (z,\theta)$ with local coordinate vanishing at the puncture
given by $s_{(z,\theta)} \circ F_1^{-1}(w,\rho) = s_{(z,\theta)} \circ
I^{-1} \circ H^{-1} (w,\rho) = s_{(z,\theta)} \circ I^{-1} \circ H^{-1}
\circ s_{H \circ I (z,\theta)}^{-1} \circ s_{H \circ I (z,\theta)}
(w,\rho)$; and last incoming puncture at $F_2(0,0) = (0,0)$ with local
coordinate vanishing at the puncture given by $F_2^{-1}(w,\rho) =
(w,\rho)$. 

In otherwords, letting
$\Theta_j^{(2)} = \Theta^{(2)}_j (\{t^k B_k, t^{k- 1/2} N_{k - 1/2} \}_{k
\in \Z}, (z,\theta)) |_{t^{1/2} = 1} \in\bigwedge_\infty$
be defined by (\ref{defining Thetas2}), the resulting
canonical supersphere is given by  
\begin{eqnarray}
\lefteqn{Q_1 \; _2\infty_0 \; Q_2} \nonumber\\
&=& (H \circ I (z,\theta); (B,N) , \hat{E}^{-1}(s_{(z,\theta)}
\circ I^{-1} \circ H^{-1} \circ s_{H \circ I (z,\theta)}^{-1}
(w,\rho)), (1, \mathbf{0})) \nonumber \\ 
&=& (H \circ I (z,\theta); (B,N) , (e^{\Theta^{(2)}_0},
\{\Theta^{(2)}_j,  \Theta^{(2)}_{j -\frac{1}{2}}\}_{j\in
\Z}), (1, \mathbf{0})). 
\end{eqnarray}
Note in particular that this implies that if $H^{-1}(w,\rho)$ is convergent
in a neighborhood of infinity, i.e., if $(B,N) \in \mathcal{H}$, and
$(z,\theta)$ is in the radius of convergence of $H \circ I$, then
$\Theta^{(2)}_j (\{t^k B_k, t^{k- 1/2} N_{k - 1/2} \}_{k \in \Z},
(z,\theta))$ is convergent for $t^{1/2} = 1$ and $(x,\varphi) =
(z,\theta)$. }\end{ex}

\subsection{$N=1$ supergeometric vertex operator superalgebras}
\label{N=1 SG-VOSA-section}

In this section we introduce the notion of {\it $N=1$ supergeometric vertex
operator superalgebra ($N=1$ SG-VOSA)}. But before we do this, we will need
some further definitions and notations.

Let $\langle \cdot , \cdot \rangle$ be the natural pairing between $V'$ and
$\bar{V}$.  For $n \in \mathbb{N}$, let 
\[\mathcal{SF}_V(n) = \mbox{Hom}_{\bigwedge_\infty}(V^{\otimes n}, \bar{V})
.\] 
For $m \in \Z$, $n \in \mathbb{N}$, and any positive integer $i \leq
m$, we define the {\it $t^{1/2}$-contraction} 
\begin{eqnarray*} 
(\cdot \; _i*_0 \; \cdot)_{t^{1/2}}: \mathcal{SF}_V(m) \times \mathcal{SF}_V(n) &
\rightarrow & \mbox{Hom} (V^{\otimes (m + n - 1)},\bar{V} [[t^{\frac{1}{2}}, t^{-
\frac{1}{2}}]])  \\ 
(f,g) &\mapsto& (f \; _i*_0 \; g)_{t^{1/2}} , 
\end{eqnarray*}
by
\begin{multline}\label{t-contraction}
(f \; _i*_0 \; g)_{t^{1/2}} (v_1 \otimes \cdots \otimes v_{m + n -
1}) \\
= \sum_{k \in \frac{1}{2} \mathbb{Z}}f(v_1 \otimes \cdots \otimes
v_{i - 1} \otimes P(k)(g(v_i \otimes \cdots \otimes v_{i + n - 1}))
\otimes v_{i + n} \otimes \cdots \otimes v_{m + n - 1}) t^k 
\end{multline}
for all $v_1,...,v_{m + n - 1} \in V$, where for any $k \in
\frac{1}{2} \mathbb{Z}$, the map $P(k) : \bar{V} \rightarrow V_{(k)}$ is the
canonical projection map.

If for arbitrary $v' \in V'$, $v_1,...,v_{m + n -1} \in V$, the formal
Laurent series in $t^{1/2}$
\[ \langle v', (f \; _i*_0 \; g)_{t^{1/2}}(v_1 \otimes \cdots \otimes
v_{m + n - 1}) \rangle \]
is absolutely convergent when $t^{1/2} = 1$, then $(f \; _i*_0
\; g)_1$ is well defined as an element of $\mathcal{SF}_V(m + n -1)$, and
we define  the {\it contraction $(f \; _i*_0 \; g)$ in $\mathcal{SF}_V(m
+ n -1)$ of $f$ and $g$} by  
\[f \; _i*_0 \; g \; = \; (f \;  _i*_0 \; g)_{t^{1/2} = 1} .\]

In \cite{B iso-thm} (resp., in \cite{B memoirs}) we define a left action
of the symmetric groups $S_n$ on $\mathcal{SF}_V(n)$ (resp., on $SK(n)$). 
These actions are needed to define the notion of $N=1$ SG-VOSA but will not
be used elsewhere in this paper.  Therefore we refer the readers to \cite{B
iso-thm} and \cite{B memoirs} for details.
 
A {\it supermeromorphic superfunction on $SK(n)$}, for $n \in \Z$, is a
superfunction $F : SK(n) \rightarrow \bigwedge_\infty$ of the form
\begin{eqnarray}
\hspace{.4in} F(Q) \! \! &=& \! \! F\bigl((z_1, \theta_1),...,(z_{n-1}, 
\theta_{n-1}); (A^{(0)},M^{(0)}), (\asqrt^{(1)}, A^{(1)}, M^{(1)}), ..., 
\label{meromorphic} \bigr.\\
& & \bigl. \hspace{2.6in}(\asqrt^{(n)}, A^{(n)}, M^{(n)})\bigr) \nonumber \\
&=& \! \! F_0\bigl( (z_1,\theta_1),...,(z_{n-1}, \theta_{n-1}); (A^{(0)},M^{(0)}), 
(\asqrt^{(1)}, A^{(1)}, M^{(1)}),..., \bigr.\nonumber \\
& & \bigl. (\asqrt^{(n)},A^{(n)}, M^{(n)})\bigr) \times \biggl(\prod^{n-1}_{i
= 1} z_i^{-s_i} \prod_{1 \leq i<j \leq n-1} (z_i - z_j - \theta_i \theta_j)^{-s_{ij}}\biggr)
\nonumber 
\end{eqnarray} 
where $s_i$ and $s_{ij}$ are nonnegative integers and 
\[F_0\bigl((z_1,\theta_1),...,(z_{n-1}, \theta_{n-1}); (A^{(0)},M^{(0)}),
(\asqrt^{(1)}, A^{(1)}, M^{(1)}),...,(\asqrt^{(n)},A^{(n)}, M^{(n)})\bigr)\]  
is a polynomial in the $z_i$'s, $\theta_i$'s, $\asqrt^{(i)}$'s,
$(\asqrt^{(i)})^{-1}$'s, $A^{(i)}_j$'s, and $M^{(i)}_{j - 1/2}$'s.
For $n=0$ a {\it supermeromorphic superfunction on $SK(0)$} is a
polynomial in the components of elements of $SK(0)$, i.e., a polynomial
in the $A^{(0)}_j$'s, and $M^{(0)}_{j - 1/2}$'s.  For $F$ of the form
(\ref{meromorphic}), we say that $F$ has a pole of order $s_{ij}$ at 
$(z_i,\theta_i) = (z_j,\theta_j)$.

Finally, the last ingredient we will need to introduce the notion of $N=1$
SG-VOSA is the $\Gamma(\asqrt, A,M,B,N)$ series which was introduced and
studied in \cite{B memoirs} and \cite{B iso-thm}.  This series is based on
the local coordinates of the incoming puncture, represented by $(\asqrt,
A,M) \in (\bigwedge_\infty^0)^\times \times \mathcal{H}$, and the outgoing
puncture, represented by $(B,N)$, being sewn when two superspheres with
tubes are sewn together. The main fact we will need for the geometric
proofs of the change of variables formulas is contained in the following
lemma which follows {}from Proposition 3.33 in \cite{B memoirs}.

\begin{lem}\label{Gamma lemma} 
Let $(A,M), (B,N) \in \mathcal{H}$ and $\asqrt \in
(\bigwedge_\infty^0)^\times$.  If either $(A,M) = \mathbf{0}$ or $(B,N) =
\mathbf{0}$, then $\Gamma(\asqrt, A,M,B,N) = 0$.
\end{lem}

\begin{defn} {\em An} $N = 1$ SG-VOSA over $\bigwedge_*$ {\em is a
$\frac{1}{2} \mathbb{Z}$-graded (by weight) $\bigwedge_\infty$-module which
is also $\mathbb{Z}_2$-graded (by sign)  
\[V = \coprod_{k \in \frac{1}{2} \mathbb{Z}} V_{(k)} = \coprod_{k \in
 \frac{1}{2} \mathbb{Z}} V_{(k)}^0 \oplus \coprod_{k \in  \frac{1}{2}
\mathbb{Z}} V_{(k)}^1 = V^0 \oplus V^1 \]  
such that only the subspace $\bigwedge_*$ of $\bigwedge_\infty$ acts
non-trivially on $V$, 
\[\dim V_{(k)} < \infty \quad \mbox{for} \quad k \in \frac{1}{2}
\mathbb{Z} ,\] 
and for any $n \in \mathbb{N}$, a map
\[ \nu_n : SK(n) \rightarrow S\mathcal{F}_V (n) \]
satisfying the following axioms: 

(1) Positive energy axiom: 
\[V_{(k)} = 0 \quad \mbox{for $k$ sufficiently small.} \]

(2) Grading axiom:  Let $v' \in V'$, $v \in V_{(k)}$, and $a \in
(\bigwedge_\infty^0)^\times$.  Then 
\[ \langle v', \nu_1(\mathbf{0},(a, \mathbf{0}))(v)\rangle =
a^{-2k} \langle v',v \rangle .\]

(3) Supermeromorphicity axiom: For any $n \in \Z$, $v' \in V'$,
and $v_1,...,v_n \in V$, the function
\[Q \mapsto \langle v', \nu_n(Q) (v_1 \otimes \cdots \otimes v_n)
\rangle \]
on $SK(n)$ is a canonical supermeromorphic superfunction (in the
sense of (\ref{meromorphic})), and if $(z_i, \theta_i)$ and $(z_j,
\theta_j)$, for $i,j \in \{1,...,n\}$, $i \neq j$, are the $i$-th and
$j$-th punctures of $Q \in SK(n)$, respectively,  then for any $v_i$
and $v_j$ in $V$, there exists $N(v_i,v_j) \in \Z$ such that for any
$v' \in V'$ and $v_k \in V$, $k \neq i,j$, the order of the pole
$(z_i, \theta_i) = (z_j, \theta_j)$ of $\langle v', \nu_n(Q) (v_1
\otimes \cdots \otimes v_n) \rangle$ is less then $N(v_i,v_j)$.  

(4) Permutation axiom:  Let $\sigma \in S_n$.  Then for any $Q
\in SK(n)$
\[\sigma(\nu_n (Q)) = \nu_n(\sigma(Q)) . \]

(5) Sewing axiom:  There exists a unique complex number $c$ (the}
central charge {\em or} rank{\em) such that if $Q_1 \in SK(m)$ and
$Q_2 \in SK(n)$ are given by 
\begin{multline*}
Q_1 = ((z_1, \theta_1),...,(z_{m-1}, \theta_{m-1}); (A^{(0)},
M^{(0)}), (a^{(1)}, A^{(1)}, M^{(1)}),\\
...,(a^{(m)}, A^{(m)}, M^{(m)})), 
\end{multline*} 
and
\[Q_2 = ((z_1', \theta_1'),...,(z_{n-1}', \theta_{n-1}'); (B^{(0)},
N^{(0)}), (b^{(1)}, B^{(1)}, N^{(1)}),...,(b^{(n)}, B^{(n)},
N^{(n)})),\] 
and if the $i$-th tube of $Q_1$ can be sewn with the 0-th tube of $Q_2$,
then for any $v' \in V'$ and $v_1,...,v_{m + n - 1} \in V$, 
\[\langle v', (\nu_m (Q_1) \; _i*_0 \; \nu_n (Q_2))_{t^{1/2}} (v_1
\otimes \cdots \otimes v_{m + n - 1}) \rangle \]
is absolutely convergent when $t^{1/2} = 1$, and
\[\nu_{m + n - 1} (Q_1 \; _i\infty_0 \; Q_2) = (\nu_m (Q_1) \;
_i*_0  \; \nu_n (Q_2)) \; e^{-\Gamma (a^{(i)}, A^{(i)}, M^{(i)},
B^{(0)}, N^{(0)})c} . \] }
\end{defn}

We denote the $N=1$ supergeometric vertex operator superalgebra defined
above by 
\[(V, \nu = \{\nu_n \}_{n \in \mathbb{N}})\]
or just by $(V, \nu)$.   

\subsection{The isomorphism between the category of $N=1$ SG-VOSAs and the
category of $N=1$ NS-VOSAs}\label{iso} 

In this section we recall {}from \cite{B iso-thm} how to construct an $N=1$
NS-VOSA {}from an $N=1$ SG-VOSA, and how to construct an $N=1$ SG-VOSA {}from
an $N=1$ NS-VOSA.  We then recall the Isomorphism Theorem {}from 
\cite{B iso-thm} which states that the two notions are equivalent by showing
that their respective categories are isomorphic.  But first we must recall
the $\iota$ function introduced in \cite{B thesis} and \cite{B iso-thm}
which maps rational superfunctions to power series expanded about certain
variables.  This $\iota$ function is a generalization of the $\iota$
function of \cite{FLM}.

Let $\bigwedge_\infty[x_1,x_2,...,x_n]_S$ be the ring of
rational functions obtained by inverting (localizing with respect to) the
set 
\[S = \biggl\{\sum_{i = 1}^{n} a_i x_i \; | \; a_i \in
\mbox{$\bigwedge_\infty^0$},
\; \mbox{not all} \; (a_i)_B = 0\biggr\} . \] 
Recall the map $\iota_{i_1 ... i_n} : \mathbb{F}[x_1,...,x_n]_S
\longrightarrow \mathbb{F}[[x_1, x_1^{-1},...,  x_n, x_n^{-1}]]$ defined
in \cite{FLM} where coefficients of elements in $S$ are restricted to the
field $\mathbb{F}$.  We extend this map to
$\bigwedge_\infty[x_1,x_2,...,x_n]_S[\varphi_1,\varphi_2,...,\varphi_n] =
\bigwedge_\infty[x_1,\varphi_1,x_2,\varphi_2,...,x_n,\varphi_n]_{S}$ in
the obvious way obtaining
\[\iota_{i_1 ... i_n} : \mbox{$\bigwedge_\infty$}
[x_1,\varphi_1,...,x_n,\varphi_n]_S
\longrightarrow \mbox{$\bigwedge_\infty$}[[x_1, x_1^{-1},..., x_n,
x_n^{-1}]][\varphi_1,...,\varphi_n] .\]
Let $\bigwedge_\infty[x_1, \varphi_1, x_2, \varphi_2,...,x_n,
\varphi_n]_{S'}$ be the ring of rational functions obtained by inverting
the set   
\[S' = \biggl\{\sum_{\stackrel{i,j = 1}{i<j}}^{n} (a_i x_i + a_{ij} \varphi_i
\varphi_j) \; | \; a_i, a_{ij} \in \mbox{$\bigwedge_\infty^0$}, \;
\mbox{not all} \; (a_i)_B = 0\biggr\}. \]
Since we use the convention that a function of even and odd
variables should be expanded about the even variables, we have 
\[\frac{1}{\sum_{\stackrel{i,j = 1}{i<j}}^{n} (a_i x_i + a_{ij} \varphi_i
\varphi_j)} = \frac{1}{\sum_{i = 1}^{n} a_i x_i} -
\frac{\sum_{\stackrel{i,j = 1}{i<j}}^{n} a_{ij} \varphi_i
\varphi_j}{(\sum_{i = 1}^{n} a_i x_i)^2} .\]
Thus $\bigwedge_\infty[x_1,\varphi_1,x_2,\varphi_2,...,x_n,\varphi_n]_{S'}
\subseteq
\bigwedge_\infty[x_1,\varphi_1,x_2,\varphi_2,...,x_n,\varphi_n]_S$,    and
$\iota_{i_1 ... i_n}$ is well defined on
$\bigwedge_\infty[x_1,\varphi_1,x_2,\varphi_2,...,x_n,\varphi_n]_{S'}$.

In the case $n = 2$, the map $\iota_{12} : \bigwedge_\infty
[x_1,\varphi_1,x_2,\varphi_2]_{S'} \longrightarrow \bigwedge_\infty
[[x_1,x_2]][\varphi_1,\varphi_2]$ is given by first expanding an
element of $\bigwedge_\infty [x_1,\varphi_1,x_2,\varphi_2]_{S'}$ as a
formal series in $\bigwedge_\infty[x_1,\varphi_1,x_2,\varphi_2]_S$ and
then expanding each term as a series in $\bigwedge_\infty [[x_1,x_2]]
[\varphi_1,\varphi_2]$ containing at most finitely many negative powers
of $x_2$ (using binomial expansions for negative powers of linear
polynomials involving both $x_1$ and $x_2$). 

Note that in the case of functions of one even variable $x_1$ and one odd
variable $\varphi_1$, the function $\iota_1$ is the identity.

Given an $N=1$ supergeometric VOSA, we construct an (algebraic) $N=1$ NS-VOSA. Let
$(V,\nu)$ be an $N=1$ SG-VOSA over $\bigwedge_*$.  We define the vacuum
$\bf{1}_\nu \in \bar{V}$ by $\mathbf{1}_\nu = \nu_0 (\mathbf{0})$;
an element $\tau_\nu \in \bar{V}$ by $\tau_\nu = (- \partial / \partial
\epsilon) \nu_0 (\mathbf{0}, M(\epsilon, 3/2))$, where $M(\epsilon,
3/2)$ is the series of odd supernumbers given by $M(\epsilon, 3/2) = \{0,
\epsilon, 0, 0,...\}$; and the vertex operator $Y_\nu (u, (x, \varphi)) =
\sum_{n \in \mathbb{Z}} u_n x^{-n-1} + \varphi \sum_{n \in \mathbb{Z}}
u_{n - 1/2} x^{-n-1}$ associated with $u \in V$ by
\begin{equation}\label{define Y}
u_n v + \theta u_{n - \frac{1}{2}} v = \mbox{Res}_z \left( z^n \nu_2
((z,\theta); \mathbf{0}, (1, \mathbf{0}), (1, \mathbf{0}) ) (v' \otimes u
\otimes v) \right) , 
\end{equation}
where $\mbox{Res}_z$ means taking the residue at the singularity $z =
0$, i.e., taking the coefficient of $z^{-1}$.

\begin{prop}\label{get a vosa} (\cite{B iso-thm}) 
The elements $\mathbf{1}_\nu$ and $\tau_\nu$ of $\bar{V}$ are in fact in
$V_{(0)}$ and $V_{(3/2)}$, respectively.  If the rank of $(V,\nu)$ is $c$,
then $(V, Y_\nu( \cdot, (x,\varphi)), \mathbf{1}_\nu, \tau_\nu)$ is an $N=1$
NS-VOSA with odd formal variables and with rank $c$ . 
\end{prop}

Given an $N=1$ NS-VOSA over $\bigwedge_*$ with rank $c \in \mathbb{C}$,
$(V, Y(\cdot,(x,\varphi)),\bf{1}, \tau)$, we construct an $N=1$
SG-VOSA by defining a sequence of maps: 
\begin{eqnarray*}
\nu_n^Y : SK(n) & \rightarrow & SF_V(n) \\
Q & \mapsto & \nu_n^Y (Q)
\end{eqnarray*}
by
\begin{multline*}
\langle v', \nu_n^Y \left( (z_1, \theta_1), ...,(z_{n-1},
\theta_{n-1}); (A^{(0)}, M^{(0)}), (a^{(1)}, A^{(1)}, M^{(1)}), 
... \right. \\
\left. ... (a^{(n)}, A^{(n)}, M^{(n)}) \right) (v_1 
\otimes \cdots \otimes v_n) \rangle 
\end{multline*} 
\begin{multline}\label{define nu} 
= \iota^{-1}_{1 \cdots n-1} \langle e^{- \sum_{j \in \Z} \left(
A^{(0)}_j L'(j) + M^{(0)}_{j - \frac{1}{2}} G'(j - \frac{1}{2})
\right)} v',\\
Y(e^{- \sum_{j \in \Z} \left(A^{(1)}_j L(j) + M^{(1)}_{j -
\frac{1}{2}} G(j - \frac{1}{2}) \right)} \cdot (a^{(1)})^{-2L(0)}
\cdot v_1, (x_1, \varphi_1)) \cdots\\
Y(e^{- \sum_{j \in \Z} \left(A^{(n-1)}_j L(j) + M^{(n-1)}_{j -
\frac{1}{2}} G(j - \frac{1}{2}) \right)} \cdot (a^{(n-1)})^{-2L(0)} 
\cdot v_{n-1}, (x_{n-1}, \varphi_{n-1})) \cdot \\
\left. e^{- \sum_{j \in \Z} \left(A^{(n)}_j L(j) +
M^{(n)}_{j - \frac{1}{2}} G(j - \frac{1}{2}) \right)} \cdot 
(a^{(n)})^{-2L(0)} \cdot v_n \rangle \right|_{(x_i,\varphi_i) =
(z_i,\theta_i)} 
\end{multline}
for $n \in \Z$, $v' \in V'$, and $v_1,...,v_n \in V$, and by
\[\langle v', \nu_0^Y (A^{(0)}, M^{(0)})\rangle = \langle e^{- \sum_{j
\in \Z} \left( A^{(0)}_j L'(j) + M^{(0)}_{j - \frac{1}{2}} G'(j -
\frac{1}{2}) \right)} v', \mathbf{1} \rangle . \]

\begin{prop}\label{get a sgvosa} (\cite{B iso-thm})
The pair $(V,\nu^Y)$ is an $N=1$ SG-VOSA over $\bigwedge_*$.    
\end{prop}

Let $\mathbf{SG}(c,*)$ be the category of $N=1$ SG-VOSAs over $\bigwedge_*$
with central charge $c$, and let $\mathbf{SV}(\varphi, c,*)$ be the category
of $N=1$ NS-VOSAs over $\bigwedge_*$ with odd formal variables and central
charge $c$ where $*$ is $L \in \mathbb{N}$ or $\infty$.  Let
$1_{SV}$ and $1_{SG}$ be the identity functors on
$\mathbf{SV}(\varphi,c,*)$ and $\mathbf{SG}(c,*)$, respectively.

\begin{thm}\label{iso theorem} (\cite{B iso-thm})
(\cite{B iso-thm}) For any $c \in \mathbb{C}$, the two categories
$\mathbf{SV}(\varphi,c,*)$ and $\mathbf{SG}(c,*)$ are isomorphic with
isomorphism given by the following functors acting on objects and morphisms 
\begin{eqnarray*}
F_{SV}: \mathbf{SV}(\varphi,c,*) &\longrightarrow& \mathbf{SG}(c,*) \\
((V, Y(\cdot,(x,\varphi)),\bf{1}, \tau), \gamma) &\mapsto& ((V, \nu^Y),
\gamma)
\end{eqnarray*}
and
\begin{eqnarray*}
F_{SG}: \mathbf{SG}(c,*) &\longrightarrow& \mathbf{SV}(\varphi,c,*)\\
((V, \nu),\gamma) & \mapsto & ((V, Y_\nu( \cdot, (x,\varphi)),
\mathbf{1}_\nu, \tau_\nu), \gamma)\\ 
\end{eqnarray*} 
which satisfy $F_{SV} \circ F_{SG} = 1_{SG}$ and $F_{SG} \circ F_{SV} =
1_{SV}$.
\end{thm}

\subsection{Geometric proof of Proposition \ref{change prop} for
superconformal change of coordinates convergent in a neighborhood of
zero}\label{section-geom1}

\begin{proof}
Assuming that $H(x,\varphi)$ given by (\ref{coordinate change}) is
convergent and bijective in a neighborhood of zero, we have
\begin{eqnarray*}
Q_1 &=& (\mathbf{0}, (\asqrt, A, M)) \in SK(1),
\end{eqnarray*}
and
\begin{eqnarray*}
Q_2 &=&((z,\theta); \mathbf{0}, (1,\mathbf{0}), (1,\mathbf{0})) \in SK(2).
\end{eqnarray*}
{}From Example \ref{example1}, we see that if $(z,\theta)$ is in the radius 
of convergence of $H$, then $Q_2$ can be sewn to $Q_1$ and 
\begin{eqnarray}
Q_1 \; _1\infty_0 \; Q_2 
&=& (H (z,\theta); \mathbf{0}, \hat{E}^{-1}(s_{(z,\theta)} \circ H^{-1}
\circ s_{H(z,\theta)}^{-1} (w,\rho)), (\asqrt, A, M)) \nonumber \\
&=& (H (z,\theta); \mathbf{0}, (\asqrt e^{\Theta^{(1)}_0},
\{\asqrt^{2j} \Theta^{(1)}_j, \asqrt^{2j - 1} \Theta^{(1)}_{j
- \frac{1}{2}}\}_{j\in \Z}), (\asqrt, A, M)). 
\end{eqnarray}
Therefore {}from the Isomorphism Theorem \ref{iso theorem}, the sewing
axiom for an $N=1$ SG-VOSA, and Lemma \ref{Gamma lemma} we have that 
\begin{eqnarray}\label{new Y equation}
\lefteqn{\left. \langle v', \gamma_H (Y(u,(x,\varphi))v)\rangle 
\right|_{(x,\varphi) = (z,\theta)}} \nonumber \\
&=& \iota_1^{-1} \langle v', e^{ - \sum_{j \in \mathbb{Z}_+} \left( A_j L(j)
+ M_{j - \frac{1}{2}} G(j-\frac{1}{2}) \right) } \cdot  \asqrt^{-2L(0)}
\cdot Y(u,(z,\theta))v\rangle
\nonumber\\ 
&=& \langle v', (\nu^Y_1(Q_1) \; _1*_0 \; \nu_2^Y(Q_2))
(u\otimes v)\rangle
\nonumber \\
&=& \langle v', (\nu_2^Y(Q_1 \; _1\infty_0 \; Q_2) (u\otimes v)\rangle
\nonumber \\
&=& \left. \langle v', Y(\gamma_{\Theta(H)} (u),H(x,\varphi))
\gamma_H(v)\rangle 
\right|_{(x,\varphi) = (z,\theta)}
\end{eqnarray}
for $u,v \in V$, $v' \in V'$, and $(z,\theta)$ in the domain of convergence
of $H$. Thus $\gamma_{\Theta(H)}$ is well defined and we have that
\[ \gamma_H (Y(u,(x,\varphi))v) =  Y(\gamma_{\Theta(H)}(u),H(x,\varphi))
\gamma_H(v) ,\]
as formal power series in $V[[x^{-1},x]] [\varphi]$, with both sides
convergent and equal if $(x,\varphi) = (z,\theta)$ is in the domain of
convergence of $H$.  Since $\gamma_{\Theta(H,t^{1/2},
(x,\varphi))}|_{t^{1/2} = 1}$ is well-defined, so is $\Theta_j^{(1)}
(t^{1/2})|_{t^{1/2} = 1}$ and thus so is $\gamma_{\Theta(H, t^{1/2},
(x,\varphi))}^{-1}|_{t^{1/2} = 1}$.  This, the fact that $\gamma_{\Theta
(H, (x,\varphi))} \circ \gamma_{\Theta(H, (x,\varphi))}^{-1} (v) = v$ for
$v \in V$, and the fact that $\gamma_H$ is invertible, give us equation
(\ref{change formula rewritten}).  
\end{proof}

\subsection{Geometric proof of Proposition \ref{change prop2} for
superconformal change of coordinates convergent in a neighborhood of
infinity}\label{section-geom2}

\begin{proof}
Let 
\[Q_1 = ((z,\theta); \mathbf{0}, (1,\mathbf{0}), (1,\mathbf{0})) \in SK(2) .
\]
Assuming that $H^{-1}(x,\varphi)$ given by (\ref{coordinate change22}) is
convergent and bijective in a neighborhood of infinity, we have
\[Q_2 = ((B,N), (1,\mathbf{0})) \in SK(1).\]
{}From Example \ref{example2}, we see that if $(z,\theta)$ is in the radius 
of convergence of $H \circ I$, then $Q_2$ can be sewn to $Q_1$ and 
\begin{eqnarray}
\lefteqn{Q_1 \; _2\infty_0 \; Q_2} \nonumber \\
&=& (H \circ I (z,\theta); (B,N) , \hat{E}^{-1}(s_{(z,\theta)} \circ I^{-1}
\circ H^{-1} \circ s_{H \circ I (z,\theta)}^{-1} (w,\rho)), (1, \mathbf{0}))
\nonumber \\  
&=& (H \circ I (z,\theta); (B,N) , (e^{\Theta^{(2)}_0}, \{\Theta^{(2)}_j, 
\Theta^{(2)}_{j -\frac{1}{2}}\}_{j\in \Z}), (1, \mathbf{0})). 
\end{eqnarray}
Therefore {}from the Isomorphism Theorem \ref{iso theorem}, the sewing
axiom for an $N=1$ SG-VOSA, and Lemma \ref{Gamma lemma} we have that 
\begin{eqnarray}\label{new Y equation2}
\lefteqn{\left. \langle v', Y(u,(x,\varphi))\xi_{H \circ I}^*(v)\rangle
\right|_{(x,\varphi) = (z,\theta)}} \nonumber \\ 
&=&  \sum_{k \in \frac{1}{2}\mathbb{Z}}  \left. \langle v', Y(u,(x,\varphi))
P_k(\xi_{H \circ I}^*(v))\rangle \right|_{(x,\varphi) = (z,\theta)} \\
&=& \langle v', (\nu_2^Y(Q_1) \; _2*_0 \; \nu_1^Y(Q_2)) (u\otimes v)\rangle  
\nonumber\\
&=& \langle v', \nu_2^Y(Q_1 \; _2\infty_0 \; Q_2) (u\otimes v)\rangle
\nonumber \\ 
&=& \left. \langle \xi_{H \circ I}(v'), Y(\xi_{\Theta(H \circ I,
(x,\varphi))} (u),H \circ I(x,\varphi))  v\rangle \right|_{(x,\varphi) =
(z,\theta)}
\nonumber \\  
&=& \left. \langle v', \xi_{H \circ I}^*(Y(\xi_{\Theta(H \circ I,
(x,\varphi))} (u), H \circ I (x,\varphi))  v)\rangle \right|_{(x,\varphi) =
(z,\theta)}.
\end{eqnarray}
for $u,v \in V$, $v' \in V'$ and $(z,\theta)$ in the radius of convergence
of $H \circ I$.  In particular, the above equality implies that 
\[\langle v', (\nu_2^Y(Q_1) \; _2*_0 \; \nu_1^Y(Q_2))_{t^{1/2}} (u\otimes
v)\rangle \]
is absolutely convergent when $t^{1/2}=1$, and is equal to
\begin{multline*}
\left. \langle v', Y(u,(x,\varphi))\xi_{H \circ I}^*(v)\rangle
\right|_{(x,\varphi) = (z,\theta)} \\
=  \left. \langle v', \xi_{H \circ I}^*(Y(\xi_{\Theta(H \circ I,
(x,\varphi))} (u), H \circ I (x,\varphi))  v)\rangle \right|_{(x,\varphi) =
(z,\theta)} ,
\end{multline*} 
implying $\xi_{\Theta(H \circ I,(x,\varphi))} = \xi_{\Theta(H \circ
I,t^{1/2}, (x,\varphi))}|_{t^{1/2} = 1}$ is well defined for $(z,\theta)$ in
the radius of convergence of $H \circ I$.  Thus, as formal power series in
$\bar{V} [[x^{-1},x]] [\varphi]$, we have
\begin{eqnarray}
Y(u,(x,\varphi))\xi_{H \circ I}^*(v) = \xi_{H \circ I}^*(Y(\xi_{\Theta(H
\circ I, (x,\varphi))} (u), H \circ I (x,\varphi))  v)
\end{eqnarray}
and both sides are well defined and equal for $(x,\varphi) = (z,\theta)$ in
the radius of convergence of $H \circ I$.  In addition, since
$\xi_{\Theta(H \circ I, t^{1/2}, (x,\varphi))}|_{t^{1/2} = 1}$ is
well-defined, so is $\Theta_j^{(2)} (t^{1/2})|_{t^{1/2} = 1}$ and thus so
is $\xi_{\Theta(H \circ I,t^{1/2}, (x,\varphi))}^{-1}|_{t^{1/2} = 1}$. 
This, the  fact that $\xi_{\Theta(H \circ I, (x,\varphi))} \circ
\xi_{\Theta(H \circ I, (x,\varphi))}^{-1} (v) = v$ for $v \in V$, and the
fact that $\xi^*_{H \circ I}$ is invertible, give us equation (\ref{change
formula2 rewritten}).
\end{proof}

\section{Formal algebraic proofs for formal superconformal changes of
variables}\label{alg-section}

\subsection{Formal algebraic proof of Proposition \ref{change prop} for
formal superconformal change of coordinates at zero}\label{alg-proof}

In this section, we let $H(x,\varphi) \in \bigwedge_\infty[[x]][\varphi]$ 
be the formal superconformal power series given by (\ref{coordinate
change}), but this time we do not assume that $H(x,\varphi)$ converges in
a neighborhood of zero. 

We will prove Proposition \ref{change prop} by a formal algebraic
argument first employed in \cite{B iso-thm} to construct an $N=1$
SG-VOSA {}from an $N=1$ NS-VOSA (i.e., to prove Proposition
\ref{get a sgvosa} above).  This result was then used in \cite{B iso-thm} to
prove the Isomorphism Theorem \ref{iso theorem}.  The difference in the
present case is that we are not assuming the convergence of $H$ and
consequently cannot use the Isomorphism Theorem \ref{iso theorem} as we did
in Section \ref{geom}.  However, as it turns out, enough of the algebraic
formalism used in \cite{B iso-thm} carries over to the current more general
setting, and we will be retracing that part of the construction  used in
\cite{B iso-thm} in this more general setting.

We first prove the following bracket formula.
  
\begin{lem}\label{bracket-Lemma}
Let $(V,Y(\cdot,(x,\varphi)), \mathbf{1}, \tau)$ be an $N=1$
NS-VOSA over a Grassmann algebra.  For $u \in V$, $t^{1/2}$ a formal even
variable, $\asqrt \in (\bigwedge_\infty^0)^\times$, and $(A_j,M_{j-1/2}) 
\in \bigwedge_\infty$, for $j \in \Z$, we have the following bracket formula
in $(\mathrm{End} \; V)[[t^{1/2}]][[x^{-1},x]][\varphi]$
\begin{equation*}
\Bigl[ \sum_{j \in \Z} \left( t^{j} \asqrt^{-2j} A_j L(j) +  t^{j
- \frac{1}{2}} \asqrt^{-2j + 1} M_{j - \frac{1}{2}} G(j -
\frac{1}{2}) \right), Y(u,(x,\varphi))
\Bigr] \hspace{.6in}
\end{equation*}
\begin{multline}\label{bracket}
= Y \Biggl( \sum_{k = -1}^{\infty} \sum_{j \in \Z}
\binom{j+1}{k+1} t^j \asqrt^{-2j} x^{j - k}  
\biggl(\!  \Bigl( A_j + 2\left(\frac{j-k}{j+1} \right)
t^{-\frac{1}{2}} \asqrt x^{-1} \varphi M_{j - \frac{1}{2}}
\Bigr)  L(k) \\
+ \; x^{-1} \Bigl( \left(\frac{j-k}{j+1} \right) t^{-\frac{1}{2}}
\asqrt M_{j - \frac{1}{2}} + \varphi \frac{(j-k)}{2}
A_j \Bigr)  G(k + \frac{1}{2})  \biggr)  u, (x,\varphi) .
\Biggr) 
\end{multline} 
Furthermore, if we apply both sides of (\ref{bracket}) to $v \in V$,
then we can set $t^{1/2} = 1$ and both sides are well-defined elements of
$V((x))[\varphi]$.  
\end{lem}

\begin{proof}
{}From (\ref{bracket relation for a vosa}), we have
\begin{multline*} 
\left[ Y(\tau, (x_1,\varphi_1)), Y(u,(x,\varphi)) \right] \\
= \mbox{Res}_{x_0} x^{-1} \delta \left( \frac{x_1 - x_0 -
\varphi_1 \varphi}{x} \right) Y(Y(\tau,(x_0, \varphi_1 -
\varphi))u,(x, \varphi)) .
\end{multline*} 
Let 
\begin{equation*}
l_t(x_1) = \sum_{j \in \Z} t^{j} \asqrt^{-2j} A_j x_1^{j + 1}
\qquad \mathrm{and} \qquad 
g_t(x_1) = \sum_{j \in \Z} t^{j - \frac{1}{2}}
\asqrt^{-2j + 1} M_{j - \frac{1}{2}} x_1^{j} .
\end{equation*}
Then 
\begin{eqnarray*}
\varphi_1 \mbox{Res}_{x_1} l_t(x_1) Y(\omega, (x_1,\varphi_1)) &=&
\varphi_1 \sum_{j \in \Z} t^j \asqrt^{-2j} A_j  L(j) \\ 
\varphi_1 \mbox{Res}_{x_1} g_t(x_1) Y(\tau, (x_1,\varphi_1)) &=&
\varphi_1 \sum_{j \in \Z} t^{j - \frac{1}{2}} \asqrt^{-2j + 1} M_{j -
\frac{1}{2}} G(j - \frac{1}{2}) ,
\end{eqnarray*}
and by the $G(-1/2)$-derivative property
(\ref{G(-1/2)-derivative}), we have 
\begin{eqnarray*}
Y(\omega, (x_1,\varphi_1)) &=& Y( \frac{1}{2} G(-\frac{1}{2}) \tau,
(x_1, \varphi_1)) \\
&=& \frac{1}{2} \Bigl( \frac{\partial}{\partial \varphi_1} + \varphi_1
\frac{\partial}{\partial x_1} \Bigr) Y(\tau, (x_1,\varphi_1))  .
\end{eqnarray*}
Thus 
\begin{eqnarray*}
& & \hspace{-.5in} \varphi_1 \biggl[ \sum_{j \in \Z} \left( t^j
\asqrt^{-2j} A_j L(j) + t^{j - \frac{1}{2}} \asqrt^{-2j + 1}
M_{j - \frac{1}{2}} G(j - \frac{1}{2}) \right) , Y(u,(x,\varphi))
\biggr] \\ 
&=& \varphi_1 \mbox{Res}_{x_1} \left[ \Bigl( \frac{1}{2}l_t(x_1)
\Bigl( \frac{\partial}{\partial \varphi_1} + \varphi_1
\frac{\partial}{\partial x_1} \Bigr) + g_t(x_1) \Bigr) Y(\tau,
(x_1,\varphi_1)), Y(u,(x,\varphi)) \right] \\
&=& \varphi_1 \mbox{Res}_{x_1} \left( \Bigl( \frac{1}{2}l_t(x_1) \Bigl(
\frac{\partial}{\partial \varphi_1} + \varphi_1
\frac{\partial}{\partial x_1} \Bigr) + g_t(x_1) \Bigr) \right. \\
& & \hspace{.7in} \left. \mbox{Res}_{x_0} x^{-1} \delta \Bigl(
\frac{x_1 - x_0 - \varphi_1 \varphi}{x} \Bigr) Y(Y(\tau,(x_0,
\varphi_1 - \varphi))u,(x, \varphi)) \right)\\
&=& Y\left(\varphi_1 \mbox{Res}_{x_1}\mbox{Res}_{x_0} \Bigl(
\frac{1}{2}l_t(x_1) \Bigl( \frac{\partial}{\partial \varphi_1} +
\varphi_1 \frac{\partial}{\partial x_1} \Bigr) + g_t(x_1) \Bigr)
\right. \\
& & \left. \hspace{1in} x^{-1} \delta \Bigl( \frac{x_1 - x_0 -
\varphi_1 \varphi}{x} \Bigr) Y(\tau,(x_0, \varphi_1 - \varphi))u,(x,
\varphi) \right) .  
\end{eqnarray*}
However using the $\delta$-function identity (\ref{delta 2 terms with
phis}), we have
\begin{multline*}
\varphi_1 \mbox{Res}_{x_1}\mbox{Res}_{x_0} \left(
\frac{1}{2}l_t(x_1) \Bigl( \frac{\partial}{\partial \varphi_1} +
\varphi_1 \frac{\partial}{\partial x_1} \Bigr) + g_t(x_1) \right) \\
x^{-1} \delta \Bigl( \frac{x_1 - x_0 - \varphi_1
\varphi}{x} \Bigr) Y(\tau,(x_0, \varphi_1 - \varphi)) 
\end{multline*}
\begin{multline*}
=  \varphi_1 \mbox{Res}_{x_1}\mbox{Res}_{x_0} \left( \frac{1}{2}
l_t(x_1) \Bigl( \frac{\partial}{\partial \varphi_1} + \varphi_1
\frac{\partial}{\partial x_1} \Bigr) + g_t(x_1) \right) \\
x_1^{-1} \delta \Bigl( \frac{x + x_0 + \varphi_1
\varphi}{x_1} \Bigr) Y(\tau,(x_0, \varphi_1 - \varphi)) 
\end{multline*}
\begin{multline*}
= \varphi_1 \mbox{Res}_{x_1}\mbox{Res}_{x_0} \left(
\frac{1}{2}l_t(x_1) \Bigl( \frac{\partial}{\partial \varphi_1} 
x_1^{-1} \delta \Bigl( \frac{x + x_0 + \varphi_1 \varphi}{x_1} \Bigr)
\Bigr) Y(\tau,(x_0,- \varphi)) \right.\\ 
 + \; \frac{1}{2}l_t(x_1) x_1^{-1} \delta \Bigl(
\frac{x + x_0}{x_1} \Bigr) \frac{\partial}{\partial \varphi_1}
Y(\tau, (x_0,\varphi_1)) \\
\left. + \; g_t(x_1) x_1^{-1} \delta \Bigl( \frac{x
+ x_0}{x_1} \Bigr) Y(\tau,(x_0,- \varphi)) \right) 
\end{multline*}
\begin{multline*}
= \varphi_1 \mbox{Res}_{x_1} \mbox{Res}_{x_0} \biggl( 
\biggl(\sum_{n \in \mathbb{Z}} n (x + x_0)^{n-1} \varphi x_1^{-n-1}
\biggr) \biggl( \frac{1}{2} \sum_{j \in \Z} t^j \asqrt^{-2j} A_j
x_1^{j + 1}
\biggr) \biggr.\\
\biggl( \sum_{k \in \mathbb{Z}} G(k + \frac{1}{2})
x_0^{-k - 2} \biggr) 
+ \biggl(\sum_{n \in \mathbb{Z}}(x + x_0)^n x_1^{-n-1} \biggr) 
\biggl( \frac{1}{2} \sum_{j \in \Z} t^j \asqrt^{-2j} A_j  x_1^{j
+ 1} \biggr)\\
\biggl( 2\sum_{k \in \mathbb{Z}} L(k) x_0^{-k - 2} \biggr)
+ \biggl(\sum_{n \in \mathbb{Z}}(x + x_0)^n x_1^{-n-1} \biggr)
\biggl( \sum_{j \in \Z} t^{j-\frac{1}{2}} \asqrt^{-2j +1}
M_{j - \frac{1}{2}} x_1^{j} \biggr) \\
\biggl. \biggl( \sum_{k \in \mathbb{Z}} G(k + \frac{1}{2})
x_0^{-k - 2} - 2\varphi L(k) x_0^{-k - 2} \biggr) \biggr)
\end{multline*}
\begin{multline*}
= \varphi_1 \mbox{Res}_{x_0} \biggl( \sum_{l \in \mathbb{N}} \sum_{k \in
\mathbb{Z}} \sum_{j \in \Z} \biggl( (j + 1) \binom{j}{l} x^{j -l} x_0^l
\frac{\varphi}{2} t^j \asqrt^{-2j} A_j G(k + \frac{1}{2})
x_0^{-k - 2} \\
+ \; \binom{j + 1}{l} x^{j-l + 1} x_0^l t^j \asqrt^{-2j}
A_j L(k) x_0^{-k - 2} \\
+ \; \binom{j}{l} x^{j-l} x_0^l t^{j-\frac{1}{2}} \asqrt^{-2j + 1}
M_{j - \frac{1}{2}} \Bigl(G(k + \frac{1}{2}) - 2\varphi L(k) \Bigr)
x_0^{-k - 2}
\biggr) \biggr) 
\end{multline*}
\begin{multline*}
= \varphi_1 \sum_{l \in \mathbb{N}} \sum_{j \in \Z} \left( \frac{(j +
1)}{2} \binom{j}{l} x^{j -l} \varphi t^j \asqrt^{-2j} A_j
G(l - \frac{1}{2}) \right. \\
+ \; \binom{j + 1}{l} x^{j-l + 1} t^j \asqrt^{-2j} A_j
L(l - 1) \\
\left. + \; \binom{j}{l} x^{j-l} t^{j-\frac{1}{2}} \asqrt^{-2j +1}
M_{j - \frac{1}{2}} (G(l - \frac{1}{2}) - 2\varphi L(l - 1)) \right) 
\end{multline*}
\begin{multline*}
= \varphi_1 \sum_{k = -1}^\infty \sum_{j \in \Z} \biggl( \binom{j +
1}{k+1} t^j \asqrt^{-2j} x^{j-k} \Bigl(A_j + 2 \left( \frac{j-k}{j + 1} \right)
x^{-1} t^{-\frac{1}{2}} \asqrt \varphi M_{j - \frac{1}{2}} \Bigr) L(k) \\
+ \; \binom{j}{k + 1} t^j \asqrt^{-2j} x^{j -k -1} \Bigl( \frac{(j + 1)}{2} 
\varphi  A_j + t^{-\frac{1}{2}} \asqrt M_{j - \frac{1}{2}} \Bigr) G(k +
\frac{1}{2})  \biggr) 
\end{multline*}
which finishes the proof of (\ref{bracket}).

Since $V$ is a positive energy representation for the Neveu-Schwarz
algebra and by the truncation condition for $Y(\cdot, (x,\varphi))$, 
when applied to $v \in V$, each side of (\ref{bracket}) is in
$V((x))[\varphi][t^{1/2}]$ and thus well defined when $t^\frac{1}{2} = 1$.
\end{proof}

We have the following immediate corollary to Lemma \ref{bracket-Lemma}.

\begin{cor}\label{theta1-corollary}
\begin{multline*}
e^{- \sum_{j \in \Z} \left( t^{j} \asqrt^{-2j} A_j L(j) +  t^{j -
\frac{1}{2}} \asqrt^{-2j + 1} M_{j - \frac{1}{2}} G(j - \frac{1}{2})
\right)}  Y(u,(x,\varphi))
\cdot  \\
e^{\sum_{j \in \Z} \left( t^{j} \asqrt^{-2j} A_j L(j) + 
t^{j - \frac{1}{2}} \asqrt^{-2j + 1} M_{j - \frac{1}{2}} G(j -
\frac{1}{2}) \right)}  
\end{multline*} 
\begin{multline*}
= \; Y \Biggl(\exp \Biggl(\! - \! \!
\sum_{k = -1}^{\infty} \sum_{j \in \Z} \binom{j+1}{k+1} t^j \asqrt^{-2j} x^{j -
k}  \\
\biggl(\!  \Bigl( A_j + 2\left(\frac{j-k}{j+1} \right) t^{-\frac{1}{2}} \asqrt
x^{-1} \varphi M_{j - \frac{1}{2}} \Bigr)  L(k) \\
+ \; x^{-1} \Bigl( \left(\frac{j-k}{j+1} \right) t^{-\frac{1}{2}}
\asqrt M_{j - \frac{1}{2}} + \varphi \frac{(j-k)}{2} A_j \Bigr)  G(k +
\frac{1}{2}) \biggl) \Biggr)  u, (x,\varphi) \Biggr)  .
\end{multline*} 
Furthermore, if we apply both sides of the equation above to $v \in V$, we 
can set $t^\frac{1}{2} = 1$.  
\end{cor}

We are now ready to prove Proposition \ref{change prop} in the case of
a formal invertible superconformal change of coordinates vanishing at zero. 

\begin{proof} 
Denote
\begin{eqnarray*}
\lefteqn{H_{t^{1/2}} (x,\varphi)} \\
&=& (\tilde{x} (t^\frac{1}{2}), \tilde{\varphi} (t^\frac{1}{2})) \\ 
&=& (t^{-\frac{1}{2}} \asqrt)^{2L_0(x,\varphi)} \cdot \exp \Bigl(\sum_{j \in
\Z} \Bigl(A_j L_j(x,\varphi) + M_{j - \frac{1}{2}} G_{j - \frac{1}{2}}
(x,\varphi) \Bigr) \Bigr) \cdot (x,\varphi) . 
\end{eqnarray*}
By Corollary \ref{theta1-corollary}, Proposition \ref{first Theta identity
in End}, equation (\ref{exponential L(-1) and  G(-1/2) property}), and
equation (\ref{conjugate by L(0)}), we have
\begin{eqnarray*}
& & \hspace{-.3in} \gamma_{H, t^{1/2}} (Y(u, (x,\varphi))v)   \\
&=& e^{- \sum_{j \in \Z} \left( A_j L(j) + M_{j - \frac{1}{2}} G(j -
\frac{1}{2}) \right)} \cdot (t^{-\frac{1}{2}} \asqrt)^{-2L(0)}
Y(u,(x,\varphi)) v  \\ 
&=& e^{- \sum_{j \in \Z} \left( A_j L(j) + M_{j - \frac{1}{2}} G(j -
\frac{1}{2}) \right)} \cdot (t^{-\frac{1}{2}} \asqrt)^{-2L(0)}
Y(u,(x,\varphi)) \cdot (t^{-\frac{1}{2}}\asqrt)^{2L(0)} \cdot \\
& & \hspace{.1in} e^{\sum_{j \in \Z} \left( A_j L(j) + M_{j - \frac{1}{2}}
G(j - \frac{1}{2}) \right)} \cdot e^{- \sum_{j \in \Z} \left( A_j L(j) + 
M_{j -\frac{1}{2}} G(j - \frac{1}{2}) \right)} \cdot \\
&& \hspace{3.3in} (t^{-\frac{1}{2}} \asqrt)^{-2L(0)} \cdot v  
\end{eqnarray*}
\begin{multline*} 
= \; \; (t^{-\frac{1}{2}} \asqrt)^{-2L(0)} e^{ - \sum_{j \in \Z}
\left( t^{j} \asqrt^{-2j} A_j L(j) + t^{j - \frac{1}{2}}
\asqrt^{-2j + 1} M_{j-\frac{1}{2}} G(j - \frac{1}{2}) \right) }
Y(u,(x,\varphi)) \\
e^{\sum_{j \in \Z} \left(t^{j} \asqrt^{-2j} A_j L(j) +  t^{j -
\frac{1}{2}} \asqrt^{-2j + 1} M_{j - \frac{1}{2}} G(j - \frac{1}{2})
\right)} (t^{-\frac{1}{2}} \asqrt)^{2L(0)} \cdot \\
e^{- \sum_{j \in \Z} \left( A_j L(j) + M_{j -\frac{1}{2}} G(j - \frac{1}{2})
\right)} \cdot (t^{-\frac{1}{2}} \asqrt)^{-2L(0)} \cdot v  
\end{multline*}
\begin{multline*}
=  \; \; (t^{-\frac{1}{2}} \asqrt)^{-2L(0)} Y \Biggl(\exp \Biggl(\! - \! \!
\sum_{k = -1}^{\infty} \sum_{j \in \Z} \binom{j+1}{k+1} t^j \asqrt^{-2j}
x^{j - k}  \\ 
\biggl(\!  \Bigl( A_j + 2\left(\frac{j-k}{j+1} \right)
t^{-\frac{1}{2}} \asqrt x^{-1} \varphi M_{j - \frac{1}{2}}
\Bigr)  L(k) \\
+ \; x^{-1} \Bigl( \left(\frac{j-k}{j+1} \right) t^{-\frac{1}{2}}
\asqrt M_{j - \frac{1}{2}} + \varphi \frac{(j-k)}{2}
A_j \Bigr)  G(k + \frac{1}{2}) \biggl) \Biggr)  u, (x,\varphi)
\Biggr) \\
(t^{-\frac{1}{2}} \asqrt)^{2L(0)} \cdot e^{- \sum_{j \in \Z} \left( A_j
L(j) + M_{j -\frac{1}{2}} G(j - \frac{1}{2})
\right)} \cdot (t^{-\frac{1}{2}} \asqrt)^{-2L(0)} \cdot v  \\
\end{multline*}
\begin{multline*}
= \; \; (t^{-\frac{1}{2}} \asqrt)^{-2L(0)}  Y \Bigl( e^{((t^{-\frac{1}{2}}
\asqrt)^{2}\tilde{x}(t^\frac{1}{2}) - x - t^{-\frac{1}{2}} \asqrt
\tilde{\varphi}(t^{\frac{1}{2}}) \varphi) L(-1) + (t^{-\frac{1}{2}}
\asqrt \tilde{\varphi} (t^\frac{1}{2}) - \varphi) G(-\frac{1}{2})} \cdot  \\
e^{- \sum_{j \in \Z} \left( \Theta^{(1)}_j (t^{\frac{1}{2}}) L(j) +
\Theta^{(1)}_{j - \frac{1}{2}} (t^{\frac{1}{2}}) G(j - \frac{1}{2}) \right)}
\cdot   e^{- 2\Theta^{(1)}_0 L(0)} u, (x,\varphi) \Bigr) \\
(t^{-\frac{1}{2}} \asqrt)^{2L(0)} \cdot e^{- \sum_{j \in \Z} \left( A_j
L(j) + M_{j -\frac{1}{2}} G(j - \frac{1}{2})
\right)} \cdot (t^{-\frac{1}{2}} \asqrt)^{-2L(0)} \cdot v  
\end{multline*}
\begin{multline*}
= \; \; (t^{-\frac{1}{2}} \asqrt)^{-2L(0)} Y\biggl( e^{- \sum_{j \in \Z}
\left( \Theta^{(1)}_j (t^{\frac{1}{2}}) L(j) + \Theta^{(1)}_{j - \frac{1}{2}}
(t^{\frac{1}{2}}) G(j - \frac{1}{2}) \right)} \cdot e^{-2\Theta^{(1)}_0
(t^{\frac{1}{2}}) L(0)} u, \biggr. \\ 
\biggl. (t^{-1} \asqrt^2 \tilde{x} (t^\frac{1}{2}), t^{-\frac{1}{2}} 
\asqrt \tilde{\varphi}(t^\frac{1}{2}))\biggr)  (t^{-\frac{1}{2}}
\asqrt)^{2L(0)} \cdot e^{- \sum_{j \in \Z} \left( A_j L(j) + M_{j
-\frac{1}{2}} G(j -\frac{1}{2}) \right)} \cdot \\
(t^{-\frac{1}{2}} \asqrt)^{-2L(0)} \cdot v 
\end{multline*}
\begin{eqnarray*}
&=&  Y\biggl((t^{-\frac{1}{2}} \asqrt)^{-2L(0)} e^{- \sum_{j \in \Z}
\left( \Theta^{(1)}_j(t^{\frac{1}{2}}) L(j) + \Theta^{(1)}_{j - \frac{1}{2}}
(t^{\frac{1}{2}}) G(j - \frac{1}{2}) \right)} \cdot e^{-2\Theta^{(1)}_0
(t^{\frac{1}{2}}) L(0)} u, \biggr. \\ 
& & \biggl.\hspace{.4in} (\tilde{x}
(t^\frac{1}{2}),\tilde{\varphi}(t^\frac{1}{2}))\biggr) 
 \cdot e^{- \sum_{j \in \Z} \left( A_j L(j) + M_{j -\frac{1}{2}} G(j -
\frac{1}{2}) \right)} \cdot  (t^{-\frac{1}{2}} \asqrt)^{-2L(0)} \cdot v \\
&=& Y(\gamma_{\Theta(H,t^{1/2}, (x,\varphi))} (u),H_{t^{1/2}}(x,\varphi))
\gamma_{H,t^{1/2}}(v) ,
\end{eqnarray*}
for $u,v \in V$.  Thus we have 
\begin{equation}\label{conclusion}
\gamma_{H, t^{1/2}} (Y(u, (x,\varphi))v) =  Y(\gamma_{\Theta(H,t^{1/2},
(x,\varphi))} (u),H_{t^{1/2}}(x,\varphi)) \gamma_{H,t^{1/2}}(v).
\end{equation}
By the truncation condition and the fact that $V$ is a positive energy
module for the Neveu-Schwarz algebra, both sides of (\ref{conclusion})
are formal power series in $V[t^{-1/2},t^{1/2}]((x))[\varphi]$, and thus
well defined at $t^{1/2} =1$. This proves the first statement of
Proposition \ref{change prop} and equation (\ref{change formula}) for
general formal superconformal change of variables $H$ vanishing at zero. 
In addition, since $\gamma_{\Theta(H,t^{1/2}, (x,\varphi))}|_{t^{1/2} = 1}$
is well-defined, so is $\Theta_j^{(1)}(t^{1/2})|_{t^{1/2} = 1}$ and thus so
is $\gamma_{\Theta(H,t^{1/2}, (x,\varphi))}^{-1}|_{t^{1/2} = 1}$.  This,
the fact that $\gamma_{\Theta(H,(x,\varphi))} \circ \gamma_{\Theta(H,
(x,\varphi))}^{-1} (v) = v$ for $v \in V$, and the fact that $\gamma_H$ is
invertible, give us equation (\ref{change formula rewritten}).
\end{proof}

\subsection{Formal algebraic proof of Proposition \ref{change prop2} for
formal superconformal change of coordinates at infinity}\label{alg-proof2}

In this section, we let $H^{-1}(x,\varphi) \in x^{-1} \bigwedge_\infty
[[x^{-1}]] [\varphi]$ be the formal superconformal power series given by
(\ref{coordinate change22}), but this time we do not assume that
$H^{-1}$ converges in a neighborhood of infinity. 

As in Section \ref{alg-proof}, we will prove Proposition \ref{change
prop2} by a formal algebraic argument first employed in \cite{B iso-thm}
to construct an $N=1$ SG-VOSA {}from an $N=1$ NS-VOSA. 

We first prove the following bracket formula.

\begin{lem}\label{bracket-Lemma2}
Let $(V,Y(\cdot,(x,\varphi)), \mathbf{1}, \tau)$ be an $N=1$
NS-VOSA over a Grassmann algebra.  For $u \in V$, $t^{1/2}$ a formal even
variable,  and $(B_j,N_{j-1/2}) \in \bigwedge_\infty$, for $j \in \Z$, we
have the following bracket formula in $(\mathrm{End} \; V) [[t^{1/2}]]
[[x^{-1},x]] [\varphi]$
\begin{equation*}
\Bigl[ \sum_{j \in \Z} \Bigl( t^j B_j L(-j) + t^{j -\frac{1}{2}} N_{j
-\frac{1}{2}} G( - j + \frac{1}{2}) \Bigr) , Y(u,(x,\varphi)) \Bigr]
\hspace{.6in}
\end{equation*}
\begin{multline}\label{bracket2}
= Y \biggl( \sum_{m = -1}^{\infty} \sum_{j \in \Z} \binom{-j + 1}{m + 1}
x^{-j - m} \biggl( \Bigl( t^j B_j + 2\varphi t^{j - \frac{1}{2}} 
N_{j - \frac{1}{2}}\Bigr) L(m)  \\
+ \; \Bigl( t^{j - \frac{1}{2}}  N_{j - \frac{1}{2}} + \varphi
x^{-1} \frac{(-j-m)}{2} t^j B_j \Bigr) G(m + \frac{1}{2})
\biggr) u, (x,\varphi) \biggr) .
\end{multline} 
\end{lem}

\begin{proof}
{}From (\ref{bracket relation for a vosa}), we have
\begin{multline*} 
\left[ Y(\tau, (x_1,\varphi_1)), Y(u,(x,\varphi)) \right] \\
= \mbox{Res}_{x_0} x^{-1} \delta \left( \frac{x_1 - x_0 -
\varphi_1 \varphi}{x} \right) Y(Y(\tau,(x_0, \varphi_1 -
\varphi))u,(x, \varphi)) .
\end{multline*} 
Let 
\[l_t(x_1) = \sum_{j \in \Z} t^j B_j x_1^{-j + 1} \qquad
\mbox{and} \qquad g_t(x_1) = \sum_{j \in \Z} t^{j - \frac{1}{2}} N_{j -
\frac{1}{2}} x_1^{-j + 1} .\]
Then 
\begin{eqnarray*}
\varphi_1 \mbox{Res}_{x_1} l_t(x_1) Y(\omega, (x_1,\varphi_1)) &=&
\varphi_1 \sum_{j \in \Z} t^j B_j L(-j) \\ 
\varphi_1 \mbox{Res}_{x_1} g_t(x_1) Y(\tau, (x_1,\varphi_1)) &=&
\varphi_1 \sum_{j \in \Z} t^{j - \frac{1}{2}} N_{j - \frac{1}{2}} G(- j +
\frac{1}{2}) ,
\end{eqnarray*}
and by the $G(-1/2)$-derivative property
(\ref{G(-1/2)-derivative}), we have 
\begin{eqnarray*}
Y(\omega, (x_1,\varphi_1)) &=& Y( \frac{1}{2} G(-\frac{1}{2}) \tau,
(x_1, \varphi_1)) \\
&=& \frac{1}{2} \Bigl( \frac{\partial}{\partial \varphi_1} + \varphi_1
\frac{\partial}{\partial x_1} \Bigr) Y(\tau, (x_1,\varphi_1))  .
\end{eqnarray*}
Thus 
\begin{eqnarray*}
& & \hspace{-.5in} \varphi_1 \biggl[ \sum_{j \in \Z} \left( t^j B_j
L(-j) + t^{j - \frac{1}{2}} N_{j -\frac{1}{2}} G( - j +
\frac{1}{2}) \right) , Y(u,(x,\varphi)) \biggr] \\
&=& \varphi_1 \mbox{Res}_{x_1} \left[ \Bigl( \frac{1}{2}l_t(x_1)
\Bigl( \frac{\partial}{\partial \varphi_1} + \varphi_1
\frac{\partial}{\partial x_1} \Bigr) + g_t(x_1) \Bigr) Y(\tau,
(x_1,\varphi_1)), Y(u,(x,\varphi)) \right] \\
&=& \varphi_1 \mbox{Res}_{x_1} \left( \Bigl( \frac{1}{2}l_t(x_1) \Bigl(
\frac{\partial}{\partial \varphi_1} + \varphi_1
\frac{\partial}{\partial x_1} \Bigr) + g_t(x_1) \Bigr) \right. \\
& & \hspace{.7in} \left. \mbox{Res}_{x_0} x^{-1} \delta \Bigl(
\frac{x_1 - x_0 - \varphi_1 \varphi}{x} \Bigr) Y(Y(\tau,(x_0,
\varphi_1 - \varphi))u,(x, \varphi)) \right)\\
&=& Y\left(\varphi_1 \mbox{Res}_{x_1}\mbox{Res}_{x_0} \Bigl(
\frac{1}{2}l_t(x_1) \Bigl( \frac{\partial}{\partial \varphi_1} +
\varphi_1 \frac{\partial}{\partial x_1} \Bigr) + g_t(x_1) \Bigr)
\right. \\
& & \left. \hspace{1in} x^{-1} \delta \Bigl( \frac{x_1 - x_0 -
\varphi_1 \varphi}{x} \Bigr) Y(\tau,(x_0, \varphi_1 - \varphi))u,(x,
\varphi) \right) .  
\end{eqnarray*}
However using the $\delta$-function identity (\ref{delta 2 terms with
phis}), we have
\begin{multline*}
\varphi_1 \mbox{Res}_{x_1}\mbox{Res}_{x_0} \left(
\frac{1}{2}l_t(x_1) \Bigl( \frac{\partial}{\partial \varphi_1} +
\varphi_1 \frac{\partial}{\partial x_1} \Bigr) + g_t(x_1) \right) \\
x^{-1} \delta \Bigl( \frac{x_1 - x_0 - \varphi_1
\varphi}{x} \Bigr) Y(\tau,(x_0, \varphi_1 - \varphi)) 
\end{multline*}
\begin{multline*}
=  \varphi_1 \mbox{Res}_{x_1}\mbox{Res}_{x_0} \left( \frac{1}{2}
l_t(x_1) \Bigl( \frac{\partial}{\partial \varphi_1} + \varphi_1
\frac{\partial}{\partial x_1} \Bigr) + g_t(x_1) \right) \\
x_1^{-1} \delta \Bigl( \frac{x + x_0 + \varphi_1
\varphi}{x_1} \Bigr) Y(\tau,(x_0, \varphi_1 - \varphi)) 
\end{multline*}
\begin{multline*}
= \varphi_1 \mbox{Res}_{x_1}\mbox{Res}_{x_0} \left(
\frac{1}{2}l_t(x_1) \Bigl( \frac{\partial}{\partial \varphi_1} 
x_1^{-1} \delta \Bigl( \frac{x + x_0 + \varphi_1 \varphi}{x_1} \Bigr)
\Bigr) Y(\tau,(x_0,- \varphi)) \right.\\ 
 + \; \frac{1}{2}l_t(x_1) x_1^{-1} \delta \Bigl(
\frac{x + x_0}{x_1} \Bigr) \frac{\partial}{\partial \varphi_1}
Y(\tau, (x_0,\varphi_1)) \\
\left. + \; g_t(x_1) x_1^{-1} \delta \Bigl( \frac{x
+ x_0}{x_1} \Bigr) Y(\tau,(x_0,- \varphi)) \right) 
\end{multline*}
\begin{multline*}
= \varphi_1 \mbox{Res}_{x_1} \mbox{Res}_{x_0} \biggl( 
\biggl(\sum_{n \in \mathbb{Z}} n (x + x_0)^{n-1} \varphi x_1^{-n-1}
\biggr) \biggl( \frac{1}{2} \sum_{j \in \Z} t^j B_j x_1^{-j + 1}
\biggr) \biggr.\\
\biggl( \sum_{k \in \mathbb{Z}} G(k + \frac{1}{2}) x_0^{-k - 2} \biggr) 
+ \biggl(\sum_{n \in \mathbb{Z}}(x + x_0)^n x_1^{-n-1} \biggr) 
\biggl( \frac{1}{2} \sum_{j \in \Z} t^j B_j x_1^{-j + 1} \biggr)\\
\biggl( 2\sum_{k \in \mathbb{Z}} L(k) x_0^{-k - 2} \biggr)
+ \biggl(\sum_{n \in \mathbb{Z}}(x + x_0)^n x_1^{-n-1} \biggr)
\biggl( \sum_{j \in \Z} t^{j-\frac{1}{2}} N_{j - \frac{1}{2}}
x_1^{-j + 1} \biggr) \\
\biggl. \biggl( \sum_{k \in \mathbb{Z}} G(k + \frac{1}{2})
x_0^{-k - 2} - 2\varphi L(k) x_0^{-k - 2} \biggr) \biggr)
\end{multline*}
\begin{multline*}
= \varphi_1 \mbox{Res}_{x_0} \biggl( \sum_{l \in \mathbb{N}} \sum_{k \in
\mathbb{Z}} \sum_{j \in \Z} \biggl( (-j + 1) \binom{-j}{l} x^{-j -l} x_0^l
\frac{\varphi}{2} t^j B_j G(k + \frac{1}{2}) x_0^{-k - 2} \biggr. \biggr.
\\ 
+ \; \binom{-j + 1}{l} x^{-j-l + 1} x_0^l t^j B_j L(k) x_0^{-k - 2}
\\
\biggl. \biggr. + \; \binom{-j + 1}{l} x^{-j-l + 1} x_0^l
t^{j-\frac{1}{2}} N_{j - \frac{1}{2}} \Bigl(G(k + \frac{1}{2}) -
2\varphi L(k) \Bigr) x_0^{-k - 2} \biggr) \biggr) 
\end{multline*}
\begin{multline*}
= \varphi_1 \sum_{l \in \mathbb{N}} \sum_{j \in \Z} \left( \frac{(-j +
1)}{2} \binom{-j}{l} x^{-j -l}  \varphi t^j B_j G(l - \frac{1}{2}) \right.
\\ 
+ \; \binom{-j + 1}{l} x^{-j-l + 1} t^j B_j L(l - 1) \\
\left. + \; \binom{-j + 1}{l} x^{-j-l + 1} t^{j-\frac{1}{2}}
N_{j - \frac{1}{2}} (G(l - \frac{1}{2}) - 2\varphi L(l - 1)) \right) 
\end{multline*}
\begin{multline*}
= \varphi_1 \sum_{m = -1}^\infty \sum_{j \in \Z} \left( \binom{-j + 1}{m+1} 
\left( t^j B_j L(m) + t^{j-\frac{1}{2}} N_{j - \frac{1}{2}} 
(G(m + \frac{1}{2}) \right) x^{-j-m} \right. \\
\left. + \; \varphi \binom{-j + 1}{m+1}\left( \frac{(-j -m )}{2} x^{-1}
t^j B_j G(m + \frac{1}{2}) + 2 t^{j-\frac{1}{2}} N_{j - \frac{1}{2}} L(m)
\right) x^{-j-m} \right)
\end{multline*}
\begin{multline*}
= \varphi_1  \sum_{m = -1}^{\infty} \sum_{j \in \Z} \binom{-j + 1}{m + 1}
x^{-j - m} \biggl( \Bigl( t^j B_j + 2\varphi t^{j - \frac{1}{2}} 
N_{j - \frac{1}{2}} \Bigr) L(m)  \\
+ \; \Bigl( t^{j - \frac{1}{2}}  N_{j - \frac{1}{2}} + \varphi
x^{-1} \frac{(-j-m)}{2} t^j B_j \Bigr) G(m + \frac{1}{2})
\biggr) 
\end{multline*}
which finishes the proof of (\ref{bracket2}).
\end{proof}

We have the following immediate corollary to Lemma \ref{bracket-Lemma2}.

\begin{cor}\label{theta2-corollary}
\begin{multline*}
e^{\sum_{j \in \Z} \left( t^{j} B_j L(-j) +  t^{j - \frac{1}{2}}
N_{j - \frac{1}{2}} G(-j + \frac{1}{2}) \right)}  Y(u,(x,\varphi))
\cdot  \\
e^{-\sum_{j \in \Z} \left( t^{j} B_j L(-j) +  t^{j - \frac{1}{2}}
N_{j - \frac{1}{2}} G(-j + \frac{1}{2})  \right)}  
\end{multline*} 
\begin{multline*}
= \; Y \Biggl(\exp \Biggl(\! \sum_{m = -1}^{\infty} \sum_{j \in \Z}
\binom{-j + 1}{m + 1} x^{-j - m} \biggl( \Bigl( t^j B_j + 2\varphi
t^{j - \frac{1}{2}}  N_{j - \frac{1}{2}} \Bigr) L(m)  \\
+ \; \Bigl( t^{j - \frac{1}{2}}  N_{j - \frac{1}{2}} + \varphi
x^{-1} \frac{(-j-m)}{2} t^j B_j \Bigr) G(m + \frac{1}{2})
\biggr) \Biggr)  u, (x,\varphi) \Biggr)  .
\end{multline*} 
\end{cor}

We are now ready to prove Proposition \ref{change prop2} in the case of
a formal invertible superconformal change of variables $H \circ I$ taking
infinity to infinity.

\begin{proof}
Denote
\begin{eqnarray*}
H_{t^{1/2}} \circ I (x,\varphi) &=& (\tilde{x} (t^\frac{1}{2}),
\tilde{\varphi}(t^\frac{1}{2})) \\
&=& \exp \Biggl(- \! \! \sum_{j \in \mathbb{Z}_+} \Bigl(t^j B_j
L_{-j}(x,\varphi) + t^{j - \frac{1}{2}} N_{j - \frac{1}{2}}
G_{-j+\frac{1}{2}}(x,\varphi) \Bigr) \! \! \Biggr) \!  \cdot (x,\varphi) .
\end{eqnarray*}
By Corollary \ref{theta2-corollary},
Proposition \ref{second Theta identity in End}, and equation
(\ref{exponential L(-1) and  G(-1/2) property}), we have
\begin{eqnarray*}
& & \hspace{-.3in} Y(u, (x,\varphi)) \xi^*_{H \circ I, t^{1/2}} (v)   \\
&=& Y(u,(x,\varphi)) e^{- \sum_{j \in \Z} \left( t^j B_j L(-j) +
t^{j - \frac{1}{2}} N_{j - \frac{1}{2}} G(-j+\frac{1}{2}) \right)} \cdot
v  \\  
&=& e^{\! - \sum_{j \in \Z} \left( t^j B_j L(-j) + t^{j - \frac{1}{2}} N_{j
- \frac{1}{2}} G(-j+\frac{1}{2}) \right)} \! \! \cdot e^{ \sum_{j \in \Z}
\left( t^j B_j L(-j) + t^{j - \frac{1}{2}} N_{j - \frac{1}{2}}
G(-j+\frac{1}{2}) \right)} \\
& & \hspace{1.3in} Y(u,(x,\varphi)) e^{- \sum_{j \in \Z} \left( t^j B_j
L(-j) + t^{j - \frac{1}{2}} N_{j - \frac{1}{2}} G(-j+\frac{1}{2}) \right)}
\cdot v  
\end{eqnarray*}
\begin{multline*}  
= \; \; e^{- \sum_{j \in \Z} \left( t^j B_j L(-j) + t^{j - \frac{1}{2}} N_{j
- \frac{1}{2}} G(-j+\frac{1}{2}) \right)}  Y \Biggl(\exp \Biggl(\! 
\sum_{m = -1}^{\infty} \sum_{j \in \Z} \binom{-j + 1}{m + 1} x^{-j - m} \\
\biggl( \! \Bigl( t^j B_j + 2\varphi t^{j - \frac{1}{2}}  N_{j -
\frac{1}{2}} \Bigr) L(m)  + \Bigl( t^{j - \frac{1}{2}}  N_{j -
\frac{1}{2}} + \varphi x^{-1} \frac{(-j-m)}{2} t^j B_j \Bigr) G(m +
\frac{1}{2}) \biggr) \! \Biggr)  \\
u, (x,\varphi) \Biggr) v  
\end{multline*}
\begin{multline*}  
= \; \; e^{- \sum_{j \in \Z} \left( t^j B_j L(-j) + t^{j - \frac{1}{2}} N_{j
- \frac{1}{2}} G(-j+\frac{1}{2}) \right)} \cdot\\
Y \Biggl( e^{\left( (\tilde{x}(t^\frac{1}{2}) - x - \tilde{\varphi}
(t^\frac{1}{2}) \varphi) L(-1) + (\tilde{\varphi}(t^\frac{1}{2}) -
\varphi) G(-\frac{1}{2}) \right)} \cdot \\  
\exp \Biggl( \! - \! \sum_{j \in \Z} \Bigl( \Theta^{(2)}_j
(t^\frac{1}{2}) L(j) + \Theta^{(2)}_{j - \frac{1}{2}} (t^\frac{1}{2}) G(j
- \frac{1}{2}) \Bigr) \! \Biggr) \cdot  \\ 
\exp \left( - 2\Theta^{(2)}_0
(t^\frac{1}{2}) L(0) \right) \cdot u, (x,\varphi) \Biggr) v  
\end{multline*}
\begin{eqnarray*}  
&=& e^{- \sum_{j \in \Z} \left( t^j B_j L(-j) + t^{j - \frac{1}{2}} N_{j -
\frac{1}{2}} G(-j+\frac{1}{2}) \right)} \cdot   \\  
& & \hspace{.5in} Y \Biggl( \exp \Biggl( \! - \! \sum_{j \in \Z} \Bigl(
\Theta^{(2)}_j (t^\frac{1}{2}) L(j) + \Theta^{(2)}_{j - \frac{1}{2}}
(t^\frac{1}{2}) G(j - \frac{1}{2}) \Bigr) \! \Biggr) \cdot  \\ 
& & \hspace{1.9in} \exp \left( - 2\Theta^{(2)}_0 (t^\frac{1}{2}) L(0)
\right)
\cdot u, (\tilde{x}(t^\frac{1}{2}),\tilde{\varphi}(t^\frac{1}{2})) \Biggr)
v  \\   
&=& \xi^*_{H \circ I,t^{1/2}} (Y(\xi_{\Theta(H \circ I,t^{1/2},
(x,\varphi))}(u), H_{t^{1/2}} \circ I (x,\varphi)) v ) ,\\ 
\end{eqnarray*}
for $u,v \in V$.  Thus we have 
\begin{equation}\label{conclusion2}
Y(u, (x,\varphi)) \xi^*_{H \circ I, t^{1/2}} (v) =  \xi^*_{H \circ
I,t^{1/2}}  (Y(\xi_{\Theta(H \circ I,t^{1/2}, (x,\varphi))}(u), H_{t^{1/2}}
\circ I (x,\varphi)) v )
\end{equation}
as formal power series in $V[[t^{1/2}]][[x^{-1},x]][\varphi]$. Setting
$t^{1/2} = 1$ on the left-hand side of (\ref{conclusion2}), we see that
both sides of (\ref{conclusion2}) are well-defined elements of
$\bar{V}[[x^{-1},x]][\varphi]$.  This proves the first statement of
Proposition \ref{change prop2} and equation (\ref{change formula2}) for
general formal superconformal change of variables $H$ at infinity. 
In addition, since $\xi_{\Theta(H \circ I,t^{1/2}, (x,\varphi))}|_{t^{1/2}
= 1}$ is well defined, $\Theta_j^{(2)}(t^{1/2})|_{t^{1/2} = 1}$ must also be
well defined, and thus so is $\xi_{\Theta(H \circ I,t^{1/2},
(x,\varphi))}^{-1}|_{t^{1/2} = 1}$.  This, the fact that for $v \in
V$, we have $\xi_{\Theta(H \circ I,(x,\varphi))} \circ \xi_{\Theta(H \circ
I, (x,\varphi))}^{-1} (v) = v$ , and the fact that $\xi^*_{H \circ I}$ is
invertible, give us equation (\ref{change formula2 rewritten}).
\end{proof}

\section{Isomorphism classes of $N=1$ NS-VOSAs arising {}from
superconformal changes of variables}\label{iso-families-section}

For $H$ a formally superconformal of the form (\ref{coordinate change}),
we note that (\ref{change formula}) is equivalent to
\begin{equation}
\gamma_H(Y(u,(x,\varphi))v) = Y(\gamma_{\Theta(H, (x,\varphi))} \circ
\gamma_H^{-1} \circ \gamma_H (u),H(x,\varphi)) \gamma_H(v) .
\end{equation} 
Thus by Remark \ref{homo remark}, we have the following corollary to
Proposition \ref{change prop}.

\begin{cor}
Let $H$, $\gamma_H$ and $\gamma_{\Theta(H, (x,\varphi))}$ be as in Section
\ref{change-zero}.  Define
\begin{eqnarray}
V_H &=& \gamma_H(V) , \\
Y_H (u,(x,\varphi)) &=& Y(\gamma_{\Theta(H, (x,\varphi))} \circ
\gamma_H^{-1} (u),  H(x,\varphi)) ,\\ 
\mathbf{1}_H &=& \gamma_H(\mathbf{1}) \\
&=& \mathbf{1} , \nonumber\\
\tau_H &=& \gamma_H(\tau)\\
&=&  \asqrt^{-3} \tau \nonumber .
\end{eqnarray}
Then $(V_H,Y_H(\cdot, (x,\varphi)), \mathbf{1}_H,\tau_H)$ is an $N=1$
NS-VOSA and is isomorphic to $(V,Y,\mathbf{1},\tau)$.
\end{cor}

Similarly for $H^{-1}$ a formally superconformal power series of the
form (\ref{coordinate change22}), we  note that since $\xi^*_{H \circ I}$ is
invertible, (\ref{change formula2}) is equivalent to
\begin{eqnarray*}
\lefteqn{(\xi^*_{H \circ I})^{-1} (Y(u, (x,\varphi))  v)}\\
&=& Y(\xi_{\Theta(H \circ I, (x,\varphi))}(u), H \circ I (x,\varphi))
(\xi^*_{H \circ I})^{-1} (v) )\\ 
&=&  Y(\xi_{\Theta(H \circ I, (x,\varphi))} \circ \xi^*_{H \circ I} \circ
(\xi^*_{H \circ I})^{-1} (u), H \circ I (x,\varphi))
(\xi^*_{H \circ I})^{-1} (v) ) .
\end{eqnarray*} 
Thus by Remark \ref{homo remark}, we have the following corollary to
Proposition \ref{change prop2}.

\begin{cor}
Let $H \circ I$, $\xi_{H \circ I}^*$ and $\xi_{\Theta(H, (x,\varphi))}$ be
as in Section \ref{change-infinity}.  Define 
\begin{eqnarray}
V_{H \circ I} &=& \coprod_{n \in \frac{1}{2}\mathbb{Z}}
(\xi^*_{H \circ I})^{-1} (V_{(n)}) , \\  
Y_{H \circ I} (u,(x,\varphi)) &=& Y(\xi_{\Theta(H \circ I, (x,\varphi))}
\circ \xi_{H \circ I}^* (u), H \circ I (x,\varphi)),\\ 
\mathbf{1}_{H \circ I} &=& (\xi^*_{H \circ I})^{-1} (\mathbf{1}) \\
&=& \mathbf{1} + N_{\frac{3}{2}} \tau +  B_2 \omega + \mbox{terms of
higher weight in $V$} ,
\nonumber\\
\tau_{H \circ I} &=& (\xi^*_{H \circ I})^{-1}(\tau) \\
&=&  \tau +  \mbox{terms of higher weight in $V$}  \nonumber .
\end{eqnarray}
Then $(V_{H \circ I},Y_{H \circ I} (\cdot,(x,\varphi)),
\mathbf{1}_{H \circ I}, \tau_{H \circ I})$ is an $N=1$ NS-VOSA and is
isomorphic to $(V,Y,\mathbf{1},\tau)$.
\end{cor}

\begin{rema}\label{H-mistake}
{\em  In \cite{H book}, Huang defines families of isomorphic
VOAs derived {}from a change of variables $f$ where $f$ is an
invertible analytic function vanishing at infinity.  However, there is a
mistake in his definition of the operator $Y_{1/f} (u_{1/f},x)v_{1/f}$ on
p. 181.  Instead of defining
\begin{eqnarray}
Y_{1/f} (u_{1/f},x)v_{1/f} &=& Y(\psi_x^{1/f}(u_{1/f}), 1/f(x))v_{1/f}\\
&=& Y(\xi_x^{-1} \circ (\xi^*_{1/f})^{-1} (u_{1/f}), 1/f(x))v_{1/f}
\nonumber
\end{eqnarray}
one should define
\begin{eqnarray}
Y_{1/f} (u_{1/f},x)v_{1/f} &=& Y(\psi_{1/f(x)}^{1/f}(u_{1/f}),
1/f(x))v_{1/f}\\ &=& Y(\xi_{1/f(x)}^{-1} \circ (\xi^*_{1/f})^{-1} (u_{1/f}),
1/f(x))v_{1/f} . \nonumber
\end{eqnarray}
The superextension of this is equivalent to using the equation
(\ref{change formula2}) to note that 
\begin{eqnarray*}
\lefteqn{\xi_{H \circ I}^*(Y( u, (x,\varphi))  v)} \\ 
&=& Y( \xi_{\Theta(H \circ I, I^{-1} \circ H^{-1} (x,\varphi))}^{-1} (u)
,I^{-1} \circ H^{-1}(x,\varphi))\xi_{H \circ I}^*(v) \nonumber\\
&=& Y( \xi_{\Theta(H \circ I, I^{-1} \circ H^{-1} (x,\varphi))}^{-1} \circ
(\xi_{H \circ I}^*)^{-1} \circ  \xi_{H \circ I}^* (u) ,I^{-1}
\circ H^{-1}(x,\varphi))\xi_{H \circ I}^*(v) .\nonumber
\end{eqnarray*}
Thus defining 
\begin{eqnarray}
V_{I^{-1} \circ H^{-1}} &=& \coprod_{n \in \frac{1}{2}\mathbb{Z}}
\xi^*_{H \circ I} (V_{(n)}) , \\  
Y_{I^{-1} \circ H^{-1}} (u,(x,\varphi)) &=&  Y( \xi_{\Theta(H \circ I, 
I^{-1} \circ H^{-1} (x,\varphi))}^{-1} \circ
(\xi_{H \circ I}^*)^{-1}  (u) ,\\
& & \hspace{1.5in} I^{-1} \circ H^{-1}(x,\varphi)),\\ 
\mathbf{1}_{I^{-1} \circ H^{-1}} &=& \xi^*_{H \circ I} (\mathbf{1}) \\
\tau_{I^{-1} \circ H^{-1}} &=& \xi^*_{H \circ I} (\tau) \\
\end{eqnarray}
gives an $N=1$ NS-VOSA isomorphic to $(V,Y,\mathbf{1},\tau)$.  One recovers
the correct analogous formulas for the non-super case given in \cite{H
book} by letting $f$ be the body component of $H^{-1}$, setting all odd
components and soul portions of supernumbers equal to zero and restricting
$V$ to $V^0$. (Recall that Huang's operator $\xi^*_{1/f}$ is the
body portion of our operator $\xi^*_{H\circ I}$ if $f^{-1}$ is the body of
$H$.)}
\end{rema}

\section{Superconformal change of variables in an annulus}\label{annulus}

Let $\Sigma_B$ be a closed annulus in the complex plane with Jordan curves
$C_1$ and $C_2$ as its boundary.  Let $\Sigma$ be the superannulus in
$\bigwedge_\infty$ with body $\Sigma_B$ and such that $\Sigma$ is closed 
under the DeWitt topology, i.e., let
\[\Sigma = \Sigma_B \times \mbox{$(\bigwedge_\infty)_S$}. \] 
Let $H^{-1}(z,\varphi)$ be an invertible orientation-preserving
superconformal map defined on $\Sigma$. 

Let $\Sigma_2$ be the open subset of $\bigwedge_\infty$ in the DeWitt 
topology whose body is the exterior of $C_2$ in the complex plane union
$\{\infty\}$.  Assume that $C_1$ and $C_2$ were chosen such that $C_1
\subset (\Sigma_2)_B$.  Assume that the interior of $C_2$ and the interior
of the body of the subset $H^{-1} (C_2 \times (\bigwedge_\infty)_S)$ both
contain zero.  Let $\Sigma_1$ be the open subset of $\bigwedge_\infty$ in
the DeWitt topology whose body is the interior of the body of $H^{-1}(C_1
\times (\bigwedge_\infty)_S)$ in the complex plane.

Let 
\[M = \Sigma_1 \sqcup \Sigma_2 / \{p = q \; \mathrm{if} \; p\in \Sigma,
q \in H^{-1}(\Sigma), \; \mathrm{and} \; H^{-1}(p) = q \} .\]
Then $M$ is a super-Riemann surface with superconformal coordinate charts
given by the two open neighborhoods $\Sigma_1$ and $\Sigma_2$, and the
superconformal coordinate transition $H^{-1}$.  It is clear that
topologically, $M$ is a supersphere.  Thus by the Uniformization Theorem for
super-Riemann surfaces proved in \cite{CR}, $M$ must be superconformally
equivalent to the super-Riemann sphere $S\hat{\mathbb{C}}$.   Let $F$ be a
superconformal isomorphism {}from $M$ to $S\hat{\mathbb{C}}$.  Then $F$
restricted to $\Sigma_1$ and $\Sigma_2$, respectively, gives invertible
superconformal functions $F_1$ and $F_2$, respectively, such that for
$(w,\rho)$ in the interior of $\Sigma$
\[F_2 (w,\rho) = F_1 \circ H^{-1} (w, \rho), \]
or equivalently
\begin{equation} \label{compatibility}
H^{-1}(w,\rho) = F_1^{-1} \circ F_2 (w, \rho) . 
\end{equation}
We can assume that 
\begin{eqnarray}
F_1 (0,0) &=& (0,0) \label{annulus-condition1} \\
F_2 (\infty, 0) &=& (\infty,0) \label{annulus-condition2}\\
\lim_{w \rightarrow \infty} \frac{\partial}{\partial \rho} F_2 (w,\rho)
&=& 1 , \label{annulus-condition3}
\end{eqnarray}
since if equations (\ref{annulus-condition1})-(\ref{annulus-condition3}) do not
hold, we can compose $F$ with a global superconformal transformation $T$ {}from
$S\hat{\mathbb{C}}$ to $S\hat{\mathbb{C}}$, such that equations 
(\ref{annulus-condition1})-(\ref{annulus-condition3}) hold for $T \circ F$;
see \cite{B memoirs}.

Since $F_1$ is superconformal and invertible in a neighborhood of zero,
and vanishing at zero, $F_1 (x,\varphi)$ can be written in the form of
(\ref{coordinate change}), and since $F_2$ is superconformal satisfying 
equations (\ref{annulus-condition2}) and (\ref{annulus-condition3}),
$F_2$ can be written in the form of the righthand side of
(\ref{coordinate change3}) with $t^{1/2} = 1$.  Thus by Propositions
\ref{change prop} and \ref{change prop2}, we have
\begin{eqnarray*}
\lefteqn{(\xi_{F_2^{-1}}^*)^{-1} \circ \gamma_{F_1} (Y(u,(x,\varphi)) v)}
\\  
&=& (\xi_{F_2^{-1}}^*)^{-1} (Y(\gamma_{\Theta( F_1 , (x,\varphi))} 
(u), F_1 (x,\varphi)) \gamma_{F_1} (v)) \\
&=& Y(\xi_{\Theta(F_2^{-1}, F_1 (x,\varphi))} \circ
\gamma_{\Theta(F_1, (x,\varphi))}  (u), F_2^{-1} \circ F_1
(x,\varphi)) (\xi^*_{F_2^{-1}})^{-1} \circ \gamma_{F_1} (v) \\
&=& Y(\xi_{\Theta(F_2^{-1}, F_1(x,\varphi))} \circ \gamma_{\Theta(F_1,
(x,\varphi))}  (u), H (x,\varphi)) (\xi^*_{F_2^{-1}})^{-1} \circ
\gamma_{F_1} (v) \\      
&=& Y(\xi_{\Theta(F_2^{-1}, F_1 (x,\varphi))} \circ \gamma_{\Theta(F_1,
(x,\varphi))} \circ \gamma_{F_1}^{-1} \circ \xi_{F_2^{-1}}^*  \circ
(\xi_{F_2^{-1}}^*)^{-1} \circ \gamma_{F_1} (u), \\ 
& & \hspace{3in} H (x,\varphi)) (\xi^*_{F_2^{-1}})^{-1} \circ
\gamma_{F_1} (v) .\\ 
\end{eqnarray*}
Thus by Remark \ref{homo remark}, we have the following Corollary to 
Propositions \ref{change prop} and \ref{change prop2}.

\begin{cor}\label{annulus-corollary}
Let $H^{-1}$ be a superconformal map defined on the superannulus $\Sigma$
as  given above, and let $F_1$ and $F_2$ be the unique superconformal
coordinates defined in neighborhoods of $0$ and $\infty$, respectively, and
satisfying  (\ref{compatibility})-(\ref{annulus-condition3}).  Let $(V,
Y(\cdot,  (x,\varphi)), \mathbf{1}, \tau)$ be an $N=1$ NS-VOSA.  We have the
following change of variables formulas in $\bar{V} [[x^{-1}, x]] [\varphi]$
\begin{multline}\label{change formula annulus}
(\xi_{F_2^{-1}}^*)^{-1} \circ \gamma_{F_1} (Y(u,(x,\varphi)) v)\\
= Y(\xi_{\Theta(F_2^{-1}, F_1 (x,\varphi))} \circ \gamma_{\Theta(F_1,
(x,\varphi))}  (u), H (x,\varphi)) (\xi^*_{F_2^{-1}})^{-1} \circ
\gamma_{F_1} (v) , 
\end{multline}
i.e.,
\begin{multline}\label{change formula annulus rewritten}
Y(u,  H (x,\varphi)) v \\
 = (\xi_{F_2^{-1}}^*)^{-1} \circ \gamma_{F_1} (Y(\gamma_{\Theta(F_1,
(x,\varphi))}^{-1} \circ \xi_{\Theta(F_2^{-1}, F_1 (x,\varphi))}^{-1}
(u),(x,\varphi)) \gamma_{F_1}^{-1} \circ \xi^*_{F_2^{-1}} (v)), 
\end{multline}
for $u,v \in V$, and for $(z,\theta)$ in the domain of convergence of $H$,
both sides of (\ref{change formula annulus}) and of (\ref{change formula
annulus rewritten}) exist for $(x,\varphi) = (z,\theta)$ and are equal.
Furthermore, define
\begin{eqnarray}
V_H &=& (\xi_{F_2^{-1}}^*)^{-1} \circ  \gamma_{F_1} (V) , \\
Y_H (u,(x,\varphi)) &=& Y(\xi_{\Theta( F_2^{-1}, F_1 (x,\varphi))} \circ
\gamma_{\Theta(F_1, (x,\varphi))} \circ  \gamma_{F_1}^{-1} \circ 
 \xi_{F_2^{-1}}^* (u),\\ 
& & \hspace{2.5in} H (x,\varphi)) ,\\ 
\mathbf{1}_H &=& (\xi_{F_2^{-1}}^*)^{-1} \circ  \gamma_{F_1} (\mathbf{1}),
\\
\tau_H &=& (\xi_{F_2^{-1}}^*)^{-1} \circ  \gamma_{F_1} (\tau).\\
\end{eqnarray}
Then $(V_H,Y_H(\cdot, (x,\varphi)), \mathbf{1}_H,\tau_H)$ is an
$N=1$ NS-VOSA and is isomorphic to $(V,Y,\mathbf{1},\tau)$.
\end{cor}

\begin{rema}\label{non-zero-remark2} 
{\em We have assumed that the annulus $\Sigma$ and its image under $H^{-1}$
circumscribe zero.  However, if it does not, we can use the superconformal
shift change of variables formulas given by (\ref{L(-1) and G(-1/2) exp 1})
which, when combined with the change of variables given in Corollary
\ref{annulus-corollary} above, provide the appropriate change of variables
formula. }
\end{rema}

\begin{rema}\label{H-mistake2}
{\em  The formula for $Y_{1/f}$ on p. 183 in \cite{H book} is
incorrect.  The formula should read
\begin{eqnarray*}
\lefteqn{Y_{1/f}(u_{1/f},x)v_{1/f}}  \\
&=& Y(\xi^{-1}_{(F^{(2)})^{-1} \circ F^{(1)}(x)}
\circ \psi^{F^{(1)}}_x \circ (\xi_{(F^{(2)})^{-1}}^*)^{-1} (u_{1/f}),
(F^{(2)})^{-1} \circ F^{(1)} (x)) v_{1/f}  \\
&=& Y(\xi^{-1}_{f(x)} \circ \psi^{F^{(1)}}_x \circ
(\xi_{(F^{(2)})^{-1}}^*)^{-1} (u_{1/f}), f (x)) v_{1/f} . 
\end{eqnarray*}
In addition, the isomorphism in Theorem 7.4.8 in \cite{H book} should be
$\xi^*_{(F^{(2)})^{-1}} \circ \varphi_{F^{(1)}}$. }
\end{rema}

\end{document}